\documentclass[11pt]{amsart}

\usepackage{latexsym}
\usepackage{amssymb}
\usepackage{amscd}
\usepackage[all]{xy}
\usepackage[mathscr]{euscript}
\usepackage{dsfont}
\usepackage{hyperref}
\usepackage{graphicx}
\usepackage{amsfonts}
\usepackage{amsmath}
\usepackage{amsthm}
\usepackage{aliascnt}

\normalsize

\RequirePackage{xspace}

\newcommand{\thought}[1]{}
\renewcommand{\thought}[1]{ \textbf{[#1]}}

\usepackage{enumerate}        
\newenvironment{roenumerate}{\begin{enumerate}[\upshape (i)]}{\end{enumerate}}

\newcommand\nc {\newcommand}
\newcommand\rnc{\renewcommand}
\renewcommand{\leq}{\leqslant}
\renewcommand{\geq}{\geqslant}

\newcount\blopone
\newcount\xone
\newcount\xtwo
\newcount\ytwo

\newtheorem{theorem}{Theorem}[section]
\newtheorem{prop}[theorem]{Proposition}
\newtheorem{com}[theorem]{Comment}
\newtheorem{redu}[theorem]{Reduction}
\newtheorem{refinement}[theorem]{Refinement}
\newtheorem{summary}[theorem]{Summary}
\newtheorem{importnota}[theorem]{Important Notation}
\newtheorem{prblm}[theorem]{Problem}
\newtheorem{notation}[theorem]{Notation}
\newtheorem{defin}[theorem]{Definition}
\newtheorem{caution}[theorem]{Caution}
\newtheorem{remark}[theorem]{Remark}
\newtheorem{reminder}[theorem]{Reminder}
\newtheorem{illustration}[theorem]{Illustration}
\newtheorem{observation}[theorem]{Observation}
\newtheorem{lemma}[theorem]{Lemma}
\newtheorem{construction}[theorem]{Construction}
\newtheorem{corollary}[theorem]{Corollary}
\newtheorem{example}[theorem]{Example}
\newtheorem{conclusion}[theorem]{Conclusion}
\newtheorem{triviality}[theorem]{Triviality}
\newtheorem{proto}[theorem]{Prototype Quasifibration}
\newtheorem{cauex}[theorem]{Cautionary Example}
\newtheorem{hypo}[theorem]{Hypothesis}
\newtheorem{subth}{ }[theorem]
\newtheorem{case}{Case}[theorem]
\newtheorem{ssubth}{ }[subth]
\newtheorem{facts}[theorem]{Facts}

\nc\tri[1]{\begin{triviality}
\label{#1}}
\nc\fac[1]{\begin{facts}
\label{#1}
\begin{em}}
\nc\cas[1]{\begin{case}
\label{#1}
\begin{em}}
\nc\rfn[1]{\begin{refinement}
\label{#1}}
\nc\prt[1]{\begin{proto}
\label{#1}}
\nc\lem[1]{\begin{lemma}
\label{#1}}
\nc\pro[1]{\begin{prop}
\label{#1}}
\nc\thm[1]{\begin{theorem}
\label{#1}}
\nc\cor[1]{\begin{corollary}
\label{#1}}
\nc\dfn[1]{\begin{defin}
\label{#1}}
\nc\sthm[1]{\begin{subth}
\label{#1}}
\nc\exm[1]{\begin{example}
\label{#1}
\begin{em}}
\nc\obs[1]{\begin{observation}
\label{#1}
\begin{em}}
\nc\plm[1]{\begin{prblm}
\label{#1}
\begin{em}}
\nc\rmk[1]{\begin{remark}
\label{#1}
\begin{em}}
\nc\rmd[1]{\begin{reminder}
\label{#1}
\begin{em}}
\nc\ntn[1]{\begin{notation}
\label{#1}
\begin{em}}
\nc\smr[1]{\begin{summary}
\label{#1}
\begin{em}}
\nc\xpl[1]{\begin{explanation}
\label{#1}
\begin{em}}
\nc\cau[1]{\begin{caution}
\label{#1}
\begin{em}}
\nc\hyp[1]{\begin{hypo}
\label{#1}}
\nc\imn[1]{\begin{importnota}
\label{#1}
\begin{em}}
\nc\rdn[1]{\begin{redu}
\label{#1}
\begin{em}}
\nc\cax[1]{\begin{cauex}
\label{#1}
\begin{em}}
\nc\cmt[1]{\begin{com}
\label{#1}
\begin{em}}
\nc\con[1]{\begin{construction}
\label{#1}
\begin{em}}
\nc\ill[1]{\begin{illustration}
\label{#1}
\begin{em}}
\nc\ssthm[1]{\begin{ssubth}
\label{#1}
\begin{em}}
\nc\cnc[1]{\begin{conclusion}
\label{#1}
\begin{em}}

\nc\elem{\end{lemma}}
\nc\erdn{\end{em}\end{redu}}
\nc\erfn{\end{refinement}}
\nc\eprt{\end{proto}}
\nc\ethm{\end{theorem}}
\nc\ecor{\end{corollary}}
\nc\edfn{\end{defin}}
\nc\esthm{\end{subth}}
\nc\epro{\end{prop}}
\nc\etri{\end{triviality}}
\nc\eexm{\end{em}
\end{example}}
\nc\eobs{\end{em}
\end{observation}}
\nc\ecmt{\end{em}
\end{com}}
\nc\efac{\end{em}
\end{facts}}
\nc\ermk{\end{em}
\end{remark}}
\nc\ermd{\end{em}
\end{reminder}}
\nc\eill{\end{em}
\end{illustration}}
\nc\eplm{\end{em}
\end{prblm}}
\nc\ecas{\end{em}
\end{case}}
\nc\ecau{\end{em}
\end{caution}}
\nc\ecax{\end{em}
\end{cauex}}
\nc\eimn{\end{em}
\end{importnota}}
\nc\entn{\end{em}
\end{notation}}
\nc\econ{\end{em}
\end{construction}}
\nc\esmr{\end{em}
\end{summary}}
\nc\expl{\end{em}
\end{explanation}}
\nc\ehyp{
\end{hypo}}
\nc\ecnc{\end{em}
\end{conclusion}}
\nc\essthm{\end{em}
\end{ssubth}}

\newtheorem{explanation}[theorem]{Explanation}

\nc\sst{\scriptstyle}
\newcommand{\comment}[1]{}
\newcommand{\ri}{\longrightarrow}
\newcommand{\sr}{\rightarrow}

\newcommand{\zz}{{\mathbb Z}}

\newcommand{\D}{{\mathbf D}}

\newcommand{\C}{{\mathbb C}}

\nc\op{^{\hbox{\rm\tiny op}}}
\nc\mth{^{\hbox{\rm\tiny th}}}
\nc\wi{\wt{\text{\it\i}}}

\nc\script{\mathscr}
\nc\z{\zeta}
\nc\bc{{\mathbb{BC}}}
\nc\ct{{\script T}}
\nc\cf{{\script F}}
\nc\cg{{\script G}}
\nc\ch{{\script H}}
\nc\ck{{\script K}}
\nc\cl{{\script L}}
\nc\cu{{\script U}}
\nc\cv{{\script V}}
\nc\ce{{\script E}}
\nc\cs{{\script S}}
\nc\car{{\script R}}
\nc\cd{{\script D}}
\nc\cc{{\script C}}
\nc\ca{{\script A}}
\nc\ci{{\script I}}
\nc\co{{\script O}}
\nc\cx{{\script X}}
\nc\cy{{\script Y}}
\nc\cz{{\script Z}}
\nc\bd{\begin{description}}
\nc\ed{\end{description}}
\nc\ctob{{\script C}at\big(\ci^{op},\ca\big)}
\nc\clim{{\ds\mathop{\rm lim}_{\ds\longleftarrow}}\,}
\nc\climi{\clim_{\!i}\,}
\nc\climn{\clim^{\!n}\,}
\nc\colim{{\ds\mathop{\rm colim}_{\ds\la}}}
\nc\colimj{{\ds\mathop{\rm colim}_{\ds\la}}{}_{j\,}}
\nc\oa{\overline{\ca}}
\nc\s{\sigma}
\nc\ta{\tau}
\nc\os{\overline\sigma}
\nc\ot{\overline\tau}
\nc\T{\Sigma}
\nc\Tm{\Sigma^{-1}}
\nc\de[1]{{\mathop{\rm deg(#1)}}}
\nc\Ad[1]{\mathop{\rm Ad}(#1)}
\nc\ad[1]{\mathop{\rm ad}(#1)}
\nc\kth{{\it K}--theory}
\nc\loc[1]{{\text{\rm Loc(#1)}}}
\nc\coloc[1]{{\text{\rm Coloc}(#1)}}

\def\der #1 {D\left(#1\right)}
\nc\prf{\begin{proof}}
\nc\eprf{\end{proof}}
\nc\ds{\displaystyle}
\nc\Tor{\text{\rm Tor}}

\nc\cb{{\script B}}
\nc\ab{{\script A}b}

\nc\be{\begin{roenumerate}}
\nc\ee{\end{roenumerate}}

\nc\cat[1]{{\script C}at\Big({\big\{#1\big\}}\op\,\,,\,\,\ab\Big)}
\nc\csab{{\script C}at\big(\cs^{op},\ab\big)}
\nc\ctab{{\script C}at\Big({\{\ct^\alpha\}}^{op},\ab\Big)}
\nc\csex{{\script E}x\big(\cs^{op},\ab\big)}
\nc\ctex{{\script E}x\Big({\{\ct^\alpha\}}^{op},\ab\Big)}
\nc\sub{\qquad\subset\qquad}
\nc\ctr[1]{{\left.\ct\left(-,#1\right)\right|}_{\cs}}
\nc\ctrf[2]{{\left.\ct\left(#1,#2\right)\right|}_{\cs}}
\nc\Ctr[1]{{\left.\ct\left(-,#1\right)\right|}_{\ct^\alpha}}
\nc\Ctrf[2]{{\left.\ct\left(#1,#2\right)\right|}_{\ct^\alpha}}

\nc\la{\longrightarrow}
\nc\nin{\noindent}
\nc\cad[1]{\text{card}(#1)}
\nc\eq{\quad=\quad}
\nc\BA{\begin{array}{c}}
\nc\EA{\end{array}}
\nc\barr{
\[
\begin{array}{cccccccccccccccc}
}
\nc\earr{
\end{array}
\]
}
\nc\as[1]{{\langle S\rangle}^{#1}}
\nc\sh{\text{\it shift}}

\nc\yy[1]{{\left.\ct\left(-,#1\right)\right|}_{\ct^c}}
\nc\vrep[2]{{\left.\ct\left(#1,#2\right)\right|}_{\ct^\alpha}}
\nc\da{\downarrow}
\nc\Hom{{\mathop{\rm Hom}}}
\nc\HHom{{\script H}{\mathop{\rm om}}}
\nc\End{{\mathop{\rm End}}}
\nc\Ext{{\mathop{\rm Ext}}}
\nc\mr\Modtc
\nc\PExt{{\mathop{\rm PExt}}}
\nc\stm{\text{\rm stmod}(kG)}
\nc\stM{\text{\rm StMod}(kG)}
\nc\e{\varepsilon}
\nc\p{\varphi}

\nc\rs{\s^{-1}A}
\nc\br{{\{\s^{-1}A\}}}
\nc\ra\ri
\nc\y[1]{\mathbf{y}#1}
\nc\x[1]{\mathbf{z}#1}
\nc\mmod[1]{#1\text{--\rm mod}}
\nc\Mod[1]{#1\text{--\rm Mod}}
\nc\Md {\ensuremath{\mathop{\textup{Mod}}}}
\rnc\mod[1]{\ensuremath{\mathop{#1\textup{--mod}}}\xspace}
\nc\Modtc{\Mod{\ct^c}}
\nc\pgldim[1]{\mathop{\rm pgldim}\,#1}
\nc\tf{{\rm [TR5]}}
\nc\tfs{{\rm [TR5$^*$]}}
\nc\Fun{\text{\rm Funct}(F\op,\ab)}
\nc\sym{\text{\rm Sym}}
\nc\sgn{\text{\rm sgn}}
\nc\Pro{\text{\rm Prod}^{}_\alpha(F\op,\ab)}
\nc\Yt[1]{{\left.\Hom_\ct^{}\left(-,#1\right)\right|}_F^{}}
\nc\dl{\delta}
\nc\Proj[1]{#1\text{--\rm Proj}}
\nc\proj[1]{#1\text{--\rm proj}}
\nc\Flat[1]{#1\text{--\rm Flat}}
\nc\Inj[1]{#1\text{--\rm Inj}}
\nc\Ima{\mathrm{Im}}
\nc\Ker{\mathrm{Ker}}
\nc\ov{\overline}
\nc\wt{\widetilde}
\nc\wh{\widehat}
\nc\ph{\varphi}
\nc\tstr{{\it t}--structure}
\nc\spec[1]{{\text{\rm Spec}(#1)}}

\newcommand{\n}{\mathfrak{n}}
\nc\EProd{\text{\rm EProd}}
\nc\ECoprod{\text{\rm ECoprod}}
\nc\Prod{\text{\rm Prod}}
\nc\ldimp{\text{\rm LDim}^{\prod}}
\nc\ldimc{\text{\rm LDim}^{\coprod}}
\nc\gen[2]{{\langle#1\rangle}^{}_{#2}}
\nc\genu[3]{{\langle#1\rangle}^{[#3]}_{#2}}
\nc\ogen[1]{\ov{\langle#1\rangle}}
\nc\ogenun[2]{\ov{\langle#1\rangle}_{#2}^{}}
\nc\ogenu[3]{\ov{\langle#1\rangle}^{[#3]}_{#2}}
\nc\ogenul[3]{\ov{\langle#1\rangle}^{(-\infty,#3]}_{#2}}
\nc\ogenuf[3]{\ov{\langle#1\rangle}^{[#3,\infty)}_{#2}}
\nc\genuf[3]{{\langle#1\rangle}^{[#3,\infty)}_{#2}}
\nc\genul[3]{{\langle#1\rangle}^{(-\infty,#3]}_{#2}}
\nc\dperf[1]{\D^{\mathrm{perf}}(#1)}
\nc\dcoh{\mathbf{D}^b_{\mathrm{coh}}}
\newcommand{\Dqc}{{\mathbf D_{\text{\bf qc}}}}
\newcommand{\Dqcmi}{{\mathbf D_{\text{\bf qc}}^-}}
\newcommand{\Dqcpl}{{\mathbf D_{\text{\bf qc}}^+}}
\newcommand{\Dqcb}{{\mathbf D_{\text{\bf qc}}^b}}
\newcommand{\Dqcp}{{\mathbf D_{\text{\bf qc}}^p}}
\newcommand{\Dqcpb}{{\mathbf D_{\text{\bf qc}}^{p,b}}}
\nc\dmcoh{\mathbf{D}^-_{\mathrm{coh}}}
\nc\dscoh{\mathbf{D}^{}_{\mathrm{coh}}}
\nc\RHHom{{\script{RH}}{\mathrm{om}}}
\nc\Coprod{\mathrm{Coprod}}
\nc\COprod{\mathrm{coprod}}
\nc\add{\mathrm{add}}
\nc\Add{\mathrm{Add}}
\nc\Smr{\mathrm{smd}}
\nc\id{\mathrm{id}}
\nc\LL{\mathbf{L}}
\nc\R{\mathbf{R}}

\nc\hoco{
\begin{picture}(40,10)
\put(20,0){\makebox(0,0)[b]{\text{\rm Hocolim}}}
\put(5,-2){\vector(1,0){30}}
\end{picture}\,\,}

\nc\holim{
\begin{picture}(40,10)
\put(20,0){\makebox(0,0)[b]{\text{\rm Holim}}}
\put(35,-2){\vector(-1,0){30}}
\end{picture}}

\begin{document}

\author{Amnon Neeman}\thanks{The research was mostly supported 
  by the Australian Research Council
  (grants  number DP150102313 and number DP200102537), but partly also
  by the Deutsche
  Forschungsgemeinschaft (SFB-TRR 358/1 2023 - 491392403).}
\address{Dipartimento di Matematica ``F.\ Enriques''\\
        Universit{\`a} degli Studi di Milano\\
        Via Cesare Saldini 50\\
	20133 Milano\\
        ITALY}
\email{amnon.neeman@unimi.it}

\title[Two Brown representability theorems]{Triangulated categories with a single
  compact generator, and two Brown representability theorems}

\begin{abstract}
  We generalize theorems of Bondal and Van den Bergh and of Rouquier. A corollary
  of our main results says the following.
  
  Let $X$ be a scheme proper
  over a an excellent, finite-dimensional noetherian ring $R$. Then the Yoneda pairing, taking an object
  $A$ in the category $\dperf X$ and an object
  $B$ in the category
  $\dcoh(X)$, to the finite $R$--module $\Hom(A,B)$,
  gives an equivalence
  of $\dcoh(X)$ with the category of finite
  $R$--linear homological functors
  $H:\dperf X\op\la\mod R$, and an equivalence of $\dperf X\op$ with the category of finite homological functors $H:\dcoh(X)\la\mod R$. Recall: a homological functor $H$ is \emph{finite} if
  $\oplus_{i=-\infty}^\infty H^i(C)$ is a
  finite $R$--module for every $C\in\dperf X$.

  Bondal and Van den Bergh proved the special case, of the assertion about $\dcoh(X)$ identifying with the finite homological functors on $\dperf X\op$,  as long as $R$ is a field
  and $X$ is projective over $R$. And the assertion about $\dperf X\op$, identifying with the finite homological functors on $\dcoh(X)$, again under the assumption that $X$ is projective over a field $R$, is due to Rouquier.

  But our theorems are far more general yet. They aren't only about schemes, they work in the abstract generality
  of triangulated categories with coproducts and a single
  compact generator, satisfying a certain approximability
  property. At the moment I only know how to prove this
  approximability for the categories $\Dqc(X)$ with $X$
  a quasicompact,
  separated scheme, for the
   homotopy category of spectra, for the category $\D(R)$ where $R$
  is a (possibly noncommutative) negatively graded dg algebra,
and for certain recollements of the above.

  The work was inspired by Jack Hall's elegant new proof of a vast
  generalization of GAGA, a proof based on representability theorems of
  the type above. The generality of Hall's result made me wonder how far
  the known representability theorems could be improved.
\end{abstract}



\maketitle

\tableofcontents

\setcounter{section}{-1}

\section{Introduction}
\label{S0}

This paper begins with the
observation that triangulated categories with coproducts, and
with a single compact
generator, have a preferred class of {\it t}--structures.
This allows us to define thick
subcategories $\ct^-$, $\ct^+$ and $\ct^b$.
A slightly subtler definition is that of $\ct^-_c$. The full
subcategory $\ct^-_c\subset\ct$
makes sense unconditionally, and it is thick
as long as there exists a compact generator $G$
and an integer $A>0$ so that $\Hom(G,\T^iG)=0$ for all $i\geq A$.
We also define a subcategory $\ct^b_c=\ct^-_c\cap\ct^b$.

In the special case where $\ct=\Dqc(X)$, with $X$ a
quasicompact, quasiseparated scheme, the preferred
class of {\it t}--structures contains the standard
\tstr, the subcategories $\ct^-$, $\ct^+$ and $\ct^b$ are
nothing other than the classical
$\Dqcmi(X)$, $\Dqcpl(X)$ and $\Dqcb(X)$,
and if $X$ is noetherian the subcategories
$\ct_c^b\subset\ct_c^-$ can be proved
to be $\dcoh(X)\subset\dmcoh(X)$. What we have learned so far is
that these standard categories 
have an intrinsic
description. There is a method
to construct them out of $\ct$ in purely triangulated-category terms.

Still in the world of triangulated categories with coproducts and
a single compact generator: the category $\ct$ may be
\emph{weakly approximable} or even \emph{approximable.}
We will define these concepts later in the introduction, and study their properties
in the body of the paper.
For now we note that the category $\Dqc(X)$ is
weakly approximable
as long as $X$ is a quasicompact, quasiseparated
scheme, and approximable if $X$ is separated.
The homotopy category of spectra is also approximable.

To show that this abstraction can be useful we will prove representability
theorems. To state them we begin with 

\dfn{D37.105}
Let $R$ be a commutative ring,
let $\ct$ be an $R$--linear triangulated category
and let $\cb\subset\ct$ be a full, replete subcategory with $\T\cb=\cb$.
A \emph{$\cb$--homological functor} is an
$R$--linear functor
$H:\cb\la\Mod R$
which takes triangles to long exact sequences.
Dually, a \emph{$\cb$--cohomological functor} is an $R$--linear functor
$\wt H:\cb\op\la\Mod R$
which takes triangles to long exact sequences. This
means that, if we have a triangle $x\la y\la z\la\T x$ with
all three of $x,y,z$ belonging to $\cb$, then $H$
(respectively $\wt H$) takes it
to a long  exact sequence in $\Mod R$.

Suppose the ring $R$ is noetherian,
and let $G\in\cb\subset\ct$ be an object. The $\cb$--(co)homological
functor $H$ is called
\emph{$G$--locally finite} if
\be
\item
  $H(\T^iG)$ 
  is a finite $R$--module for all $i\in\zz$.
\item
  $H\big(\T^iG\big)=0$ for $i\ll 0$.
\setcounter{enumiv}{\value{enumi}}
\ee
The $\cb$--(co)homological functor $H$
is \emph{$G$--finite} if, in addition to
the above, we have
\be
\setcounter{enumi}{\value{enumiv}}
\item
  $H\big(\T^iG\big)=0$ 
  for $i\gg 0$.
\setcounter{enumiv}{\value{enumi}}
\ee
\edfn

\rmk{R-2.3}
Let $\ct$  be an $R$--linear triangulated category,
$\cb$ a full, replete subcategory with $\T\cb=\cb$, and $H$ a
$\cb$--(co)homological functor. If $H$ is $G$--locally finite
(respectively finite)
for every $G\in\cb$ we leave out the $G$, and just say that
$H$ is locally finite (respectively finite).

Note that if $H$ is $G$--locally finite (respectively finite) then
it is also $G'$--locally finite for any $G'$ obtainable from
$G$ by forming in $\cb$  finite direct sums, direct summands,
suspensions or triangles.
Thus local finiteness (respectively finiteness) can be checked
on any classical generator.
\ermk

\rmk{R29.102345}
The careful reader will observe that 
Definition~\ref{D37.105} is \emph{not} self-dual.
If we replace $\ct$ by $\ct\op$, and take
$\cb\subset\ct\op$
to be a full, replete subcategory with $\T\cb=\cb$,
then for a (co)-homological functor $H$ on $\cb$
to be $G$--finite is of course
equivalent to the (co)homological
functor $H$ being $G$--finite on
$\cb\op\subset\ct$.
But local finiteness changes, because on $\ct\op$ the
shift functor $\T$ is replaced by $\Tm$.

There exists a way to remedy this, by viewing
local finiteness as continuity with respect to certain
metrics. This will be discussed extensively in
subsequent articles.
\ermk

Our main theorem says that

\thm{T1.-1}
Let $R$ be a commutative, noetherian ring, and
$\ct$ an $R$--linear triangulated category with coproducts.
Assume $\ct$
is approximable, and suppose further that
there exists in $\ct$ a compact generator $G$ such that
$\Hom\big(G,G[n]\big)$ is a finite $R$--module for all $n\in\zz$.
Consider the two functors
\[
\cy:\ct^-_c\la\Hom_R^{}\big([\ct^c]\op\,,\,\Mod R\big),\qquad
\wt\cy:\big[\ct^-_c\big]\op\la\Hom_R^{}\big(\ct^b_c\,,\,\Mod R\big)
\]
defined by the formulas $\cy(B)=\Hom(-,B)$ and
$\wt\cy(A)=\Hom(A,-)$. Note that, in these formulas,
we permit all $A,B\in\ct^-_c$. But the $(-)$ in the formula
$\cy(B)=\Hom(-,B)$ is assumed to belong to $\ct^c$,
whereas the $(-)$ in the formula $\wt\cy(A)=\Hom(A,-)$ must
lie in $\ct^b_c$.
Now consider the following composites
\[\xymatrix@C+20pt@R-20pt{
\ct^b_c \ar@{^{(}->}[r]^i & \ct^-_c
\ar[r]^-{\cy} &
\Hom_R^{}\big([\ct^c]\op\,,\,\Mod R\big) \\
\big[\ct^c\big]\op \ar@{^{(}->}[r]^{\wi} & \big[\ct^-_c\big]\op
\ar[r]^-{\wt\cy} &
\Hom_R^{}\big(\ct^b_c\,,\,\Mod R\big)
}\]
We assert:
\be
\item
  The functor $\cy$ is full, and the essential image consists
  of the
  locally finite cohomological functors. The composite $\cy\circ i$
  is fully faithful, and the essential image consists of the
  finite cohomological functors.
\item
  With the notation as in Reminder~\ref{R1.-9}(xii),
  assume further that there exists an integer $N>0$ and an object $G'\in\ct^b_c$
  with $\ct=\ogen {G'}_N^{(-\infty,\infty)}$.
  Then
  the functor $\wt\cy$ is full, and the essential image consists
  of the
  locally finite homological functors. The composite $\wt\cy\circ \wi$
  is fully faithful, and the essential image consists of the
  finite homological functors.
\ee
\ethm

\rmk{R0.987656789}
The original versions of the result were posted
on the archive as two articles, with the first proving
Theorem~\ref{T1.-1}(i)
and the second devoted to 
Theorem~\ref{T1.-1}(ii). The reason they could not
be merged, at the time, was that the proof of
Theorem~\ref{T1.-1}(ii) depended on very useful
little results in joint work with Jesse Burke and
with Bregje Pauwels, see
\cite[Lemma~3.6 and 3.9(iv)]{Burke-Neeman-Pauwels18}.
And the joint article with Burke and Pauwels
built in turn on the theory developed
to prove 
Theorem~\ref{T1.-1}(i).
This issue has since been solved, through the
inscrutable vissicitudes of the publication
process: as it happens the joint work
with Burke and Pauwels appeared first, and
the two articles for which I happen to be
the sole author could
therefore be merged.

But I should take this opportunity to highlight
the value of the joint work with
Burke and Pauwels, especially
\cite[Lemma~3.9(iv)]{Burke-Neeman-Pauwels18}.
At the time we proved it it seemed,
to the three of us, to
be a small, technical lemma. But by now this
little lemma has
found diverse and surprising applications.
\ermk

\nin
From Theorem~\ref{T1.-1}(i) we deduce:

\cor{C0.-978}
Let $\ct$ be as in Theorem~\ref{T1.-1}(i),
but assume further that $\ct^c$ is contained in $\ct^b_c$.
Let $(\ct^{\leq0},\ct^{\geq0})$ be one of
the preferred {\it t}--structures.

Let $\cl:\ct^b_c\la\cs$ be an $R$--linear triangulated functor. Then the functor $\cl$ has a right adjoint
if and only if the following three conditions hold:
\be
\item
  For any pair of objects $(t,s)$, with $t\in\ct^c$ and $s\in\cs$,
  the $R$--module
  $\Hom\big(\cl(t),s\big)$ is finite.
\item
  For any object $s\in\cs$ there exists an integer $A>0$ with
  $\Hom\big(\cl(\ct^b_c\cap\ct^{\leq-A})\,,\,s\big)=0$.
\item
  For any object $t\in\ct^c$ and any object $s\in\cs$ there exists an integer
  $A$ so that $\Hom\big(\cl(\T^mt),s\big)=0$ for all $m\leq-A$.
\ee
\ecor

\nin
In the special example of $\ct=\Dqc(X)$, Theorem~\ref{T1.-1} specializes to

\exm{E1.-3}
Let $X$ be a scheme proper over a noetherian
ring $R$.
Then it is separated and quasicompact,
hence the category $\ct=\Dqc(X)$ is approximable.
But properness also guarantees that, for any compact
generator $G\in\ct$ and any $i\in\zz$, the $R$--module
$\Hom(\T^iG,G)$ is finite. The conditions of 
Theorem~\ref{T1.-1}(i) are satisfied, and
the conclusion is:
\be
\item
The functor $\cy\circ i$ gives  an
equivalence from the category
$\dcoh(X)$ to the category of finite cohomological functors
${\dperf X}\op\la\Mod R$.
\item
On the larger category $\dmcoh(X)$, the functor $\cy$
is full and the essential image is the category of locally finite
cohomological functors ${\dperf X}\op\la\Mod R$.
\setcounter{enumiv}{\value{enumi}}
\ee
If we furthermore assume that the ring $R$ is
of finite Krull dimension and excellent, 
then from Aoki~\cite[Theorem~5.1]{Aoki20}
we know that there exists an object $G'\in\ct^b_c=\dcoh(X)$
and an integer $N>0$ with $\Dqc(X)=\ogenun{G'}N$.
Thus the hypotheses of Theorem~\ref{T1.-1}(ii)
are also satisfied, and we obtain:
\be
\setcounter{enumi}{\value{enumiv}}
\item
The functor $\wt\cy\circ \wi$ gives  an
equivalence from the category
${\dperf X}\op$ to the category of
finite homological functors
$\dcoh(X)\la\Mod R$.
\item
On the larger category $\dmcoh(X)\op$, the functor $\wt\cy$
is full and the essential image is the
category of locally finite
homological functors $\dcoh(X)\la\Mod R$.
\ee
\eexm

\plm{P1.-infty}
The results in Theorem~\ref{T1.-1} and of
Example~\ref{E1.-3} have an obvious symmetry, which
leads one to wonder if there is a way to understand
and unify them.

Only very recently have I made any substantial
progress on this, and the work is still very
much in progress.
\eplm

\rmk{R1.-7}
If $R$ is a field and $X$ is projective over $R$, then the
part of 
Example~\ref{E1.-3} concerning
the objects in the image of $\cy\circ i$  and of
$\wt\cy\circ\wi$ is
known---see Bondal and Van den Bergh~\cite[Theorem~A.1]{BondalvandenBergh04} for $\cy\circ i$ and
Rouqier~\cite[Corollary 7.51(ii)]{Rouquier08}
for $\wt\cy\circ\wi$.
Neither of the old theorems 
says anything about the functors $\cy\circ i$
and $\wt\cy\circ\wi$ being fully faithful.

About the existing proofs: what Rouquier
presents in \cite[Corollary 7.51(ii)]{Rouquier08}
is the outline of how a proof might go, which
unfortunately I haven't been able fill out.
And aside from this one published sketch,
there is nothing
in the existing literature resembling a proof
of Example~\ref{E1.-3}(iv).

The existing proofs of variants
of Example~\ref{E1.-3}(ii),
including the current one, proceed in two steps.
Starting with a finite $\ct^c$--cohomological functor $H$ one
first proves that $H\cong\cy(t)$ for some $t\in\ct$, and then shows
that $t$ must actually belong to $\ct^b_c$.
Bondal and Van den
Bergh~\cite[Theorem~A.1]{BondalvandenBergh04} and
Jack Hall~\cite[Proposition~4.1]{Hall18} rely on suitable
special features that allow the functor
$H:\ct^c=\dperf X\la\Mod R$ to extend
to a cohomological functor on all of $\ct=\Dqc(X)$, and then use
the usual Brown representability theorem for $\Dqc(X)$.
For Bondal and Van den Bergh
the key is forming the double dual---this works since $R$ is assumed
a field, and a finite-dimensional
vector space over $R$ is canonically isomorphic to its double
dual. Jack Hall relies on the fact that his
functors come from morphisms of ringed
spaces $c:\mathfrak{X}\la X$, and
formal properties then provide  adjoints
\[\xymatrix@C+40pt{
\Dqc(X)\ar@<0.5ex>[r]^-{\text{natural}} &\ar@<0.5ex>[l]^-{\R Q_X^{}} \D(X)
\ar@<0.5ex>[r]^-{\LL c^*}&\ar@<0.5ex>[l]^-{\R c_*} \D(\mathfrak{X})
}\]

 We should recall one more result in the literature: although
Ben-Zvi, Nadler and Preygel~\cite[Section~3]{BenZvi-Nadler-Preygel16}
is not technically either a special case or a generalization
of  Theorem~\ref{T1.-1}, the
reader is nonetheless encouraged to look at it---there are
interesting parallels.
Enhancements
play a role in~\cite{BenZvi-Nadler-Preygel16}, as well as
the construction of an explicit generator and estimates
similar to those of
\cite[Theorem~4.1]{Lipman-Neeman07}.

What's different here is the generality. Let
$H$ be \emph{any}
locally finite $\ct^c$--cohomological functor.
Under hypotheses
weaker
than approximability (see
Proposition~\ref{P37.1005} for the precise statement)
we prove that $H\cong\cy(t)$ where $t\in\ct$ is some object---the
existence of $t$ is formal, not special 
to narrow classes of $\ct$'s or $H$'s.
And by combining a careful analysis of the
proof of Proposition~\ref{P37.1005},
with the theory developed in Section~\ref{S17}, we will
deduce---\emph{under only
the approximability hypothesis}---that $t$ must belong to $\ct^-_c$.
\ermk

\rmk{R1.-8}
The work was inspired at the time
by the lovely new proof of a vast
generalization of GAGA,
to be found in Jack Hall~\cite{Hall18}.
More precisely: it was inspired by the original idea, which is
to be found in \cite[Section~2, the section
labeled ``A simple case'']{Hall18}. As Hall's paper became
more general it developed a different tack
and, by the published version
\cite{Hall22}, the short section with the simple,
central idea went missing. 
One of the points of the current paper is that our representability
theorems obviate the need to do much to pass from
Hall's original, simple idea to a full-blown proof. 
In Example~\ref{EA.3}
the reader can find this spelt out: Appendix~\ref{SA}
is all of two pages long and gives a full proof of GAGA.

The condensed summary of Appendix~\ref{SA} is as follows.
Let $X$ be a scheme proper over the complex numbers $\C$.
With $R=\C$ we
apply Corollary~\ref{C0.-978}, with $\ct=\Dqc(X)$ 
[and hence with $\ct^b_c=\dcoh(X)$], with 
$\cs=\dcoh(X^{\text{\rm an}})$, and
where the functor $\cl:\dcoh(X)\la\dcoh(X^{\text{\rm an}})$ 
is the analytification.
The hypotheses of Corollary~\ref{C0.-978} are trivial to check, and hence
$\cl$ has a right adjoint 
$\car:\dcoh(X^{\text{\rm an}})\la\dcoh(X)$. And then,
proving that $\cl$ and $\car$ are quasi-inverses,
reduces to checking that the unit $\eta:\id\la\car\cs$
and counit $\e:\cl\car\la\id$
of adjunction are
isomorphisms, and standard-enough techniques
make this 
an easy exercise.
For details 
see Appendix~\ref{SA}.

It isn't often that one achieves such a thing: find a more
elegant, much shorter proof of a theorem by a giant like
Jean-Pierre Serre. After all: Serre isn't only a giant as a 
mathematician,
he is also a master of crisp and elegant exposition.
 Of course most of the credit goes to 
Jack Hall, he had the beautiful key idea. In this article we 
provide the technical, triangulated category framework
allowing for a straightforward and direct
passage from the simple idea to 
a complete proof.
\ermk

\rmk{R1.infty+1}
In the time since the first versions of the current
manuscript appeared on the archive, people have found
applications of the representability theorems
presented in Example~\ref{E1.-3} quite different
from the one that originally motivated the author.
Perhaps the most spectacular of these is
Bondarko~\cite{Bondarko22}, which uses the full
strength of the respresentability theorems
to extablish a bijection between semiorthogonal
decompositions of $\dperf X$ and of $\dcoh(X)$.
This amounts to a major improvement on earlier,
wonderful results
by Karmazyn, Kuznetsov
and Shinder~\cite{Karmazyn-Kuznetsov-Shinder22}.
\ermk

We have already mentioned that
part of the interest of the paper is that natural objects,
like the subcategories 
$\dcoh(X)\subset\dmcoh(X)$
of the category $\ct=\Dqc(X)$,
have an intrinsic
description. The definitions are not hard to give, we
include them in the Introduction. Before all else we recall
some standard notation.

\rmd{R1.-9}
Let $\ct$ be a triangulated category. We define
\be
\item
  If $\ca\subset\ct$ is a full subcategory, then $\Smr(\ca)$ is
  the full subcategory of all direct summands of objects of $\ca$.
\item
  If $\ca\subset\ct$ is a full subcategory, then $\add(\ca)$ is
  the full subcategory of all finite direct sums of objects of $\ca$.
\item
  If $\ct$ has small coproducts and $\ca\subset\ct$  is a
  full subcategory, then $\Add(\ca)$ is
  the full subcategory of all coproducts of objects of $\ca$.
\item
  If $\ca,\cb$ are two full subcategories of
  $\ct$, then $\ca\star\cb$ is
  the full subcategory of all objects $y\in\ct$
  such that there exists a triangle $a\la y\la b\la$ with $a\in\ca$ and
  $b\in\cb$.
\item
  Given an object $G\in\ct$ and two integers $A\leq B$, let
  $\cc\subset\ct$ be the full subcategory with objects
  $\{\T^{-i}G\mid A\leq i\leq B\}$.
  For integers $n>0$ we define the subcategories
  $\COprod_n\big(G[A,B]\big)$, inductively on the integer $n$, by the formulas
\begin{eqnarray*}
\COprod_1^{}\big(G[A,B]\big)&=&\ds\add(\cc)\ ,\\
\COprod_{n+1}^{}\big(G[A,B]\big)&=&
\COprod_1^{}\big(G[A,B]\big)\star\COprod_n^{}\big(G[A,B]\big)\ ,\\
\COprod\big(G[A,B]\big)&=&{\textstyle\bigcup_{n=1}^\infty}\,\,\COprod_n^{}\big(G[A,B]\big)\ .
\end{eqnarray*}
\item
  Given an object $G\in\ct$ and three integers $A\leq B$ and
  $n>0$,
  we define the subcategories $\genu Gn{A,B}$ by the formula
  $\genu Gn{A,B}=\Smr\big[\COprod_n^{}\big(G[A,B]\big)\big]$.
\item
  We adopt the following conventions:
  \[\begin{array}{ccc}\genul GnB=\bigcup_{A}\genu Gn{A,B}, &
  \genuf GnA=\bigcup_B^{}\genu Gn{A,B},&
  \gen Gn=\bigcup_{A\leq B}^{}\genu Gn{A,B},\\*[7pt]
  \gen G{}=\bigcup_{n>0}^{}\gen Gn, &
  \genu G{}{A,B}=\bigcup_{n>0}^{}\genu Gn{A,B}, &
  \genul G{}B=\bigcup_{A}^{}\genu G{}{A,B},\\*[7pt]
& \genuf G{}A=\bigcup_{B}^{}\genu G{}{A,B}.&
  \end{array}
\]
\item
  Suppose $\ct$ has coproducts, let $G$ be an object, and let
   $A\leq B$ be two integers. We define
  $\cc\subset\ct$ to be the full subcategory with objects
  $\{\T^{-i}G\mid A\leq i\leq B\}$.
  For integers $n>0$ we define the subcategories
  $\Coprod_n\big(G[A,B]\big)$, inductively on the integer $n$, by the formulas
\begin{eqnarray*}
\Coprod_1^{}\big(G[A,B]\big)&=&\ds\Add(\cc)\ ,\\
\Coprod_{n+1}^{}\big(G[A,B]\big)&=&
\Coprod_1^{}\big(G[A,B]\big)\star\Coprod_n^{}\big(G[A,B]\big)\ .
\end{eqnarray*}
In other words the difference between $\Coprod$ and $\COprod$
is that in $\Coprod$ we allow infinite coproducts in
the formation of $\Coprod_1$. The inductive procedure is unaltered.
\item
  We allow $A$ and $B$ to be infinite in (viii). For example
  $\Coprod_1^{}\big(G(-\infty,B]\big)$ is defined
  to be $\Add(\cc)$ with $\cc=\{\T^{-i}G\mid i\leq B\}$.
\item
  Let $A\leq B$ be integers, possibly infinite. Then $\Coprod\big(G[A,B]\big)$
  is the smallest full subcategory $\cs\subset\ct$, closed under coproducts,
  with $\cs\star\cs\subset\cs$, and with $\T^{-i}G\in\cs$ for
  $A\leq i\leq B$.
\item
  For triples of
  integers $A\leq B$ and $n>0$, we let
  $\ogenu Gn{A,B}=\Smr\big[\Coprod_n\big(G[A,B]\big)\big]$.
  In this formula we also allow $A$ and $B$ to be infinite. 
\item
  For pairs of integers $A\leq B$ we
  let $\ogenu G{}{A,B}=\Smr\big[\Coprod\big(G[A,B]\big)\big]$.
  In this formula we also allow $A$ and $B$ to be infinite, but as
  it happens for infinite $A$ we obtain nothing new. The
  categories
  \[\Coprod\big(G(-\infty,B]\big),\qquad
  \Coprod\big(G(-\infty,\infty)\big)\]
  are closed under coproducts and (positive) suspensions, and therefore
  contain all direct summands of their objects.
\ee
\ermd

The following lemma is an easy consequence of the definitions.

\lem{L0.-9238}
Suppose $G,H$ are objects in a triangulated category $\ct$. We show
\be
\item
  If $H\in\gen G{}$ then there exists an integer $A>0$ with
  $H\in\genu GA{-A,A}$.
\item
  If $\gen G{}=\gen H{}$ then there exists 
an integer $A>0$ with $H\in\genu GA{-A,A}$ and $G\in\genu HA{-A,A}$.
\ee
\elem

\prf
For (i) the assumption is $H\in\gen G{}=\cup_{A>0}^{}\genu  GA{-A,A}$, hence
$H$ belongs to one of the sets in the union. For (ii) observe that
$\gen G{}=\gen H{}$ implies $H\in\gen G{}$ and $G\in\gen H{}$ and apply (i).
\eprf

Now we come to the first new definition.

\dfn{D0.1}
Suppose we are given two {\it t}--structures on a triangulated
category $\ct$, that is we are given two pairs of
subcategories $(\ct_1^{\leq0},\ct_1^{\geq0})$ and
$(\ct_2^{\leq0},\ct_2^{\geq0})$ satisfying the conditions in
\cite[D\'efinition~1.3.1]{BeiBerDel82}.
These {\it t}--structures are \emph{equivalent}
if and only if there exists an integer $A>0$ with
$\ct_1^{\leq-A}\subset\ct_2^{\leq0}\subset\ct_1^{\leq A}$.
\edfn

\obs{O0.2}
For any \tstr\ $(\ct^{\leq0},\ct^{\geq0})$ we
have $\ct^{\leq0}=\Tm\big({^\perp\ct^{\geq0}}\big)$ and
$\ct^{\geq0}=(\T\ct^{\leq0})^\perp$. It immediately
follows that two {\it t}--structures
$(\ct_1^{\leq0},\ct_1^{\geq0})$ and
$(\ct_2^{\leq0},\ct_2^{\geq0})$
are equivalent if and only if there exists an integer $A>0$ with
$\ct_1^{\geq A}\subset\ct_2^{\geq0}\subset\ct_1^{\geq-A}$.
\eobs

\obs{O0.2.5}
Recall that, for any \tstr\ $(\ct^{\leq0},\ct^{\geq0})$, the categories
$\ct^-$, $\ct^+$ and $\ct^b$ are defined by
\[
\ct^-=\cup_{m>0}^{}\ct^{\leq m}\ ,\qquad
\ct^+=\cup_{m>0}^{}\ct^{\geq -m}\ ,\qquad
\ct^b=\ct^-\cap\ct^+\ .
\]
If
$(\ct_1^{\leq0},\ct_1^{\geq0})$ and
$(\ct_2^{\leq0},\ct_2^{\geq0})$
are equivalent {\it t}--structures we note
\be
\item
  $\ct_1^-=\ct_2^-$, $\ct^+_1=\ct^+_2$ and $\ct^b_1=\ct^b_2$.
\item
  If $\ct^-$ [respectively $\ct^+$, respectively $\ct^b$] contains a
  compact generator $G\in\ct^c$, then
  $\ct^-$ [respectively $\ct^+$, respectively $\ct^b$] contains all
  of
  $\ct^c$.
\ee
\eobs

\prf
We prove (i) and (ii) for $\ct^-$ and leave $\ct^+$ and $\ct^b$
to the reader. To prove (i) observe that the inclusions
$\ct_1^{\leq-A}\subset\ct_2^{\leq0}\subset\ct_1^{\leq A}$ imply
\[
\bigcup_{m>0}^{}\ct_1^{\leq-A+m}\sub\bigcup_{m>0}^{}\ct_2^{\leq m}\sub
\bigcup_{m>0}^{}\ct_1^{\leq A+m}
\]
that is $\ct_1^-\subset\ct_2^-\subset\ct_1^-$.

For the proof of (ii) the assumption is
that $G\in\ct^-$.
This makes $\ct^-\subset\ct$ a thick subcategory containing
$G$, hence $\ct^c=\gen G{}\subset\ct^-$.
\eprf

\exm{E0.3}
Let $\ct$ be a triangulated category with coproducts. Given any
compact object $G\in\ct$, from Alonso, Jerem{\'{\i}}as and
Souto~\cite[Theorem~A.1]{Alonso-Jeremias-Souto03} we
learn that $\ct$ has a unique \tstr\ $(\ct_G^{\leq0},\ct_G^{\geq0})$
\emph{generated} by $G$.
In the notation of Reminder~\ref{R1.-9}, the aisle $\ct_G^{\leq0}$
of this \tstr\ is
nothing other than $\ct_G^{\leq0}=\ogenul G{}0$.
It follows formally that both $\ct_G^{\leq0}$ and $\ct_G^{\geq0}$ are closed
under coproducts and direct summands---the closure under direct
summands is true for any aisle and co-aisle of a \tstr, the closure
of $\ct_G^{\leq0}$ under coproducts is also true for any aisle,
while
the fact that $\ct_G^{\geq0}$ is closed under coproducts may be found in
\cite[Proposition~A.2]{Alonso-Jeremias-Souto03}; it comes from the
compactness of the object $G$.

If $G,H$ are two compact objects of $\ct$ with $\gen G{}=\gen H{}$,
Lemma~\ref{L0.-9238}(ii) tells us that
there exists an integer $A>0$ with $H\in\genu GA{-A,A}$
and $G\in\genu HA{-A,A}$. Hence
$\ogenul H{}{-A}\subset \ogenul G{}0\subset \ogenul H{}A$, that is
$\ct_H^{\leq-A}\subset\ct_G^{\leq0}\subset\ct^{\leq A}_H$.
Thus the {\it t}--structures generated by $G$ and $H$
are equivalent. This leads us to
\eexm

\dfn{D0.5}
If the compactly
generated triangulated category $\ct$ has a
single compact object $G$
that generates it, then
the \emph{preferred equivalence class of {\it t}--structures} is
the one containing the \tstr\
$\big(\ct_G^{\leq0},\ct_G^{\geq0}\big)$  generated by $G$. 
\edfn

\rmk{R0.7}
For any compact generator $G$ we have that $\gen G{}=\ct^c$, the full
subcategory of all compact objects. Any two compact generators $G,H$ satisfy
$\gen G{}=\ct^c=\gen H{}$, and Example~\ref{E0.3} says that $G$ and $H$ generate
equivalent {\it t}--structures. Thus the preferred equivalence class
of {\it t}--structures does not depend on the choice of compact generator.

Now
\cite[Proposition~A.2]{Alonso-Jeremias-Souto03} guarantees
that, in the preferred equivalence class, there will exist some
{\it t}--structures with
$\ct^{\leq0}$ and $\ct^{\geq0}$ both closed under coproducts---just
take $\big(\ct_G^{\leq0},\ct_G^{\geq0}\big)$ for a compact
generator $G$.
The reader should note that this property is \emph{not} stable under
equivalence. 
In general there will be {\it t}--structures in the preferred
equivalence class where $\ct^{\geq0}$
is not closed in $\ct$ under coproducts.

From Observation~\ref{O0.2.5}(i) we learn that, as long as we
stick to the preferred equivalence class of
{\it t}--structures, the categories $\ct^-$,
$\ct^+$ and $\ct^b$ are intrinsic.
\ermk

And now for the next formal construction.

\dfn{D0.13}
Suppose $\ct$ is a triangulated category with coproducts
and let $(\ct^{\leq0},\ct^{\geq0})$ be a \tstr.

An object $F$ belongs to the subcategory $\ct_c^-\subset\ct$ if, for
any 
integer $m>0$, there exists a triangle
$E\la F\la D$ with $E\in\ct^c$ and $D\in\ct^{\leq-m}$.

The subcategory $\ct^b_c$ is defined by $\ct^b_c=\ct^-_c\cap\ct^b$.
\edfn

\rmk{R0.14}
Note that the definition of $\ct^-_c$ depends on the choice of a \tstr,
but not much---equivalent {\it t}--structures lead to the same
$\ct^-_c$. For any choice of \tstr\ the category
$\ct^-_c$ contains $\ct^c$. After all if $F$ is compact
then the triangle $F\stackrel\id\la F\la 0$ has $F\in\ct^c$
and $0\in\ct^{\leq-m}$, for every $m$ and every \tstr.
\ermk

\rmk{R0.13.99}
Assume the  \tstr\ $(\ct^{\leq0},\ct^{\geq0})$
is such that there is a compact generator
$G$ contained in $\ct^-$;
any \tstr\ in the preferred equivalence class
is an example, after all $G\in\ogenul G{}0=\ct_G^{\leq0}\subset\ct^-$.
Observation~\ref{O0.2.5}(ii)
gives that $\ct^c\subset\ct^-$, and 
 Definition~\ref{D0.13}
tells us that, for any integer $m>0$,
\[
\ct^-_c\sub\ct^c*\ct^{\leq-m}\sub \ct^-*\ct^-\eq\ct^-\ .
\]
\ermk

Still in gorgeous generality we will prove

\pro{P0.14.5}
Let $\ct$ be a triangulated category with coproducts, and
let $(\ct^{\leq0},\ct^{\geq0})$ be a \tstr. If there
exists an integer $A>0$ and a compact generator
$G\in\ct$ with $\Hom\big(\T^{-A}G,\ct^{\leq0}\big)=0$
then $\ct^b_c\subset\ct^-_c$ are triangulated subcategories of $\ct$.
If furthermore $G\in\ct^-$,
then $\ct^b_c\subset\ct^-_c\subset\ct$ are thick
subcategories of $\ct^-$.
\epro

\rmk{R0.14.7}
We are most interested in the special case 
where the \tstr\ $(\ct^{\leq0},\ct^{\geq0})$
is in the preferred equivalence
class and $\ct^b_c\subset\ct^-_c$ are independent of choices.

Suppose there exists a compact generator $G$ and an integer $A>0$,
so that $\Hom(G,\T^iG)=0$ for
all $i\geq A$. Define the full subcategory $\cs$ by
\[
\cs\eq \{S\in\ct\mid\Hom(\T^{-A}G,S)=0\}.
\]
The compactness of $G$
says that $\cs$ is closed under coproducts, by hypothesis
$\cs$ contains $\T^iG$ for all $i\geq0$, while obviously $\cs$ is
closed under direct summands and $\cs*\cs\subset\cs$.
Therefore $\cs$ contains $\ogenul G{}0=\ct^{\leq0}_G$.
We deduce that $\Hom\big(\T^{-A}G,\ct^{\leq0}_G)=0$.
Since $G$ is obviously in $\ct^{\leq0}_G\subset\ct^-$,
Proposition~\ref{P0.14.5} informs us that $\ct^b_c\subset\ct^-_c$ are
thick subcategories of $\ct^-$.
\ermk

For the structure defined so far we needed very little. To go
further it turns out to be useful to estimate how much
effort it takes to approximate an object in $\ct^-$ by
a compact generator $G$. This leads us to

\dfn{D0.9}
Let $\ct$ be a triangulated category with
coproducts.
The category $\ct$ is called
\emph{weakly approximable} if
there exists a compact generator $G$,  a
\tstr\ $(\ct^{\leq0},\ct^{\geq0})$
and an integer $A>0$ so that 
\be
\item
$\T^AG\in\ct^{\leq0}$ and $\Hom\big(\T^{-A}G,\ct^{\leq0}\big)=0$.
\item
Every object $F\in\ct^{\leq 0}$ admits a triangle
$E\la F\la D$ with $E\in\ogenu G{}{-A,A}$ and
$D\in\ct^{\leq-1}$.
\setcounter{enumiv}{\value{enumi}}
\ee
The category $\ct$ is called \emph{approximable} if the integer $A$
can be chosen to further satisfy
\be
\setcounter{enumi}{\value{enumiv}}
\item
In the triangle 
$E\la F\la D$ of (ii) above we may  
strengthen the condition on $E$, we may assume
$E\in\ogenu GA{-A,A}\subset\ogenu G{}{-A,A}$.
\setcounter{enumiv}{\value{enumi}}
\ee
\edfn

The following are easy to prove, they will be part of a string of
formal consequences of approximability, see Section~\ref{S17}.

\fac{F0.11}
Let $\ct$ be a triangulated category with coproducts.
If $\ct$ is weakly approximable then
\be
\item
  The \tstr\ $(\ct^{\leq0},\ct^{\geq0})$, which is part of
  Definition~\ref{D0.9} and is assumed to satisfy some
  hypotheses, must belong to the preferred equivalence
  class.
\item
For any compact generator $G$ and any \tstr\
$(\ct^{\leq0},\ct^{\geq0})$ in the preferred equivalence class
  there must exist an integer $A$, depending on $G$ and
on the \tstr\ $(\ct^{\leq0},\ct^{\geq0})$, which satisfies 
Definition~\ref{D0.9}(i) and (ii). If $\ct$ is approximable the
integer $A$ may be chosen to satisfy (iii) as well.
\ee
Thus in proving that $\ct$ is (weakly) approximable we can choose our compact
generator and \tstr\ to suit our convenience. Once we know the category
is approximable, it follows that the convenient \tstr\ is in the
preferred class, and any compact generator and any \tstr\ in
the preferred equivalence class fulfill the approximability
criteria.
\efac

\fac{F0.15}
As stated in the first few paragraphs of the introduction [before we
  presented the definitions] we will prove that,
if $X$ is a quasicompact, separated scheme,
then $\ct=\Dqc(X)$ is approximable  and the standard
\tstr\ is in the preferred equivalence class.
If $X$ is only quasiseparated, then the reader is
referred to \cite[Theorem~3.2~(iii) and (iv)]{Neeman22A}
for the fact that
$\ct=\Dqc(X)$ is weakly approximable and
the standard
\tstr\ is in the preferred equivalence class.
In the special case where $X$ is noetherian then
$\ct^b_c\subset\ct^-_c$ are just $\dcoh(X)\subset\dmcoh(X)$.
For non-noetherian $X$ the description
of $\ct^b_c\subset\ct^-_c$ is slightly more complicated,
but still classical---see Remark~\ref{R19.8976}
for details.
The fact that the standard \tstr\ is in the preferred
equivalence class tells us that $\ct^-=\Dqcmi(X)$, $\ct^+=\Dqcpl(X)$
and $\ct^b=\Dqcb(X)$.

Another example
is the homotopy category $\ct$ of spectra.
In this case we can take $\ct^{\leq0}\subset\ct$ to be the
subcategory of connective spectra---the \tstr\ this defines
is in the preferred equivalence class. The category $\ct$ turns out
to be approximable, and the subcategory $\ct^-_c$ is the category
of spectra $X$ whose stable homotopy groups $\pi_i(X)$
are finitely generated
$\zz$--modules and $\pi_i(X)=0$ for $i\ll0$.
And $\ct^b_c\subset\ct^-_c$ is the subcategory where all but finitely
many of the $\pi_i(X)$ vanish.

The representability we prove in Theorem~\ref{T1.-1} applies to this
example but the result is not new. There is a theorem of
Adams~\cite{Adams71} which says that every cohomological functor $H$
on $\ct^c$ is the restriction of a representable one on $\ct$,
and it is easy to show that finiteness or local finiteness
of $H$ translate to saying that the representing object must lie in
$\ct^b_c$ or $\ct^-_c$. But the theorem of Adams does not generalize
to $\dperf X\subset\Dqc(X)$; see
\cite{Neeman4,Christensen-Keller-Neeman99}.

If $X$ is a quasicompact, quasiseparated scheme and
$Z\subset X$ is a closed subset with
quasicompact complement, then the category 
$\ct={\mathbf D_{\text{\bf qc},Z}}(X)$, the subcategory of
$\Dqc(X)$ of all complexes supported
on $Z$,
turns out to be
weakly approximable but not (in general) approximable.
The proof of the weak approximability is
in \cite[Theorem~3.2(iv)]{Neeman22A}, and
an example showing that approximability fails in
general is given in \cite[Remark~8.1]{Neeman22A}.
In \cite[Theorem~3.2(iii)]{Neeman22A} the reader will
learn that the standard \tstr\ is in the preferred
equivalence class. If $X$ is
noetherian the categories $\ct_c^-$
and $\ct_c^b$ are (respectively) the intersections of
$\dmcoh(X)$ and $\dcoh(X)$ with the category 
$\mathbf
D_{\text{\bf qc},Z}(X)$;
the reader can find this in
\cite[second paragraph of Section~8]{Neeman22A}.
\efac

The definitions have all been made and the reader can go back
to the statements of Theorem~\ref{T1.-1}
and Corollary~\ref{C0.-978}, which are now precise.
Note that in both results $\ct$ has to be approximable, weakly
approximable is not enough.

We have discussed what we know, but should point out that there are many
more potential examples. After all: let $R$ be a commutative
ring and let $T$ be a dg $R$--algebra. Then the category
$\ct=\D(\Mod T)$ is a triangulated category with coproducts and
a single
compact generator $T$. It has a preferred equivalence
class of {\it t}--structures, one can define the
intrinsic subcategories $\ct^-$, $\ct^+$, $\ct^b$, $\ct^-_c$ and $\ct^b_c$,
and in general I have no idea what they are. If $H^i(T)=0$ for $i\gg0$ then the
subcategories $\ct^-_c$ and $\ct^b_c$ are thick, this follows
from Remark~\ref{R0.14.7}. If $H^i(T)=0$ for all $i\geq1$ we are in the trivial
case (see Remark~\ref{R19.303}),
where it's easy to prove the category $\D(\Mod T)$ to be approximable
and work out explicitly what
are $\ct^-$, $\ct^+$, $\ct^b$, $\ct^-_c$ and $\ct^b_c$.
More generally: if $T$ is a dg algebra, with $H^i(T)=0$
for all $i\geq2$, then Bondarko
and Vostokov~\cite[Corollary~4.3]{Bondarko-Vostokov20}
proves that the derived category $\D(\Mod T)$
is weakly approximable.
So far the only other general result, producing further examples
of approximable
triangulated categories, is \cite[Theorem~4.1]{Burke-Neeman-Pauwels18}.
It says that, under reasonable hypotheses, the recollement of two approximable
triangulated categories is approximable. 
But for $T$ a general dga, satisfying $H^i(T)=0$ for $i\gg0$,
I have no idea when the categories
$\D(\Mod T)$ are approximable or weakly approximable.
In view of
Theorem~\ref{T1.-1} and
Corollary~\ref{C0.-978} it would be interesting to find out, especially
since the categories $\D(\Mod T)$ are of
so much current active interest---their
study is at the core of noncommutative algebraic geometry. Who knows, there
might be a noncommutative generalization of GAGA.

\medskip

\nin
    {\bf Acknowledgements.}\ \ The author is grateful to an
    anonymous referee, for his careful reading of an earlier version, for pointing out many errors, and for suggesting valuable expository improvements.

\section{Basics}
\label{S-1}

Since {\it t}--structures will play a big part in the article
we begin with a quick reminder of some elementary facts.

\rmd{R-1.1}
In this section $\ct$ will be a triangulated category and
$(\ct^{\leq0},\ct^{\geq0})$
will be a \tstr\ on $\ct$. The category $\ca=\ct^{\leq0}\cap\ct^{\geq0}$
is abelian, it is called the heart of the \tstr. The functor
$\ch(-)=\big[(-)^{\leq0}\big]^{\geq0}$ is a homological functor $\ch:\ct\la\ca$.
We will let $\ch^\ell$ be the functor
$\ch^\ell(-)=\ch\big[\T^\ell(-)\big]=\T^\ell\big[(-)^{\leq\ell}\big]^{\geq\ell}$.
\ermd

\lem{L-1.3}
Let $\ct$ be a triangulated category and let $(\ct^{\leq0},\ct^{\geq0})$
be a \tstr\ on $\ct$. If $F$ is an object of $\ct^-$, and
$\ch^\ell(F)=0$ for all $\ell>-i$, then $F$ belongs to 
$\ct^{\leq -i}$.
\elem

\prf
We are given that $F$ belongs to $\ct^-=\cup_n^{}\ct^{\leq n}$, hence
$F\in\ct^{\leq n}$ for some $n$ and the map
$F^{\leq n}\la F$ is an isomorphism. But
now the triangle $F^{\leq\ell-1}\la F^{\leq\ell}\la \T^{-\ell}\ch^\ell(F)$
informs us that, as long as $\ell>-i$, the map
$F^{\leq\ell-1}\la  F^{\leq\ell}$
is also an isomorphism.
Composing the string of isomorphisms
$F^{\leq-i}\la F^{\leq-i+1}\la\cdots\la F^{\leq n}\la F$
we have
that $F^{\leq-i}\la F$
is an isomorphism---therefore $F\in\ct^{\leq-i}$.
\eprf

\lem{L-1.5}
If there is an integer $A$ and a generator $G\in\ct$ with
$\Hom\big(G,\ct^{\leq-A}\big)=0$, then
\be
\item
  Any object $F\in\ct^-$,
with $\ch^\ell(F)=0$ for all $\ell$, must vanish.
\item
If $f:E\la F$ is a morphism in $\ct^-$ such that $\ch^\ell(f)$ is an
isomorphism for every $\ell\in\zz$, then $f$ is an isomorphism.
\ee
\elem

\prf
To prove (i) assume $\ch^\ell(F)=0$ for
all $\ell$;
Lemma~\ref{L-1.3} says that $F$ belongs to $\cap_\ell^{}\ct^{\leq\ell}$.
But then $\Hom(\T^iG,F)=0$ for all $i\in\zz$, and as $G$ is a generator
this implies $F=0$.

(ii) follows by applying (i) to the mapping cone of $f$.
\eprf

\lem{L-1.7}
Suppose the category $\ct$ has coproducts, and the \tstr\ is such that
both $\ct^{\leq0}$ and $\ct^{\geq0}$ are closed under the coproducts of $\ct$.
Then:
\be
\item
  The functors $(-)^{\leq0}$ and $(-)^{\geq0}$ both respect coproducts.
\item
  The heart $\ca\subset\ct$ is closed in $\ct$ under coproducts,
  and the
  functor $\ch:\ct\la\ca$ respects coproducts.
\item
  The abelian category
  $\ca$ satisfies [AB4], that is coproducts are exact.
\item
  If $E_1\la E_2\la E_3\la\cdots$ is a sequence of objects and morphisms
  in $\ct$, then 
there is a short exact
sequence in the heart $\ca$ of the \tstr
\[\xymatrix{
  0 \ar[r] & 
  \colim\ch^\ell(E_i) \ar[rr]& &\ch^\ell\Big(\hoco E_i\Big)\ar[rr] &&
  \colim^1\ch^{\ell+1}(E_i)\ar[r] & 0
}\]
\ee
\elem

\prf
Suppose we are given in $\ct$
a collection of objects $\{E_\lambda^{},\,\lambda\in\Lambda\}$. For
each $\lambda$ we have a canonical triangle
$E_\lambda^{\leq0}\la E_\lambda^{}\la E_\lambda^{\geq1}\la\T E_\lambda^{\leq0}$.
The coproduct of these triangles is a triangle
\[\xymatrix{
  \ds \bigoplus_{\lambda\in\Lambda}E_\lambda^{\leq0}
  \ar[r] &\ds \bigoplus_{\lambda\in\Lambda} E_\lambda^{}\ar[r]
  &\ds \bigoplus_{\lambda\in\Lambda} E_\lambda^{\geq1}
  \ar[r] &\ds \bigoplus_{\lambda\in\Lambda}\T E_\lambda^{\leq0}
}\]
By hypothesis $\oplus^{}_{\lambda\in\Lambda}E_\lambda^{\leq0}$
belongs to $\ct^{\leq0}$ and  $\oplus^{}_{\lambda\in\Lambda}E_\lambda^{\geq1}$
belongs to $\ct^{\geq1}$, and the triangle above
must be canonically isomorphic
to
\[\xymatrix{
  \ds \left(\bigoplus_{\lambda\in\Lambda}E_\lambda^{}\right)^{\leq0}
  \ar[r] &\ds \bigoplus_{\lambda\in\Lambda} E_\lambda^{}\ar[r]
  &\ds \left(\bigoplus_{\lambda\in\Lambda} E_\lambda^{}\right)^{\geq1}
  \ar[r] &\ds \left(\bigoplus_{\lambda\in\Lambda}\T E_\lambda^{}\right)^{\leq0}
}\]
This proves (i).

Since $\ct^{\leq0}$ and $\ct^{\geq0}$ are closed in $\ct$ under coproducts
so is their intersection $\ca=\ct^{\leq0}\cap\ct^{\geq0}$. By (i) we know that
the functors $(-)^{\leq0}$ and $(-)^{\geq0}$ both
respect coproducts, hence so does their
composition $\ch(-)=\big[(-)^{\leq0}\big]^{\geq0}$. This proves (ii).

The category $\ct$ has coproducts and its subcategory $\ca$ is closed
under these coproducts, hence $\ca$ has coproducts---it satisfies
[AB3]. Now suppose we are given a set
$\{f_\lambda^{}:A_\lambda^{}\sr B_\lambda^{},\,\lambda\in\Lambda\}$
of morphisms
in $\ca$. Complete these to triangles
$A_\lambda^{}\la B_\lambda^{}\la C_\lambda^{}\la\T A_\lambda^{}$ and form
the coproduct
\[\xymatrix{
  \ds\bigoplus_{\lambda\in\Lambda}A_\lambda^{}
  \ar[rr]^-{\bigoplus_{\lambda\in\Lambda}f_\lambda^{}} & &
 \ds\bigoplus_{\lambda\in\Lambda}B_\lambda^{}
  \ar[r] &
 \ds\bigoplus_{\lambda\in\Lambda}C_\lambda^{}
 \ar[r] &
  \ds\bigoplus_{\lambda\in\Lambda}\T A_\lambda^{}
}\]
which is a triangle. The long exact sequence obtained by
applying $\ch$
to this triangle tells us that the kernel
of the map $\bigoplus_{\lambda\in\Lambda}^{}f_\lambda^{}$ is
$\ch^{-1}\big(\bigoplus_{\lambda\in\Lambda}^{}C_\lambda^{}\big)$,
but (ii) informs us that this is
$\bigoplus_{\lambda\in\Lambda}^{}\ch^{-1}(C_\lambda^{})$, which
is $\bigoplus_{\lambda\in\Lambda}^{}\mathrm{Ker}(f_\lambda^{})$.
The right exactness of coproducts is formal,
 completing the proof of (iii).

 Finally (iv) follows by applying the functor $\ch$ to
 the triangle
\[\xymatrix{
  \ds\bigoplus_{i=1}^\infty E_i\ar[r] &
  \ds\bigoplus_{i=1}^\infty E_i\ar[r] &
\hoco E_i \ar[r] &
  \ds\bigoplus_{i=1}^\infty \T E_i
}\]
and using (ii) to compute the long exact sequence.
\eprf

\rmk{R-1.9}
Remark~\ref{R0.7} tells us that, if $\ct$ is a triangulated category
with coproducts and a single compact generator, then the preferred
equivalence class contains {\it t}--structures $(\ct^{\leq0},\ct^{\geq0})$
with $\ct^{\leq0}$ and $\ct^{\geq0}$ both closed under coproducts.
This is the situation in which we will apply Lemma~\ref{L-1.7}.
Note also that Remark~\ref{R0.7} warns us that not every \tstr\ in
the preferred equivalence class needs to satisfy this property.

We will mostly use Lemma~\ref{L-1.7}(iv) in the special case where
the sequences $\ch^\ell(E_1)\la\ch^\ell(E_2)\la\ch^\ell(E_3)\la\cdots$
eventually stabilize for every $\ell$.
When this happens the $\colim^1$ terms all vanish,
and the natural map is an isomorphism
$\colim\, \ch^\ell(E_i)\la
\ch^\ell\big(\hoco E_i\big)$.
\ermk

\rmk{R-1.909}
We should note that, in the special case where the
\tstr\ is compactly generated (for example the \tstr\
generated by a single compact object), then much more
is now known. As the reader can see in 
Saor{\'{\i}}n and
{\v{S}}{\v{t}}ov{\'{\i}}{\v{c}}ek~\cite[Theorem~8.31]{Saorin-Stovicek23}, the heart of a compactly generated
\tstr\ has to be a locally finitely
presented Grothendieck abelian category.

The proof of this recent theorem, by Saor{\'{\i}}n and
{\v{S}}{\v{t}}ov{\'{\i}}{\v{c}}ek, is substantially more
involved than the simple-minded, short argument we
presented in this section. And for us, in this
article, the results of the current section
amply cover what we will need.
\ermk

\section{The fundamental properties of approximability}
\label{S17}

\lem{L17.1}
Let $\ct$ be a triangulated category
with a \tstr\ $(\ct^{\leq0},\ct^{\geq0})$, and let $\cs\subset\ct$ be a full
subcategory with $\T\cs=\cs$. Assume $\ca$ is also a full subcategory
of $\ct$, and define $\ca(m)$ inductively by
\be
\item
  $\ca(1)=\ca$.
\item
  $\ca(m+1)=\ca(m)\star\T^{m}\ca$.
\ee
Suppose every object in $F\in\cs\cap\ct^{\leq0}$ admits a
triangle $E_1\la F\la D_1$,
with $E_1\in\ca$ and $D_1\in\cs\cap\ct^{\leq-1}$. Then
we can construct a sequence $E_1\la E_2\la E_3\la\cdots$, with a
map from the sequence to $F$ and so that, if we complete $E_m\la F$
to a triangle $E_m\la F\la D_m$,
then $E_m\in\ca(m)$ and $D_m\in\cs\cap\ct^{\leq-m}$.
\elem

\prf
We are given the case $m=1$; assume we have constructed the
sequence as far as an integer $m>0$, and we want to extend it to $m+1$.
Take
any object $F\in\cs\cap\ct^{\leq0}$, and by the inductive hypothesis
construct the sequence up to $m$. In particular
choose a triangle $E_m\la F\la D_m$ with
$E_m\in\ca(m)$ and $D_m\in\cs\cap\ct^{\leq-m}$.
Now apply the case $m=1$ to $\T^{-m}D_m$; we produce
a triangle $E'\la D_m\la D_{m+1}$ with 
$D_{m+1}\in\cs\cap\ct^{\leq-m-1}$ and
$E'\in\T^m\ca$.
Form an octahedron from the
composable morphisms $F\la D_m\la D_{m+1}$, that is
\[\xymatrix{
E_m\ar[d] \ar@{=}[r] & E_m\ar[d] & \\
E_{m+1}\ar[d]  \ar[r] & F\ar[d]\ar[r] & D_{m+1}\ar@{=}[d]\\
E' \ar[r] & D_m\ar[r] & D_{m+1}
}\]
The object $D_{m+1}$ belongs to $\cs\cap\ct^{\leq-m-1}$ by construction.
The triangle $E_m\la E_{m+1}\la E'$ tells us that
$E_{m+1}\in\ca(m)\star\T^m\ca=\ca(m+1)$, and we have factored the map
$E_m\la F$ as $E_m\la E_{m+1}\la F$ so that, in the triangle
$E_{m+1}\la F\la D_{m+1}$, we have $E_{m+1}\in\ca(m+1)$ and
$D_{m+1}\in\cs\cap\ct^{\leq-m-1}$.
\eprf

\cor{C17.3}
Let $\ct$ be a triangulated category with coproducts, let $G\in\ct$
be an object, and let $(\ct^{\leq0},\ct^{\geq0})$ be a \tstr. The following
is true.
\sthm{C17.3.1}
Suppose every object $F\in\ct^{\leq0}$ admits a triangle $E_1\la F\la D_1$,
with $E_1\in\ogenu G{}{-A,A}$ and $D_1\in\ct^{\leq-1}$. Then
we can extend to a sequence $E_1\la E_2\la E_3\la\cdots$, with a
map from the sequence to $F$ and so that, if we complete $E_m\la F$
to a triangle $E_m\la F\la D_m$,
then $E_m\in\ogenu G{}{1-m-A,A}$ and $D_m\in\ct^{\leq-m}$.
\esthm
\sthm{C17.3.2}
Suppose every object $F\in\ct^{\leq0}$ admits a triangle $E_1\la F\la D_1$,
with $E_1\in\ogenu G{A}{-A,A}$ and $D_1\in\ct^{\leq-1}$.
Then
we can extend to a sequence $E_1\la E_2\la E_3\la\cdots$, with a
map from the sequence to $F$ and so that, if we complete $E_m\la F$
to a triangle $E_m\la F\la D_m$,
then $E_m\in\ogenu G{mA}{1-m-A,A}$ and $D_m\in\ct^{\leq-m}$.
\esthm
\sthm{C17.3.3}
For a full subcategory $\cs\subset\ct$ with $\T\cs=\cs$,
suppose every object $F\in\cs\cap\ct^{\leq0}$ admits a triangle $E_1\la F\la D_1$,
with $E_1\in\genu G{}{-A,A}$ and $D_1\in\cs\cap\ct^{\leq-1}$.
 Then
we can extend to a sequence $E_1\la E_2\la E_3\la\cdots$, with a
map from the sequence to $F$ and so that, if we complete $E_m\la F$
to a triangle $E_m\la F\la D_m$,
then $E_m\in\genu G{}{1-m-A,A}$ and $D_m\in\cs\cap\ct^{\leq-m}$.
\esthm
\sthm{C17.3.4}
For a full subcategory $\cs\subset\ct$ with $\T\cs=\cs$,
suppose every object $F\in\cs\cap\ct^{\leq0}$ admits a triangle $E_1\la F\la D_1$,
with $E_1\in\genu G{A}{-A,A}$ and $D_1\in\cs\cap\ct^{\leq-1}$.
 Then
we can extend to a sequence $E_1\la E_2\la E_3\la\cdots$, with a
map from the sequence to $F$ and so that, if we complete $E_m\la F$
to a triangle $E_m\la F\la D_m$,
then $E_m\in\genu G{mA}{1-m-A,A}$ and $D_m\in\cs\cap\ct^{\leq-m}$.
\esthm
\ecor

\prf
In each case we apply Lemma~\ref{L17.1} with a suitable choice of $\ca$ and
$\cs$.

To prove (\ref{C17.3.1}) let $\cs=\ct$ and let $\ca=\ogenu G{}{-A,A}$.
By induction we see that $\ca(m)\subset\ogenu G{}{1-m-A,A}$ and the result
follows.

To prove (\ref{C17.3.2}) let $\cs=\ct$ and let $\ca=\ogenu G{A}{-A,A}$.
By induction we see that $\ca(m)\subset\ogenu G{mA}{1-m-A,A}$ and the result
follows.

To prove (\ref{C17.3.3})  let $\ca=\genu G{}{-A,A}$.
By induction we see that $\ca(m)\subset\genu G{}{1-m-A,A}$ and the result
follows.

To prove (\ref{C17.3.4}) let $\ca=\genu G{A}{-A,A}$.
By induction we see that $\ca(m)\subset\genu G{mA}{1-m-A,A}$ and the result
follows.
\eprf

\lem{L17.5}
Suppose $\ct$ is a compactly generated triangulated category, $G$ is a
compact generator and  $(\ct^{\leq0},\ct^{\geq0})$ a \tstr. Suppose
there exists an integer
$B$ with $\Hom(\T^{-B}G,\ct^{\leq0})=0$.

With any sequence $E_1\la E_2\la E_3\la\cdots$ mapping to $F$,
and such that in the triangles $E_m\la F\la D_m$ we have $D_m\in\ct^{\leq-m}$,
the (non-canonical) map
$\hoco E_m\la F$ is an isomorphism.
\elem

\prf
For any $n\geq0$ we have $\ct^{\leq-n}\subset\ct^{\leq0}$,
hence
$\Hom(\T^{-B}G,\ct^{\leq-n})=0$. By shifting we deduce that
$\Hom(\T^{-\ell} G,\ct^{\leq-m})=0$ as long as $m+\ell\geq B$.

The triangle $E_m\la F\la D_m$, with $D_m\in\ct^{\leq-m}$, tells us that if
$m>\max(1,B-\ell)$ then the functor $\Hom(\T^{-\ell} G,-)$ takes the map
$E_m\la F$ to an isomorphism. Now \cite[Lemma~2.8]{Neeman96}, applied
to the compact object $G\in\ct$ and the map from the sequence
$\{E_m\}$ to $F$, tells us that $\Hom(\T^{-\ell} G,-)$ takes
the map $\hoco E_m\la F$ to an isomorphism. But $G$ is a generator,
hence the map $\hoco E_m\la F$ must be an isomorphism.
\eprf

\pro{P17.7}
Suppose the triangulated
category $\ct$, the generator $G$ and
the \tstr\ $(\ct^{\leq0},\ct^{\geq0})$ are as in the hypotheses of
weakly approximable categories of Definition~\ref{D0.9}.
We remind the reader: $\ct$ has coproducts,
$G$ is a compact generator, and there is an integer $A>0$ so that 
\be
\item
$\T^AG\in\ct^{\leq0}$ and $\Hom(\T^{-A}G,\ct^{\leq0})=0$.
\item
Every object $F\in\ct^{\leq 0}$ admits a triangle
$E\la F\la D$ with $E\in\ogenu G{}{-A,A}$ and
$D\in\ct^{\leq-1}$.
\setcounter{enumiv}{\value{enumi}}
\ee
Then the \tstr\ $(\ct^{\leq0},\ct^{\geq0})$ is in
the preferred equivalence class.
\epro

\prf
By (i) we have $\T^AG\in\ct^{\leq0}$, hence $\T^mG\in\ct^{\leq0}$ for all
$m\geq A$. Therefore $\ct^{\leq0}$ contains $\ogenul G{}{-A}=\ct_G^{\leq-A}$.
It remains to show an inclusion in the other direction.

But (\ref{C17.3.1}) constructed, for every object $F\in\ct^{\leq0}$, a sequence
$E_1\la E_2\la E_3\la\cdots$ with
$E_m\in\ogenu G{}{1-m-A,A}\subset\ogenul G{}A$.
In Lemma~\ref{L17.5} we proved that $F$ is isomorphic
to $\hoco E_m$. There exists a triangle
\[\xymatrix@C+40pt{
  \ds\bigoplus_{m=1}^\infty E_m \ar[r] &
  F \ar[r] &
\ds\T\left[\bigoplus_{m=1}^\infty E_m\right] 
}\]
where the outside terms obviously lie in $\ogenul G{}A=\ct_G^{\leq A}$. Hence
$F\in\ct_G^{\leq A}$, and since $F\in\ct^{\leq0}$ is arbitrary
we conclude that $\ct^{\leq0}\subset\ct_G^{\leq A}$.
\eprf

\lem{L17.9}
Let $\ct$ be a compactly generated triangulated category, let
$G$ be a compact generator, and let $(\ct_1^{\leq0},\ct_1^{\geq0})$
and $(\ct_2^{\leq0},\ct_2^{\geq0})$ be two equivalent {\it t}--structures.
Let $A>0$ be an integer so that, with $k=1$, the conditions
\be
\item
$\T^AG\in\ct_k^{\leq0}$ and $\Hom(\T^{-A}G,\ct_k^{\leq0})=0$.
\item
Every object $F\in\ct_k^{\leq 0}$ admits a triangle
$E\la F\la D$ with $E\in\ogenu G{}{-A,A}$ and
$D\in\ct_k^{\leq-1}$.
\setcounter{enumiv}{\value{enumi}}
\ee
both hold. Then, after increasing the integer $A$ if
necessary, 
(i) and (ii) will also hold for $k=2$.
Furthermore if (iii) below holds for $k=1$
\be
\setcounter{enumi}{\value{enumiv}}
\item
In the triangle 
$E\la F\la D$ of (ii) above we may  
strengthen the condition on $E$, we may assume
$E\in\ogenu GA{-A,A}\subset\ogenu G{}{-A,A}$.
\setcounter{enumiv}{\value{enumi}}
\ee
then the integer $A$ may be chosen large enough so that
(iii) will hold for $k=2$.
\elem

\prf
Because the {\it t}--structures are equivalent we may choose an
integer $B$ so that $\ct_2^{\leq-B}\subset\ct_1^{\leq0}\subset\ct_2^{\leq B}$.
Hence
$\Hom\big(\T^{-A-B}G,\ct_2^{\leq0}\big)\cong\Hom\big(\T^{-A}G,\ct_2^{\leq-B}\big)=0$,
where the vanishing is because $\ct_2^{\leq-B}\subset\ct_1^{\leq0}$ and
$\Hom\big(\T^{-A}G,\ct_1^{\leq0}\big)=0$. Also
$\T^AG\in\ct_1^{\leq0}\subset\ct_2^{\leq B}$ implies
$\T^{A+B}G\in\ct_2^{\leq 0}$.
This proves (i) for $k=2$, as long as we replace $A$ by $A+B$.

If  $F$ is an object in $\ct_2^{\leq0}\subset\ct_1^{\leq B}$ we may,
using (ii) in combination with
(\ref{C17.3.1}) applied to $\T^BF\in\ct_1^{\leq 0}$, construct a triangle
$E_{2B+1}^{}\la F\la D_{2B+1}^{}$ with
$E_{2B+1}^{}\in\ogenu G{}{-B-A,B+A}$ and
$D_{2B+1}^{}\in\ct_1^{\leq-B-1}\subset\ct_2^{\leq-1}$.
Thus (ii) also holds for $k=2$, as long as $A$
is replaced by $A+B$.

It remains to prove the assertion (iii) for $k=2$, assuming it holds
for $k=1$. By (\ref{C17.3.2}) applied to $\T^BF\in\ct_1^{\leq 0}$, we
may construct the triangle
$E_{2B+1}^{}\la F\la D_{2B+1}^{}$ with
$E_{2B+1}^{}\in\ogenu G{(2B+1)A}{-B-A,B+A}$ and
$D_{2B+1}^{}\in\ct_1^{\leq-B-1}\subset\ct_2^{\leq-1}$.
Thus assertion (iii) holds, but we must replace $A$ by
$\wt A=\max\big[A+B,A(2B+1)\big]$.
\eprf

\pro{P17.11}
Suppose $\ct$ is a weakly approximable triangulated category, $H$ is a
compact generator, and $(\ct_1^{\leq0},\ct_1^{\geq0})$ is any \tstr\ in
the preferred equivalence class. Then there exists an integer
$A>0$ so that
\be
\item
$\T^AH\in\ct_1^{\leq0}$ and $\Hom(\T^{-A}H,\ct_1^{\leq0})=0$.
\item
Every object $F\in\ct_1^{\leq 0}$ admits a triangle
$E\la F\la D$ with $E\in\ogenu H{}{-A,A}$ and
$D\in\ct_1^{\leq-1}$.
\setcounter{enumiv}{\value{enumi}}
\ee
If the category $\ct$ is approximable then the integer $A$
may be chosen to further satisfy
\be
\setcounter{enumi}{\value{enumiv}}
\item
In the triangle 
$E\la F\la D$ of (ii) above we may  
strengthen the condition on $E$, we may assume
$E\in\ogenu HA{-A,A}\subset\ogenu H{}{-A,A}$.
\setcounter{enumiv}{\value{enumi}}
\ee
\epro

\prf
The definition of weakly approximable categories gives us
a compact generator $G$, a \tstr\ $(\ct^{\leq0},\ct^{\geq0})$
and an integer $A$
satisfying (i) and (ii), plus (iii) if $\ct$ is
approximable. Proposition~\ref{P17.7} guarantees
that $(\ct^{\leq0},\ct^{\geq0})$ is in the preferred
equivalence class of
{\it t}--structures. By assumption so is
$(\ct_1^{\leq0},\ct_1^{\geq0})$, hence
the {\it t}--structures $(\ct^{\leq0},\ct^{\geq0})$
and
$(\ct_1^{\leq0},\ct_1^{\geq0})$
are equivalent. By Lemma~\ref{L17.9} we can, by
modifying the integer $A$, also have the conditions
(i), (ii) and [when appropriate] (iii) hold for
the the \tstr\ $(\ct_1^{\leq0},\ct_1^{\geq0})$ and the
compact generator $G$. Thus
we may assume that the {\it t}--structures are the
same. We have a single \tstr\
$(\ct^{\leq0},\ct^{\geq0})=(\ct_1^{\leq0},\ct_1^{\geq0})$,
and two compact generators $G$ and $H$. There exists an integer
$A$ that works for $G$ and the \tstr\  $(\ct^{\leq0},\ct^{\geq0})$,
and we need to produce an integer that works for
$H$ and the \tstr\ $(\ct^{\leq0},\ct^{\geq0})$.

We are given that $G$ and $H$ are compact generators of $\ct$,
hence
$\gen G{}=
\ct^c =\gen H{}$,
and Lemma~\ref{L0.-9238}(ii) allows us to choose an integer
$B>0$ with
$G\in \genu HB{-B,B}$ and 
$H\in \genu GB{-B,B}$.
By (i) for $G$ we know that $\T^AG\in\ct^{\leq0}$ and
$\Hom\big(\T^{-A}G,\ct^{\leq0}\big)=0$. It immediately follows
that $\genu GB{-A-2B,-A}\subset\ct^{\leq0}$ and
that $\Hom\Big(\genu GB{A,A+2B},\ct^{\leq0}\Big)=0$,
and as $\T^{A+B}H\in\genu GB{-A-2B,-A}$ and
$\T^{-A-B}H\in\genu GB{A,A+2B}$ we deduce
that $\T^{A+B}H\in\ct^{\leq0}$ and that
$\Hom\big(\T^{-A-B}H,\ct^{\leq0}\big)=0$. This established (i) for $H$,
if we replace $A$ by $A+B$.

Now for (ii) and (iii): for any $F\in\ct^{\leq0}$ we
know that there exists a triangle $E\la F\la D$ with
$D\in\ct^{\leq-1}$, with $E\in\ogenu G{}{-A,A}$, and if
$\ct$ is approximable we may even choose
$E$ to lie in $\ogenu GA{-A,A}$. But $G$ belongs
to $\genu HB{-B,B}$, and therefore
\[
\ogenu G{}{-A,A}\subset\ogenu H{}{-A-B,A+B}\qquad\text{while}
\qquad \ogenu G{A}{-A,A}\subset\ogenu H{AB}{-A-B,A+B}\ .
\]
Thus (ii) and [when appropriate] (iii) hold for $H$ if $A$
is replaced by $\max(A+B,AB)$.
\eprf

\rmk{R17.13}
We have so far proved Facts~\ref{F0.11}: Proposition~\ref{P17.7}
amounts to \ref{F0.11}(i) and Proposition~\ref{P17.11}
to \ref{F0.11}(ii). The remainder of the section will be
devoted to the basic properties of the subcategory
$\ct^-_c$ of Definition~\ref{D0.13}.
\ermk

\lem{L17.17}
Suppose $\ct$ is a triangulated category with coproducts and
let $(\ct^{\leq0},\ct^{\geq0})$ be a \tstr.
Assume there
exists a compact generator $G$ and an integer $A>0$ so that
$\Hom\big(\T^{-A}G,\ct^{\leq0}\big)=0$.

Then for any compact object $H\in\ct$ there exists an integer
$B>0$, depending on $H$, with
$\Hom\big(\T^{-B}H,\ct^{\leq0}\big)=0$.
\elem

\prf
Let $H\in\ct$ be a compact object.
The fact that $G$ is a compact generator gives
the equality in $H\in\ct^c=\gen G{}$;
Lemma~\ref{L0.-9238}(i) allows us to deduce that
$H\in\genu G{}{-C,C}$ for some $C>0$.
Thus $\T^{-A-C}H\in\genu G{}{A,A+2C}$, and as
$\Hom\Big(\genu
G{}{A,A+2C}\,,\,\ct^{\leq0}\Big)=0$
the Lemma follows, with $B=A+C$.
\eprf

\lem{L17.19}
Suppose $\ct$ is a compactly generated triangulated category
 and
let $(\ct^{\leq0},\ct^{\geq0})$ be a \tstr.
Assume there
exists a compact generator $G$ and an integer $A>0$ so that
$\Hom\big(\T^{-A}G,\ct^{\leq0}\big)=0$.

Then the subcategory $\ct^-_c\subset\ct$ is triangulated.
\elem

\prf
It is clear that $\ct^-_c$ is closed under all suspensions and is
additive. We must show that, if $R\la S\la T\la\T R$ is a triangle so
that $R$ and $T$ belong to $\ct^-_c$, then $S$ must also belong to
$\ct^-_c$.

Choose any integer $m>0$. Because $T$ belongs to $\ct^-_c$ we may
choose a triangle $T'\la T\la T''$ with $T'\in\ct^c$ and
$T''\in\ct^{\leq-m}$. Since $T'$ is compact, Lemma~\ref{L17.17} says
that
we may choose an integer $B>0$ with $\Hom\big(T',\ct^{\leq-B}\big)=0$. 

Now $R$ belongs to $\ct^-_c$, allowing us to choose a triangle $R'\la
R\la R''$ with 
$R'\in\ct^c$ and $R''\in\ct^{\leq-m-B}$. We have a diagram
\[\xymatrix{
T'\ar[r] & T\ar[r] \ar[d] & T'' \\
\T R' \ar[r] & \T R\ar[r]  &\T R''
}\]
The composite from top left to bottom right is a map $T'\la \T R''$,
with $\T R''\in\ct^{\leq-m-B-1}\subset\ct^{\leq-B}$. Since
$B>0$ was chosen so that $\Hom\big(T',\ct^{\leq-B}\big)=0$ the map
$T'\la\T R''$ must vanish, hence the composite $T'\la \T R$ must factor
through $\T R'\la\T R$. We produce a commutative
square
\[\xymatrix{
T'\ar[r]\ar[d] & T \ar[d] \\
\T R' \ar[r] & \T R
}\]
which we may complete to a $3\times3$ diagram where the rows and
columns are triangles
\[\xymatrix{
R' \ar[r] \ar[d] & R\ar[r]  \ar[d] & R''\ar[r] \ar[d] &\T R' \ar[d] \\
S' \ar[r] \ar[d] & S\ar[r]  \ar[d] & S''\ar[r] \ar[d] &\T R' \ar[d] \\
T'\ar[r] \ar[d] & T\ar[r] \ar[d] & T'' \ar[r]\ar[d]  &\T T' \ar[d] \\
\T R' \ar[r] & \T R\ar[r]  &\T R''\ar[r]  &\T^2 R''
}\]
Because $R'$ and $T'$ are compact, the triangle $R'\la S'\la T'$
tells us that $S'$ must be compact. Also
$T''\in\ct^{\leq-m}$ and $R''\in\ct^{\leq-m-B}\subset\ct^{\leq-m}$,
and the triangle $R''\la S''\la T''$ implies that
$S''\in\ct^{\leq-m}$.
The triangle $S'\la S\la S''$ now does the job for $S$.
\eprf

\pro{P17.21}
Suppose $\ct$ is a compactly generated triangulated category
 and
let $(\ct^{\leq0},\ct^{\geq0})$ be a \tstr.
Assume there
exists a compact generator $G$ and an integer $A>0$ so that
$\T^AG\in\ct^{\leq0}$ and $\Hom\big(\T^{-A}G,\ct^{\leq0}\big)=0$.

Then the subcategory $\ct^-_c$ is thick.
\epro

\prf
We already know that $\ct^-_c$
is triangulated, we need to prove it
closed under direct summands. Suppose
therefore that $S\oplus T$ belongs to $\ct^-_c$, we must prove that
so does $S$.

Consider the map $0\oplus\id:S\oplus T\la S\oplus T$. Complete to
a triangle
\[\xymatrix@C+40pt{
S\oplus T\ar[r]^-{0\oplus\id} & S\oplus T \ar[r] & S\oplus\T S
}\]
By Lemma~\ref{L17.19} we deduce that $S\oplus\T S$ belongs to $\ct^-_c$.
Induction on $n$ allows us to prove that, for any $n\geq0$,
the object $S\oplus\T^{2n+1} S$
belongs to $\ct^-_c$. To spell it out: we
have proved the case $n=0$ above.
For any $n$ we know that $\T^{2n+1}(S\oplus\T S)\cong\T^{2n+2}S\oplus\T^{2n+1}S$
lies in $\ct^-_c$, and induction on $n$ allows us to assume that so
does $S\oplus \T^{2n+1}S$.
The triangle
\[\xymatrix@C+40pt{
  \T^{2n+2}S\oplus \T^{2n+1}S \ar[r]^-{0\oplus\id} & S\oplus \T^{2n+1}S \ar[r] &
  S\oplus\T^{2n+3} S
}\]
then informs us that $S\oplus\T^{2n+3} S$ belongs to $\ct^-_c$.

By Remark~\ref{R0.13.99} the category $\ct^-_c$ is contained in $\ct^-$,
and the object $S\oplus\T S$ must belong to $\ct^{\leq\ell} $ for some $\ell>0$.
Hence $S$ belongs to $\ct^{\leq\ell}$ and, for every  integer
$m>0$, we have that $\T^{\ell+m}S\in\ct^{\leq-m}$. Choose
an integer $n\geq0$ with $2n+2\geq\ell+m$; then
$\T^{2n+2}S\in\ct^{\leq-m}$.
Since the object $S\oplus\T^{2n+1}S$ belongs to $\ct^-_c$ we may choose
a triangle $K\la S\oplus\T^{2n+1}S\la P$ with
$K\in\ct^c$ and $P\in\ct^{\leq-m}$. Now form the octahedron on
the composable morphisms $K\la S\oplus\T^{2n+1}S\la S$. We
obtain
\[\xymatrix@C+40pt{
K \ar[r]\ar@{=}[d] & S\oplus\T^{2n+1}S \ar[r]\ar[d] & P\ar[d]\\
K \ar[r] &           S\ar[r]\ar[d] &           Q\ar[d] \\      
         &           \T^{2n+2}S \ar@{=}[r] &  \T^{2n+2}S
}\]
The triangle $P\la Q\la  \T^{2n+2}S$, together with
the fact that both $P$ and $ \T^{2n+2}S$ belong to $\ct^{\leq-m}$,
tell us that $Q$ must belong to $\ct^{\leq-m}$. Now
the triangle $K\la S\la Q$ does the trick for
$S$.
\eprf

The next few results work out how $\ct^-_c$ behaves when $\ct$
is approximable or weakly approximable.

\lem{L17.23}
Let us fix a weakly approximable [or approximable]
triangulated category $\ct$. Choose a compact
generator $G$ and a \tstr\ $(\ct^{\leq0},\ct^{\geq0})$
in the preferred equivalence class. Choose an
integer $A$ as in Proposition~\ref{P17.11},
and
let $m>2A+1$ be an integer.

Then for any $K\in\ct^c\cap\ct^{\leq0}$
there exists an object $L$ and a triangle
$E\la K\oplus L\la D$ with $E\in\genu G{}{1-m-A,A}$ and
$D\in\genul G{}{-m+A}$.

If $\ct$ is approximable we may further assume $E\in\genu G{mA}{1-m-A,A}$.
\elem

\prf
Because $K$ belongs to $\ct^{\leq0}$ the result (\ref{C17.3.1}) permits us
to construct, for every integer $m>0$, a triangle $E_m\la K\la D_m$ with
$E_m\in\ogenu G{}{1-m-A,A}$ and
$D_m\in\ct^{\leq-m}$. If the category $\ct$ is approximable
we may assume $E_m\in\ogenu G{mA}{1-m-A,A}$.

The object $K$ is assumed compact and Lemma~\ref{L17.17} produces
for us a positive integer $B$, which we
may assume $\geq m+2A$, with $\Hom\big(K,\ct^{\leq-B}\big)=0$.
Since we chose $B\geq m+2A$ we have 
$\Hom\Big(\T\ogenu G{}{1-m-A,A}\,,\,\ct^{\leq-B}\Big)=0$,
and in particular 
$\Hom\big(\T E_m,\ct^{\leq -B}\big)=0$. From the triangle
$K\la D_m\la \T E_m$ and the fact that $\Hom\big(K,\ct^{\leq -B}\big)$
and $\Hom\big(\T E_m,\ct^{\leq -B}\big)$ both vanish we deduce
that $\Hom\big(D_m,\ct^{\leq-B}\big)=0$.

From $D_m\in\ct^{\leq -m}$ we construct a triangle
$E'\la D_m\la Q$ with $E'\in\ogenu G{}{1-B-A,-m+A}$
and $Q\in\ct^{\leq-B}$. Since $\Hom\big(D_m,\ct^{\leq-B}\big)=0$ the map
$D_m\la Q$ must vanish, hence
the map $E'\la D_m$ must be a split epimorphism.
Since $E'$ belongs to $\ogenu G{}{1-B-A,-m+A}\subset\ogenul G{}{-m+A}$ so
does its direct summand $D_m$.

We have learned that $K$ belongs to
$\ogenu G{}{1-m-A,A}*\ogenul G{}{-m+A}$,
and if $\ct$ is approximable $K$ even
belongs to the smaller $\ogenu G{mA}{1-m-A,A}*\ogenul G{}{-m+A}$.
Now set
\begin{eqnarray*}
  \cx_1^{}&=&\Coprod\big(G[1-m-A,A]\big)\ ,\\
  \cx_2^{}&=&\Coprod_{mA}^{}\big(G[1-m-A,A]\big)\ ,\\
  \cz&=&\Coprod\big(G(-\infty,-m+A]\big)\ .
\end{eqnarray*}
Then $\cz=\Smr(\cz)$ is closed under direct summands so
$\cz=\ogenul G{}{-m+A}$, while
\[
\ogenu G{}{1-m-A,A}=\Smr(\cx_1^{})\quad\text{ and }\quad\ogenu G{mA}{1-m-A,A}=\Smr(\cx_2^{})\ .
\]
We are given that $K$ belongs to $\Smr(\cx_i^{})*\cz\subset\Smr(\cx_i^{}*\cz)$
with $i=1$ or 2, depending
on whether $\ct$ is approximable.
Choose an object $K'$ in one of
the categories $\cx_i^{}*\cz$ above, so that $K$ is a direct summand
and $K\stackrel f\la K'\la K$ is a pair of morphisms
composing to the identity. Now put
\begin{eqnarray*}
  \ca_1^{}&=&\COprod\big(G[1-m-A,A]\big)\ ,\\
  \ca_2^{}&=&\COprod_{mA}^{}\big(G[1-m-A,A]\big)\ ,\\
  \cc&=&\COprod\big(G(-\infty,-m+A]\big)\ .
\end{eqnarray*}
By \cite[Lemma~1.8]{Neeman17} any morphism from an object in $\ct^c$, to
any of $\cx_1^{}$, $\cx_2$ or $\cz$, factors (respectively)
through an object in $\ca_1^{}$, $\ca_2^{}$ or $\cc$.
Furthermore: we have that $\ct^c\subset\ct$ is a triangulated
subcategory
and contains $\cc$.
Now the map $f:K\la K'$ is a morphism from
$K\in\ct^c$ to an object $K'\in\cx_i^{}*\cz$, with
$i=1$ or $i=2$.  By \cite[Lemma~1.6 and Remark~1.7]{Neeman17}
the map $f$ factors as $K\la K''\la K'$ where $K''\in\ca_i^{}*\cc$,
with $i=1$ or $2$. Since the composite
$K\la K''\la K'\la K$ is the identity we deduce that
$K$ is a direct summand of the object
in $K''\in\ca_i^{}*\cc$, proving the Lemma.
\eprf

\lem{L17.25}
Let us fix a weakly approximable [or approximable]
triangulated category $\ct$. Choose a compact
generator $G$ and a \tstr\ $(\ct^{\leq0},\ct^{\geq0})$
in the preferred equivalence class. 

There exists an integer $B>0$ so that,
for any object $K\in\ct^c\cap\ct^{\leq0}$,
there exists  a triangle
$E\la K\la D$ with $E\in\genu G{}{-B,B}$ and
$D\in\ct^{\leq-1}$.

If $\ct$ is approximable we may further assume $E\in\genu G{B}{-B,B}$.
\elem

\prf
Choose an
integer $A$ as in Proposition~\ref{P17.11}.
We apply Lemma~\ref{L17.23} to the object $K$, with $m=4A+1$,
and obtain a triangle
$E\la K\oplus L\la D$ with $E\in\genu G{}{-5A,A}$ and
$D\in\genul G{}{-3A-1}$; if the category is
approximable we may even assume $E\in\genu G{(4A+1)A}{-5A,A}$.
Now $\genul G{}{-3A-1}\subset\ct^c$ hence the object
$D$ is compact, and it belongs to $\ct^{\leq-2A-1}$
since $G\in\ct^{\leq A}$. Applying
Lemma~\ref{L17.23} to the object $\T^{-2A-1}D$ and with $m=6A$
we obtain a triangle $E'\la D\oplus M\la D'$ with
$E'\in\genu G{}{-9A,-A-1}$ and
$D'\in\genul G{}{-7A-1}$; if the category is
approximable we may even assume $E'\in\genu G{6A^2}{-9A,-A-1}$.
Now complete the composable 
maps
$K\oplus L\oplus M\la D\oplus M\la D'$ to an
octahedron
\[\xymatrix@C+40pt{
E\ar@{=}[d]\ar[r] & E''\ar[r]\ar[d] & E'\ar[d]\\
E\ar[r] & K\oplus L\oplus M\ar[r]\ar[d] & D\oplus M\ar[d] \\
        & D'\ar@{=}[r] & D'
}\]
We know that $E\in\genu G{}{-5A,A}$ and $E'\in\genu G{}{-9A,-A-1}$,
and the triangle $E\la E''\la E'$ tells us that
$E''$ belongs to
\[
\genu G{}{-5A,A} *\genu G{}{-9A,-A-1}\sub
\genu G{}{-9A,A}\ ;
\]
if $\ct$ is approximable  $E''$ belongs to $\genu G{10A^2+A}{-9A,A}$.

Now
the object $D'$ belongs to $\genul G{}{-7A-1}\subset\ct^{\leq-6A-1}$.
The object $E'$ belongs to $\genu G{}{-9A,-A-1}\subset\ct^{\leq-1}$
and the triangle $E'\la D\oplus M\la D'$ guarantees that $D\oplus M$
and therefore its direct summand $M$ belongs to $\ct^{\leq-1}$.
Summarizing we have
\be
\item
  The object $E$ belongs to $\genu G{}{-5A,A}$,
  the object $E''$ belongs to $\genu G{}{-9A,A}$,
  the object $D$ belongs to $\genul G{}{-3A-1}\subset\ct^{-2A-1}$,
  the object $M$ belongs to $\ct^{\leq-1}$ and
  the object $D'$ belongs to $\ct^{\leq-6A-1}$.
\item
  If the category $\ct$ is approximable then the objects $E$ and $E''$
  were chosen so that $E\in\genu G{4A^2+A}{-5A,A}$
  and  $E''\in\genu G{10A^2+A}{-9A,A}$.
\ee
Now consider the following diagram
\[\xymatrix@C+40pt{
  E\ar[r] & K\oplus L\ar[r] \ar[d] & D\\
  E''\ar[r] & K\oplus L\oplus M\ar[r] & D'
}\]
where the vertical map is the direct sum of $\id:L\la L$ with the zero map.
The composite from top left to bottom right is a morphism
$E\la D'$, with $E\in\genu G{}{-5A,A}$ and
$D'\in\ct^{\leq-6A-1}$, hence must vanish. Therefore the composite
$E\la K\oplus L\la K\oplus L\oplus M$ must factor through
$E''\la K\oplus L\oplus M$.
We deduce a commutative square
\[\xymatrix@C+40pt{
  E\ar[r]\ar[d] & K\oplus L \ar[d] \\
  E''\ar[r] & K\oplus L\oplus M
}\]
which we may complete to a $3\times3$ diagram whose rows and columns
are triangles
\[\xymatrix@C+40pt{
  E\ar[r]\ar[d] & K\oplus L\ar[r] \ar[d] & D\ar[d]\\
  E''\ar[r]\ar[d] & K\oplus L\oplus M\ar[r]\ar[d] & D'\ar[d]\\
  \wt E\ar[r] & K\oplus \T K\oplus M\ar[r] & D''
}\]
The triangle $D'\la D''\la \T D$, together with the fact that
$D\in\ct^{-2A-1}$ and $D'\in\ct^{-6A-1}$, tell us that $D''\in\ct^{\leq-2A-2}$.
The triangle $E''\la \wt E\la \T E$, combined with
the fact that $E\in\genu G{}{-5A,A}$
and  $E''\in\genu G{}{-9A,A}$,
tell us that $\wt E\in\genu G{}{-9A,A}$; if $\ct$ is approximable
we have that $E\in\genu G{4A^2+A}{-5A,A}$
and  $E''\in\genu G{10A^2+A}{-9A,A}$ and therefore $\wt E$
belongs to $\genu G{14A^2+2A}{-9A,A}$. Now complete the composable
maps $\wt E\la K\oplus \T K\oplus M\la K$ to an octahedron
\[\xymatrix@C+40pt{
  \wt E\ar[r]\ar@{=}[d] & K\oplus \T K\oplus M\ar[r]\ar[d] & D''\ar[d]\\
  \wt E\ar[r]           & K                \ar[r]\ar[d]   & \wt D\ar[d] \\
    & \T^2K\oplus\T M\ar@{=}[r] &\T^2 K\oplus\T M
}\]
We have that $\T^2K$ and $\T M$ both belong to $\ct^{\leq-2}$ and
$D''$ belongs to $\ct^{\leq-2A-2}$. Hence $\wt D\in\ct^{\leq-2}$, and
the triangle $\wt E\la K\la \wt D$ satisfies the assertion of
the Lemma.
\eprf

\pro{P17.27}
Let us fix a weakly approximable [or approximable]
triangulated category $\ct$. Choose a compact
generator $G$ and a \tstr\ $(\ct^{\leq0},\ct^{\geq0})$
in the preferred equivalence class. Choose an integer $B>0$
as in Lemma~\ref{L17.25}.

Then 
for any object $F\in\ct^-_c\cap\ct^{\leq0}$
there exists  a triangle
$E\la F\la D$ with $E\in\genu G{}{-B,B}$ and
$D\in\ct^{\leq-1}$.

If $\ct$ is approximable we may further assume $E\in\genu G{B}{-B,B}$.
\epro

\prf
Because $F$ belongs to $\ct^-_c$ we may choose a triangle
$K\la F\la D_1$ with $K\in\ct^c$ and $D_1\in\ct^{\leq-1}$. The
triangle $\Tm D_1\la K\la F$, coupled with the fact that
both $\Tm D_1$ and $F$ belong to $\ct^{\leq0}$, tell us that
$K\in\ct^{\leq0}$. Thus $K\in\ct^c\cap\ct^{\leq0}$.

We may therefore apply Lemma~\ref{L17.25}; there exists a triangle
$E\la K\la D_2$ with $E\in\genu G{}{-B,B}$ and $D_2\in\ct^{\leq-1}$.
If $\ct$ is approximable the object $E$ may be chosen
in $\genu GB{-B,B}$.
Now complete the composable maps $E\la K\la F$ to an octahedron
\[\xymatrix{
  E \ar[r]\ar@{=}[d] & K \ar[r] \ar[d] & D_2\ar[d]\\
  E\ar[r] & F \ar[d]\ar[r] & D\ar[d] \\
     & D_1\ar@{=}[r] & D_1
}\]
The triangle $D_2\la D\la D_1$, coupled with the fact
that $D_2$ and $D_1$ both
lie in $\ct^{\leq-1}$, tell us that $D\in\ct^{\leq-1}$.
And
the triangle $E\la F\la D$ satisfies the assertion of the Proposition.
\eprf

\cor{C17.29}
Let $\ct$ be a weakly approximable triangulated category.
Let $G$ be a compact generator and let $(\ct^{\leq0},\ct^{\geq0})$ be
a \tstr\ in the preferred equivalence class.
Choose an integer $B>0$
as in Lemma~\ref{L17.25}.

For any object
$F\in\ct^-_c\cap\ct^{\leq0}$
there exists a sequence of objects $E_1\la E_2\la E_3\la\cdots$
mapping to $F$, and so that $E_m\in\genu G{}{1-m-B,B}$ and
in each triangle $E_m\la F\la D_m$ we have $D_m\in\ct^{\leq-m}$.
For
any such sequence the
non-canonical map $\hoco E_m\la F$ is an isomorphism.

If the category is approximable we may construct the $E_m$ to lie
in $\genu G{mB}{1-m-B,B}$.
\ecor

\prf
The fact that any such sequence would deliver a non-canonical
isomorphism $\hoco E_m\la F$ is contained in Lemma~\ref{L17.5}.
We need to prove the existence of the sequence.

In Proposition~\ref{P17.27} we constructed a triangle $E_1\la F\la D_1$
with $E_1\in\genu G{}{-B,B}$ and $D_1\in\ct^{\leq-1}$.
But $\genu G{}{-B,B}\subset\ct^c$, and in
Remark~\ref{R0.14} we noted that $\ct^c\subset\ct^-_c$.
In the triangle $E_1\la F\la D_1$ we have that both $E_1$ and $F$
lie in $\ct^-_c$, while Lemma~\ref{L17.19} proved that the
category $\ct^-_c$ is triangulated. Therefore $D_1\in\ct^-_c\cap\ct^{\leq-1}$.
If we let $\cs=\ct^-_c$ we are
in the situation of Corollary~\ref{C17.3}, more
specifically the hypotheses of (\ref{C17.3.3}) hold; if $\ct$
is approximable the hypotheses of (\ref{C17.3.4}) hold.
The current
Corollary
is simply the conclusions
of (\ref{C17.3.3}) and (\ref{C17.3.4}).
\eprf

\section{Products and homotopy inverse limits in weakly approximable
  triangulated categories}
\label{SS9999}

This short
section assembles together a few facts about products
an homotopy inverse limits in weakly approximable
triangulated categories. We are ready to prove these
facts now, but will not use them until much later
in the article, in the proof of Theorem~\ref{T1.-1}(ii).
See Section~\ref{SS2}.

It should perhaps be highlighted that 
\cite[Lemmas~3.6  and 3.9(iv)]{Burke-Neeman-Pauwels18}
play a key
role in the proofs of the present section.
It was thus the joint work, with
Burke and Pauwels, that produced the key lemmas
allowing me to prove Theorem~\ref{T1.-1}(ii).
See also Remark~\ref{R0.987656789}.

\lem{LL2.1}
Suppose $\ct$ is a weakly approximable triangulated category, and let
$\big(\ct^{\leq0},\ct^{\geq0}\big)$ be a \tstr\ in
the preferred equivalence class. Suppose
$
\cdots\la Z_3\la Z_{2}\la Z_1\la Z_0
$
is an inverse sequence of objects in $\ct$ so that, for any integer
$n$, the functor $(-)^{\geq n}$ takes it to a sequence that is eventually
stable.

Let $Z=\holim Z_i$ and consider the natural map $f_i:Z\la Z_i$. Then
for $i\gg0$ the functor $(-)^{\geq n}$ takes the map $f_i$ to an isomorphism.

More precisely: there exists an integer $L>0$ so that, whenever the
sequence
$
\cdots\la Z_3^{\geq-L}\la Z_{2}^{\geq-L}\la Z_1^{\geq-L}\la Z_0^{\geq-L}
$
is constant, the map $f_i^{\geq0}:Z^{\geq0}\la Z_i^{\geq0}$ is an
isomorphism.
\elem

\prf
Let $G$ be a compact generator for $\ct$ and suppose
$\big(\ct^{\leq0},\ct^{\geq0}\big)$ is a \tstr\ in the preferred equivalence
class. By \cite[Lemma~3.9(iv)]{Burke-Neeman-Pauwels18}
there exists an integer $A>0$ so
that, if $\Hom(\T^{i}G,X)=0$ for all $i\leq A$, then $X$ must belong to
$\ct^{\leq0}$. Choose and fix such an integer $A$.
Then choose an integer $B>0$ with $\Hom\big(G,\ct^{\leq-B}\big)=0$.

By shifting and passing to a subsequence
it suffices to prove the ``moreover'' assertion, and we assert
that $L=A+B+1$ works. Suppose therefore that
the sequence
$
\cdots\la Z_3^{\geq-A-B-1}\la Z_{2}^{\geq-A-B-1}\la Z_1^{\geq-A-B-1}\la Z_0^{\geq-A-B-1}
$
is constant, and
consider the commutative diagram in which the rows are triangles
\[\xymatrix@C+30pt{
\ds\prod_{i=0}^\infty  Z_i^{\leq-A-B-2}\ar[r]^-{}
& \ds\prod_{i=0}^\infty Z_i \ar[r]^-{}\ar[d]_-{1-\mathrm{shift}}
& \ds\prod_{i=0}^\infty Z_i^{\geq-A-B-1}
\ar[d]^-{1-\mathrm{shift}}\\
\ds\prod_{i=0}^\infty  Z_i^{\leq-A-B-2}\ar[r]^-{}
& \ds\prod_{i=0}^\infty Z_i \ar[r]^-{}
& \ds\prod_{i=0}^\infty Z_i^{\geq-A-B-1}
}\]
We may complete it to a $3\times3$ diagram of triangles
\[\xymatrix@C+30pt{
X \ar[r]\ar[d] & Z\ar[r]\ar[d]& Y\ar[d] \\ 
\ds\prod_{i=0}^\infty  Z_i^{\leq-A-B-2}\ar[r]^-{}\ar[d]
& \ds\prod_{i=0}^\infty Z_i \ar[r]^-{}\ar[d]_-{1-\mathrm{shift}}
& \ds\prod_{i=0}^\infty Z_i^{\geq-A-B-1}
\ar[d]^-{1-\mathrm{shift}}\\
\ds\prod_{i=0}^\infty  Z_i^{\leq-A-B-2}\ar[r]^-{}
& \ds\prod_{i=0}^\infty Z_i \ar[r]^-{}
& \ds\prod_{i=0}^\infty Z_i^{\geq-A-B-1}
}\]
where the middle column is the definition of $Z=\holim Z_i$.
The inverse sequence $Z_i^{\geq-A-B-1}$ is constant,
hence the object $Y$ is canonically isomorphic
to $Z_i^{\geq-A-B-1}$.
What we know about the object $X$ is that it sits in the triangle
given by the left column
\[\xymatrix@C+30pt{
\ds\prod_{i=0}^\infty(\Tm Z_i)^{\leq-A-B-1}\ar[r] &
X\ar[r] &
\ds\prod_{i=0}^\infty Z_i^{\leq-A-B-2}
}\]
By the choice of $B$ we have that $\Hom(\T^iG,-)$ kills the two outside terms
whenever $i\leq A+1$, hence the choice of $A$ gives that the two outside terms
lie in $\ct^{\leq-1}$. Therefore $X\in\ct^{\leq-1}$, and it
follows that the truncation
$(-)^{\geq0}$ takes the map $Z\la Y=Z_i^{\geq-A-B-1}$ to an isomorphism. 
\eprf

\pro{PP2.95}
Suppose $\ct$ is a weakly approximable triangulated category,
and let $\big(\ct^{\leq0},\ct^{\geq0}\big)$ be a \tstr\
in the preferred equivalence
class. Let
$Y$ be any object in $\ct$,
assume
$
\cdots\la Z_3\la Z_{2}\la Z_1\la Z_0
$
is an inverse sequence of objects in $\ct$, and
let $f_*:Y\la Z_*$ be a map from $Y$ to the inverse system.
Suppose that, for any integer
$i>0$, there exists an integer $N>0$ so that
\[
n\geq N\quad\Longrightarrow\quad f_n^{\geq-i}:Y^{\geq-i}\la Z_n^{\geq-i}\text{ is an
  isomorphism.}
\]
If $Z=\holim Z_n$ is the homotopy inverse limit, then the
(non-canonical) map $f:Y\la Z$  
is an isomorphism.

Applying to the inverse system $Z_n=Y^{\geq-n}$, we learn that $\ct$
is left-complete with respect to any \tstr\ in the preferred
equivalence class.
\epro

\prf
By
\cite[Lemma~3.6]{Burke-Neeman-Pauwels18}
the \tstr\ is nondegenerate, hence it suffices to prove that the morphism
$\ph^{\geq-i}:Y^{\geq-i}\la Z^{\geq-i}$ is an isomorphism for every $i\in\zz$.
By shifting it suffices to prove this for $i=0$.
Let $L>0$ be the integer in the ``moreover'' part of Lemma~\ref{LL2.1},
choose an integer $N$ so large that $f_n^{\geq-L}:Y^{\geq-L}\la Z_n^{\geq-L}$
is an isomorphism
whenever $n\geq N$, and
apply the functor $(-)^{\geq0}$ to the commutative triangle below
\[\xymatrix@C+40pt@R-20pt{
Y\ar[dr]^-{f^{}_N}\ar[dd]_-{f} &    \\
                & Z_N^{\geq-L} \\
Z\ar[ru]_-{\ph_N^{}} &
}\]
Lemma~\ref{LL2.1} teaches us that $\ph_N^{\geq0}$ is an isomorphism.
Since $f_N^{\geq0}$ is also an isomorphism, the
commutativity forces $f^{\geq0}$
to be an isomorphism.
\eprf

\lem{LL29.-1}
Let $\ct$ be a weakly approximable triangulated category,
and choose a \tstr\  $\big(\ct^{\leq0},\ct^{\geq0}\big)$
in the preferred equivalence class.
 Suppose
$
\cdots\la Z_3\la Z_{2}\la Z_1\la Z_0
$
is an inverse sequence of objects in $\ct^-_c$ so that, for any integer
$n$, the functor $(-)^{\geq n}$ takes it to a sequence that is eventually
stable.

Then $Z=\holim Z_i$ belongs to $\ct^-_c$. 
\elem

\prf
Since $\ct$ is weakly approximable, the ``moreover''
part of Lemma~\ref{LL2.1} provides
an integer $L>0$ such that, given an inverse sequence
$\cdots \la Z_3\la Z_2\la Z_1\la Z_0$ in $\ct$ with
$\cdots \la Z_3^{\geq-L}\la Z_2^{\geq-L}\la Z_1^{\geq-L}\la Z_0^{\geq-L}$
all isomorphisms, we have
that the map $\big[\holim Z_*\big]^{\geq0} \la Z_0^{\geq0}$ is an
isomorphism. Also: since $\ct$ is weakly approximable it
has a compact generator $G$,
and Corollary~\ref{C17.29} permits
us to choose an integer $B>0$ such that
\be
\item
$\Hom\big(G,\ct^{\leq-B}\big)=0$.
\item
Every object $X\in\ct^-_c$ admits a triangle $W\la X\la D$, with
$W\in\genuf G{}{-B}$ and $D\in\ct^{\leq0}$.
\ee

OK: we have chosen the integers $L$ and $B$ and it's time to get
to work. Suppose we are given in $\ct^-_c$ a sequence
$\cdots \la Z_3\la Z_2\la Z_1\la Z_0$ satisfying the hypotheses, put
$Z=\holim Z_*$, we need to show that $Z\in\ct^-_c$. 
Choose any integer $m>0$. We
choose an integer $i>0$ so that the maps in the subsequence
$\cdots \la Z_{i+2}^{\geq-m-L-2B+1}\la Z_{i+1}^{\geq-m-L-2B+1}\la Z_i^{\geq-m-L-2B+1}$
are all isomorphisms. By (ii) we may choose a triangle
$W\la Z_i\la D$ with $W\in\genuf G{}{-m-B}$ and $D\in\ct^{\leq-m}$.
By the choice of $L$ the map
$Z^{\geq-m-2B+1}\la Z_i^{\geq-m-2B+1}$ is an isomorphism,
and in the triangle $Z\la Z_i\la \wt D$ we have
$\wt D\in\ct^{\leq-m-2B}$. Therefore the composite
$W\la Z_i\la \wt D$ is a map from $W\in\genuf G{}{-m-B}$ to
$\wt D\in\ct^{\leq-m-2B}$  and vanishes by (i). We may therefore
factor $W\la Z_i$ as $W\la Z\la Z_i$. Completing to an octahedron
\[\xymatrix{
W\ar[r]\ar@{=}[d] & Z\ar[r]\ar[d] & D'\ar[d] \\
W\ar[r]           & Z_i\ar[r]\ar[d] & D\ar[d] \\
                  & \wt D\ar@{=}[r] & \wt D
}\]
produces a triangle $W\la Z\la D'$, with $W\in\ct^c$.
The triangle $\Tm\wt D\la D'\la D$, coupled with the facts
that $\Tm\wt D\in\ct^{\leq-m-2B+1}$ and $D\in\ct^{\leq-m}$,
guarantee that
$D'\in\ct^{\leq-m}$.
\eprf

\section{Examples}
\label{S19}

In Sections~\ref{S17} and \ref{SS9999} we developed
some abstraction, and it's high time
to look at examples. We begin with the trivial ones.

\exm{E19.25}
Let $R$ be an associative (possibly noncommutative) ring, 
and put $\ct=\D(R)$. Note: our 
$\D(R)$ is an abbreviation for what, in
the 
more explicit notation standard in representation theory, 
usually goes by the name $\D(\Mod R)$. 
It is the derived category whose objects are the 
(possibly unbounded) cochain complexes of left $R$-modules. The
category $\ct$ has coproducts and $R$ is a compact generator.
Let $(\ct^{\leq0},\ct^{\geq0})$ be the standard \tstr. Then
$\Hom\big(\Tm R,\ct^{\leq0}\big)=0$ and $\T R\in\ct^{\leq0}$.

Let $F$ be any object in $\ct^{\leq0}$. That is, we take
a cochain complex $F$ with $H^\ell(F)=0$ for all $\ell>0$.
Such a complex has a free resolution; it is
isomorphic in $\ct$ to the cochain complex
\[\xymatrix{
  \cdots\ar[r] & F^{-3}\ar[r] & F^{-2}\ar[r] & F^{-1}\ar[r] & F^{0}\ar[r] &
   0\ar[r] & 0\ar[r] &0\ar[r] &\cdots
}\]
with $F^i$ free $R$--modules. The brutal truncation produces for us
a short exact sequence of cochain complexes
\[\xymatrix{
  \cdots\ar[r] & 0\ar[r]\ar[d] & 0\ar[d]\ar[r]
  & 0\ar[d]\ar[r] & F^{0}\ar@{=}[d]\ar[r] &
   0\ar@{=}[d]\ar[r] & 0\ar@{=}[d]\ar[r] &0\ar[r]\ar@{=}[d] &\cdots \\
  \cdots\ar[r] & F^{-3}\ar[r]\ar@{=}[d] & F^{-2}\ar@{=}[d]\ar[r]
  & F^{-1}\ar@{=}[d]\ar[r] & F^{0}\ar[d]\ar[r] &
   0\ar@{=}[d]\ar[r] & 0\ar@{=}[d]\ar[r] &0\ar[r]\ar@{=}[d] &\cdots \\
  \cdots\ar[r] & F^{-3}\ar[r] & F^{-2}\ar[r] & F^{-1}\ar[r] & 0\ar[r] &
   0\ar[r] & 0\ar[r] &0\ar[r] &\cdots
}\]
and this is a triangle $E\la F\la D$ with $E\in\ogenu R1{0,0}$ and
$D\in\ct^{\leq-1}$. The category is approximable, the standard \tstr\
is in preferred equivalence class, and $\ct^-_c$ is
just $\D^-(\proj R)$,
the category of bounded-above complexes which
admit projective resolutions by finitely generated projectives.

The category $\ct^b_c$ is the intersection $\ct^-_c\cap\ct^b$;
it consists of the objects in $\D^-(\proj R)$
with only finitely many nonzero cohomology groups.
We have an inclusion $\D^b(\proj R)\subset\ct^b_c$,
and for $R$ general I don't know much
about the difference $\ct^b_c-\D^b(\proj R)$.
When $R$ is noetherian we have $\ct^b_c=\D^b(\mod R)$, which
is usually much larger than $\D^b(\proj R)$.
\eexm

\exm{E19.64}
A very similar analysis works when $\ct$ is the homotopy category of spectra.
The sphere $S^0$ is a compact generator.
Consider the \tstr\ where $\ct^{\leq0}$ is the category of connective
spectra---these are the spectra $F$ with $\pi_i(F)=0$ when $i<0$.
Then $\Hom\big(\Tm S^0,\ct^{\leq0}\big)=0$ and $\T S^0\in\ct^{\leq0}$.
And any object $F\in\ct^{\leq0}$ admits a triangle
$E\la F\la D$ with $E\in\ogenu {S^0}1{0,0}$ and $D\in\ct^{\leq-1}$; this
just says that we may choose a bouquet of zero-spheres $E$ and a map
$E\la F$ which is surjective on $\pi_0$. The category is
approximable and the \tstr\ above is in the preferred
equivalence class.
\eexm

\rmk{R19.303}
Examples~\ref{E19.25} and \ref{E19.64} should be viewed as the baby case.
If $\ct$ has a compact generator $G$, such that $\Hom(G,\T^iG)=0$ for
all $i>0$, then $\ct$ is approximable. Just take the \tstr\
$\big(\ct_G^{\leq0},\ct_G^{\geq0})\big)$; then
$\T G\in\ct^{\leq0}_G$ and $\Hom\big(\Tm G,\ct_G^{\leq0}\big)=0$,
and every object $F\in\ct_G^{\leq0}$ admits a triangle
$E\la F\la D$ with $E\in\ogenu G1{0,0}$ and $D\in\ct_G^{\leq-1}$.
\ermk

Now we turn to more complicated examples of approximable
and weakly approximable categories. We begin
with

\rmd{R19.1}
If $X$ is a quasicompact, quasiseparated scheme then $\ct=\Dqc(X)$ has
a single compact generator, see Bondal and Van den
Bergh~\cite[Theorem~3.1.1(ii)]{BondalvandenBergh04}.
Let $G$ be any such compact generator;
\cite[Theorem~3.1.1(i)]{BondalvandenBergh04}
tells us that $G$ is a perfect complex.
\ermd

Since the categories $\Dqc(X)$ have single
compact generators, one can wonder which of
them is approximable or weakly
approximable. To address this question we begin with

\lem{L19.501}
Let $X$ be a separated, quasicompact scheme, let $\ct=\Dqc(X)$ be
its derived category, and let $(\ct^{\leq0},\ct^{\geq0})$ be the
standard \tstr. Then there is a compact generator $G'\in\ct$ and
an integer $A>0$, so that every object $F\in\ct^{\leq0}$ admits
a triangle $E\la F\la D$ with $E\in\ogenu {G'}A{-A,A}$ and $D\in\ct^{\leq-1}$.
\elem

\prf
Absolute noetherian approximation,
that is
Thomason and Trobaugh~\cite[Theorem~C.9]{ThomTro}
or~\cite[Tags 01ZA combined with 01YX]{stacks-project},
allows us to choose a separated
scheme $Y$, of finite type over $\zz$, and an affine
map $f:X\la Y$. 
From \cite[Definition~5.2 and Theorem~5.8]{Neeman17}, in the
special case where $\cs=\Dqc(Y)^{\leq0}$ for
the standard \tstr\ and $m=0$, we learn:
\be
\item
There exists a compact generator $G\in\Dqc(Y)$ and an integer $A>0$, so that
every $F'\in\Dqc(Y)^{\leq0}$ admits a triangle $E'\la F'\la D'$ with
$E'\in\ogenu GA{-A,A}$ and $D'\in\Dqc(Y)^{\leq-1}$.
\setcounter{enumiv}{\value{enumi}}
\ee
Hall and Rydh~\cite[Lemma~8.2]{Hall-Rydh13} tells us that
$G'=\LL f^*G$
is a compact generator for $\Dqc(X)$;
this is our choice of $G'$ for the Lemma.
Now take any object $F\in\Dqc(X)^{\leq0}$. Since
$f:X\la Y$ is affine we have $\R f_*F\cong f_*F\in\Dqc(Y)^{\leq0}$, and (i)
above permits us to find a triangle
$E'\la f_*F\la D'$ with $E'\in\ogenu GA{-A,A}$ and $D'\in\Dqc(Y)^{\leq-1}$.
Applying the functor $\LL f^*$, and remembering that
$\LL f^*\Dqc(Y)^{\leq0}\subset\Dqc(X)^{\leq0}$, we deduce
\be
\setcounter{enumi}{\value{enumiv}}
\item
There is in $\Dqc(X)$ a triangle
$\LL f^*E'\la \LL f^*f_* F\la\LL f^*D'$, with $\LL f^*D'\in\Dqc(X)^{\leq-1}$
and $\LL f^*E'\in\LL f^*\ogenu GA{-A,A}\subset\ogenu{\LL f^*G}A{-A,A}$.
\setcounter{enumiv}{\value{enumi}}
\ee
But the counit of adjunction gives a map $\e:\LL f^*f_*F\la F$,
and the fact that the maps
$f_*F\stackrel{\eta f_*}\la f_*\LL f^* f_*F \stackrel{ f_*\e}\la f_*F$
compose to the identity tells us that the functor
$f_*$ takes $\e:\LL f^*f_*F\la F$ to a split epimorphism.
In particular $f_*\e$ induces an epimorphism on
cohomology sheaves and, because $f$ is affine, this means that
$\e$ induces an epimorphism of cohomology sheaves already over
$X$. We have a morphism
$\e:\LL f^*f_*F\la F$ in $\Dqc(X)^{\leq0}$ and, if we complete
it to a triangle, the long exact sequence of cohomology
sheaves gives
\be
\setcounter{enumi}{\value{enumiv}}
\item
In the triangle
$\LL f^* f_*F \stackrel{ \e}\la F \la D''$ we have $D''\in\Dqc(X)^{\leq-1}$.
\setcounter{enumiv}{\value{enumi}}
\ee
Next we form the octahedron
\[\xymatrix{
\LL f^*E'\ar[d]\ar@{=}[r] &\LL f^*E'\ar[d] & \\
\LL f^*f_*F\ar[d]\ar[r]^-\e & F\ar[d] \ar[r] & D''\ar@{=}[d]\\
\LL f^*D' \ar[r] &D\ar[r] & D''
}\]
and (ii) tells us that $\LL f^*E'\in\ogenu{\LL f^*G}A{-A,A}$ and
$\LL f^*D'\in\Dqc(X)^{\leq-1}$, while (iii) gives that $D''\in\Dqc(X)^{\leq-1}$.
The triangle $\LL f^*D'\la D\la D''$ tells us that $D\in\Dqc(X)^{\leq-1}$,
and the triangle $\LL f^*E'\la F\la D$ does the trick.
\eprf

\exm{E19.2}
Assume $X$ is separated and quasicompact,
and let the \tstr\ on $\ct=\Dqc(X)$ be the standard one.
Lemma~\ref{L19.501} finds a generator $G'$ and
an integer $A>0$ so that, for every object $F\in\ct^{\leq0}$, there exists
a triangle $E\la F\la D$ with $E\in\ogenu {G'}A{-A,A}$ and $D\in\ct^{\leq-1}$.
From Reminder~\ref{R19.1} we know that the compact generator $G'\in\Dqc(X)$ is a perfect complex and we may,
after increasing the
integer $A$ if necessary, guarantee that
$\Hom\big(\T^{-A}G',\ct^{\leq0}\big)=0$ and
$\T^{A}G'\in\ct^{\leq0}$. 
Putting this together we have that $\ct$ is approximable; it
satisfies
Definition~\ref{D0.9} for
the compact generator $G'$, the standard \tstr, and the integer $A>0$.

And now Proposition~\ref{P17.7}
informs us that
the standard \tstr\ is in the preferred
equivalence class.
\eexm

\rmk{R19.9096}
In the recent
\cite[Theorem~3.2~(iii) and (iv)]{Neeman22A}
the reader can find a proof that, in
the generality where $X$ is
only assumed quasicompact and quasiseparated,
the category $\ct=\Dqc(X)$ is weakly approximable and
the standard \tstr\ is in the preferred equivalence
class.

In the present article we care much more about the case
where $\ct$ is not only weakly approximable, but
genuinely approximable. And, at least with the current
machinery, to prove the approximability of
$\Dqc(X)$ we need to assume $X$ separated.
The methods used in
\cite{Neeman22A}, to prove the weak approximability
of $\Dqc(X)$ for general $X$, do not modify
in a straightforward way
to prove approximability. In
\cite[Remarks~8.1 and 8.2]{Neeman22A}
there is a discussion: these techniques prove the
weak approximability of categories that are decidedly
not approximable.
\ermk

\rmk{R19.8976}
Let $X$ be a quasicompact, quasiseparated scheme, and
let $\ct=\Dqc(X)$. The category $\ct$
has a single compact generator
by Reminder~\ref{R19.1}, and we
now know that the standard \tstr\ is in the preferred
equivalence class. One immediately computes that the
preferred $\ct^-$, $\ct^+$ and $\ct^b$ are given by
\[
\ct^-=\Dqcmi(X),\qquad\ct^+=\Dqcpl(X),\qquad
\ct^b=\Dqcb(X).
\]
Next we turn our attention to the computation of
$\ct^-_c$ and its subcategory $\ct^b_c$.

To facilitate the discussion
it might be helpful to imagine
two classes of objects in $\ct=\Dqc(X)$.
\be
\item
  The objects belonging to $\ct^-_c=\Dqc(X)^-_c$, for
  the $\ct^-_c$ defined
  using the standard \tstr.
\item
  The objects locally with this property. That is the objects
  $P\in\Dqc(X)$
  such that, for any open immersion $j:\spec R\la X$, the
  object $\LL j^*P\in\Dqc\big(\spec R\big)\cong\D(R)$ is
  in $\D(R)^-_c$. See Example~\ref{E19.25} for a
  description of $\D(R)^-_c$.
\ee
The objects satisfying (ii) are classically
called \emph{pseudocoherent,} they were first studied
in Illusie's expos\'es~\cite{Illusie71B,Illusie71A} in SGA6.
Now
\cite[Theorem~4.1]{Lipman-Neeman07} is precisely the statement
that the objects satisfying (ii) all satisfy (i).
It is trivial to check that the objects satisfying (i)
must satisfy (ii); this means that, for the
standard \tstr\ on $\ct=\Dqc(X)$,
the subcategory $\ct^-_c$ is just
$\Dqcp(X)\subset\Dqc(X)$,
the subcategory of pseudocoherent complexes.
If $X$ happens to be noetherian then pseudocoherence
simplifies to something more familiar: for noetherian $X$
we have $\ct^-_c=\Dqcp(X)=\dmcoh(X)$. And the subcategory
$\ct^b_c\subset\ct^-_c$ is also explicit, it
is the full subcategory $\Dqcpb(X)\subset\Dqcp(X)$ whose
objects are the pseudocoherent
complexes with bounded cohomology. When
$X$ is noetherian this simplifies to
$\ct^b_c=\Dqcpb(X)=\dcoh(X)$.
\ermk

\smr{S19.1098}
Until now we have shown
\be
\item
  The category $\Dqc(X)$ is
  weakly approximable for all quasicompact, quasiseparated
  $X$, and approximable as long as $X$ is separated.
\item
  The standard \tstr\ is in the preferred equivalence
  class.
\item
  We have computed the intrinsic subcategories $\Dqc(X)^-$,
  $\Dqc(X)^+$, $\Dqc(X)^b$, $\Dqc(X)^-_c$ and $\Dqc(X)^b_c$.
\setcounter{enumiv}{\value{enumi}}
\ee  
The computations in (iii) yielded no big surprises,
the intrinsic subcategories turned out to be classical,
old friends.

When we apply Corollary~\ref{C17.29} we discover something new.
\be
\setcounter{enumi}{\value{enumiv}}
\item
  Let $X$ be a quasicompact, quasiseparated scheme,
  and let $G$ be a compact generator of $\Dqc(X)$.
There exists
 an integer $B>0$ so that, for any
 integer $m>0$ and any
object $F\in\Dqcp(X)\cap\Dqc(X)^{\leq0}$,
there is a triangle $E_m\la F\la D_m$ with
$E_m\in\genu G{}{1-m-B,B}$ and
$D_m\in\Dqc(X)^{\leq-m}$.
\item
  If $X$ is separated we may choose $B$
  to guarantee that
$E_m\in\genu G{mB}{1-m-B,B}$.
\setcounter{enumiv}{\value{enumi}}
\ee
\esmr

\section{An easy representability theorem}
\label{S74}

It's time to start proving representability theorems. The main theorems
in the article are a little technical to prove---they hinge 
on taking homotopy colimits
carefully. For this reason I thought it best to illustrate the
methods in a simple case, which involves no homotopy
colimits. In this section we give a simple proof of an old
result of Rouquier, generalizing an even older result of Bondal and
Van den Bergh.

But first we recall the following little lemma, well-known
in some circles.

\rmd{R74.9834}
Let $\cs$ be a triangulated category, let $\cb\subset\cs$
be a full, replete subcategory with $\T\cb=\cb$, and let
$\ca$ be an abelian category. Assume we are given
a short exact sequence
$0\la H'\la H\la H''\la 0$ of functors $\cb\op\la\ca$.

If two of the three
functors $H',H,H'':\cb\op\la\ca$ are $\cb$--cohomological
as in Definition~\ref{D37.105},
then so is the third.
\ermd

\prf
Let $x\la y\la z\la\T x$ be a triangle in $\cs$, such that
the objects $x,y,z$ all belong to $\cb$. Then with $\wt H$
being any of $H'$, $H$ or $H''$, the sequence
\[\xymatrix{
\cdots\ar[r] &\wt H(\T x)\ar[r] &
\wt H(z)\ar[r] &\wt H(y)\ar[r] &\wt H(x)\ar[r] & 
\wt H(\Tm z)\ar[r]&\cdots
}\]
is a cochain complex in the abelian category $\ca$.
Thus we have three cochain
complexes, one for each of $H'$, $H$ or $H''$. Because
two of the functors are $\cb$--cohomological, two
of these three cochain complexes must be acyclic.
But then the short exact sequence 
$0\la H'\la H\la H''\la 0$ tells us that the third
complex is also acyclic.
\eprf
  
Next we prove

\lem{L74.1}
Let $R$ be a commutative, noetherian ring, let
$\cs$ be an $R$--linear triangulated
category, let $G\in\cs$ be an object and,
with the notation of Definition~\ref{D37.105}
and Remark~\ref{R-2.3}, assume
$\Hom(-,G)$ is a $G$--finite cohomological functor.
Let $k\geq0$ be an integer and
$H$ a finite $\gen G{2^k}$--cohomological functor. Then
there exists an object
$F\in\gen G{2^k}$ and an epimorphism $\Hom(-,F)|_{\gen G{2^k}}\la H(-)$.
\elem

\prf
The proof is by induction on $k$. Suppose $k=0$:
by hypothesis $H(\T^{-i}G)$ is a finite $R$--module, and vanishes for $i$
outside a
bounded interval $[-A,A]$. For each $i$ with $-A\leq i\leq A$ choose
a finite set of generators $\{f_{ij},\,j\in J_i\}$ for
the $R$--module $H(\T^{-i}G)$. Yoneda tells us that
each $f_{ij}$ corresponds to a natural transformation
$\ph_{ij}:\Hom(-,\T^{-i}G)|_{\gen G1}\la H(-)$. Define
$F$ to be $F=\oplus_{i=-A}^A\oplus_{j\in J_i}^{}\T^{-i}G$,
and let $\ph:\Hom(-,F)|_{\gen G1}\la H(-)$ be the
composite
\[\xymatrix@C+30pt{
  \Hom(-,F)|_{\gen G1}\ar@{=}[r] &
  \ds\bigoplus_{i=-A}^A\,\,\bigoplus_{j\in J_i}\Hom(-\,,\,\T^{-i}G)|_{\gen G1}
  \ar[r]^-{(\ph_{ij})} & H(-)
}\]
Obviously $F$ belongs to $\gen G1$ and $\ph$ is surjective.

Now suppose we know the Lemma for $k\geq0$, and let $H$ be a
finite $\gen G{2^{k+1}}$--cohomological functor. Then
the restriction of $H$ to  $\gen G{2^k}$ is a finite
$\gen G{2^k}$--cohomological functor, and induction
permits us to find an object $F_1\in\gen G{2^k}$
and an epimorphism $\ph:\Hom(-,F_1)|_{\gen G{2^k}}\la H(-)|_{\gen G{2^k}}$.
Complete the map $\ph$ to a short exact sequence
\[\xymatrix@C+20pt{
  0\ar[r] & H'(-)\ar[r]^-{\s} &\Hom(-,F_1)|_{\gen G{2^k}} \ar[r]^-{\ph} &
  H(-)|_{\gen G{2^k}}\ar[r] & 0\ .
}\]
By Reminder~\ref{R74.9834} the functor $H'$ is $\gen G{2^k}$--cohomological, it is clearly finite, and induction
applies again to tell us that there exists an object $F_2\in\gen G{2^k}$
and an epimorphism $\rho:\Hom(-,F_2)|_{\gen G{2^k}}\la H'(-)$.
Combining the results we deduce an 
exact sequence of functors
\[\xymatrix@C+20pt{
  \Hom(-,F_2)|_{\gen G{2^k}}  \ar[r]^-{\s\rho} &
  \Hom(-,F_1)|_{\gen G{2^k}} \ar[r]^-{\ph} &
  H(-)|_{\gen G{2^k}}\ar[r] & 0\ .
}\]
Because $F_1$ and $F_2$ both lie in $\gen G{2^k}$
the functors $\Hom(-,F_i)|_{\gen G{2^k}}$ are
representable for $i\in\{1,2\}$---Yoneda's lemma applies.
The natural transformation
$\s\rho:\Hom(-,F_2)|_{\gen G{2^k}}\la \Hom(-,F_1)|_{\gen G{2^k}}$
is $\Hom(-,\alpha)|_{\gen G{2^k}}$ for some morphism $\alpha:F_2\la F_1$,
the natural transformation $\ph:\Hom(-,F_1)|_{\gen G{2^k}}\la H(-)|_{\gen G{2^k}}$
corresponds to some element $y\in H(F_1)$, and the vanishing
of the composite $\ph\s\rho$ says that $H(\alpha):H(F_1)\la H(F_2)$
must take $y$ to zero. 

Complete $\alpha: F_2\la F_1$ to a triangle
$F_2\stackrel\alpha\la F_1\stackrel\beta\la F\stackrel\gamma\la \T F_2$.
As $F_1$ and $\T F_2$ belong to $\gen G{2^k}$, the triangle
tells us that $F$ must be in $\gen G{2^k}*\gen G{2^k}\subset\gen G{2^{k+1}}$.
And now we remember that $H$ is actually a $\gen G{2^{k+1}}$--cohomological
functor. The sequence
$H(F)\stackrel{H(\beta)}\la H(F_1)\stackrel{H(\alpha)}\la H(F_2)$ is
exact, and the vanishing of $H(\alpha)(y)$ says that there exists an
element $x\in H(F)$ so that $H(\beta):H(F)\la H(F_1)$ takes
$x\in H(F)$ to $y\in H(F_1)$. By Yoneda $x$ corresponds to a natural
transformation $\psi:\Hom(-,F)|_{\gen G{2^{k+1}}}\la H(-)$.
The fact that $H(\beta)x=y$ translates, via Yoneda,
to the assertion that the composite
\[\xymatrix@C+20pt{
  \Hom(-,F_1)|_{\gen G{2^{k+1}}}  \ar[rr]^-{\Hom(-,\beta)|_{\gen G{2^{k+1}}}} &&
  \Hom(-,F)|_{\gen G{2^{k+1}}} \ar[r]^-{\psi} &
  H(-)
}\]
restricts to be $\ph$ on the category $\gen G{2^k}$.

We assert that $\psi$ is surjective.
Take any object $C'\in\gen G{2^{k+1}}$; we need to
show the surjectivity of the map $\psi:\Hom(C',F)\la H(C')$.
Now $\gen G{2^{k+1}}=\Smr\big[\gen G{2^{k}}*\gen G{2^{k}}\big]$,
so there exists an object $C''\in\gen G{2^{k+1}}$ with
$C'\oplus C''\in\gen G{2^{k}}*\gen G{2^{k}}$. Put $C=C'\oplus C''$ and
it clearly suffices to prove the surjectivity
of $\psi:\Hom(C,F)\la H(C)$. Choose a triangle 
$A\stackrel{\alpha'}\la B\stackrel{\beta'}\la C\stackrel{\gamma'}\la \T A$
with $A,B\in\gen G{2^{k}}$.

Take any $z\in H(C)$. The map $H(\beta'):H(C)\la H(B)$ takes
$z$ to an element $H(\beta')z\in H(B)$.
But the map $\ph:\Hom(B,F_1)\la H(B)$ is surjective,
hence there is an element $g\in\Hom(B,F_1)$ with $\ph(g)=H(\beta')(z)$.
Now $0=H(\alpha')H(\beta')z=H(\alpha')\ph(g)=\ph\Hom(\alpha',F_1)g$,
where the last equality is by the naturality of the map
$\ph:\Hom(-,F_1)|_{\gen G{2^k}}\la H(-)|_{\gen G{2^k}}$. Therefore
$\Hom(\alpha',F_1)g=g\alpha'\in\Hom(A,F_1)$ lies in the kernel
of $\ph:\Hom(A,F_1)\la H(A)$. The exact
sequence
\[\xymatrix@C+50pt{
  \Hom(-,F_2)|_{\gen G{2^k}}  \ar[r]^-{\Hom(-,\alpha)|_{\gen G{2^k}}} &
  \Hom(-,F_1)|_{\gen G{2^k}} \ar[r]^-{\ph} &
  H(-)|_{\gen G{2^k}}
}\]
tells us that there is an $f\in\Hom(A,F_2)$ with
$g\alpha'=\Hom(A,\alpha)f=\alpha f$.

In concrete terms we have produced a commutative diagram
\[\xymatrix@C+20pt{
  A\ar[r]^{\alpha'}\ar[d]^f & B\ar[r]^{\beta'}\ar[d]^g & C
  \ar[r]^{\gamma'} & \T A\\
F_2\ar[r]^{\alpha} & F_1\ar[r]^{\beta} & F
  \ar[r]^{\gamma}& \T F_2
}\]
where the rows are triangles, which we may complete to
a morphism of triangles
\[\xymatrix@C+20pt{
  A\ar[r]^{\alpha'}\ar[d]^f & B\ar[r]^{\beta'}\ar[d]^g & C
  \ar[r]^{\gamma'}\ar[d]^h & \T A\ar[d]^{\T f}\\
F_2\ar[r]^{\alpha} & F_1\ar[r]^{\beta} & F
  \ar[r]^{\gamma}& \T F_2
}\]
We now have an element $h\in\Hom(C,F)$ with $h\beta'=\beta g$,
or in more complicated terms $\Hom(\beta',F)h=\Hom(B,\beta)g$. This
buys us the third equality below
\begin{eqnarray*}
  H(\beta')z &=& \ph(g) \\
  &=& \psi\Hom(B,\beta)(g) \\
  &=&\psi\Hom(\beta',F)h \\
  &=&H(\beta')\psi(h)
\end{eqnarray*}
The first equality is the choice of $g$, the second is because
on the category $\gen G{2^k}$ we have
$\ph=\psi\circ\Hom(-,\beta)$, and the fourth equality is the naturality of
$\psi$. We deduce that $H(\beta'):H(C)\la H(B)$ annihilates
$z-\psi(h)$.

But the exact sequence
$H(\T A)\stackrel{H(\gamma')}\la H(C)\stackrel{H(\beta')}\la H(B)$
tells us there exists a $z'\in H(\T A)$ with
$H(\gamma')z'=z-\psi(h)$. Since $\T A$ belongs
to $\gen G{2^k}$ the map $\ph:\Hom(\T A,F_1)\la H(\T A)$ must be
surjective; there is a $\lambda\in\Hom(\T A,F_1)$ with
$\ph(\lambda)=z'$. Therefore
\begin{eqnarray*}
  z-\psi(h) &=& H(\gamma')\ph(\lambda) \\
  &=& H(\gamma')\psi\Hom(\T A,\beta)(\lambda)\\
  &=&\psi\Hom(\gamma',F)\Hom(\T A,\beta)(\lambda)\\
  &=&\psi(\beta\lambda\gamma')
\end{eqnarray*}
where the second equality is the fact that
on the category $\gen G{2^k}$ we have $\ph=\psi\circ\Hom(-,\beta)$,
the third is the naturality of $\psi$, and the fourth is obvious.
Hence $z=\psi(h+\beta\lambda\gamma')$ is in the image of $\psi$.
\eprf

 When $R$ is a field the theorem below is due to Bondal and Van den
Bergh~\cite[Theorem~1.3]{BondalvandenBergh04}, and
in the generality below it may be found in
Rouquier~\cite[Theorem~4.16 and Corollary~4.18]{Rouquier08}. We include
this new proof because it contains the simple ideas, whose more technical
adaptation will yield the theorems of Sections~\ref{S37} and \ref{SS3}.

\thm{T74.3}
Let $R$ be a noetherian, commutative ring, let $\cs$ be an $R$--linear
triangulated category, and assume that $G\in\cs$ is a
strong generator---we remind the reader, this means that
there exists some integer $n>0$ with $\gen Gn=\cs$. Assume $H$
is a finite cohomological functor, as is $\Hom(-,G)$. Then
there exists a cohomological functor $H'$ with $H\oplus H'$
representable. If $\cs$ is Karoubian, meaning idempotents split,
then $H$ is representable.
\ethm

\prf
Choose an integer $k$ with $\gen G{2^k}=\cs$. Applying
Lemma~\ref{L74.1} to the $\gen G{2^k}$--cohomological functor
$H$ we can find an epimorphism $\Hom(-,F)\la H(-)$.
Complete to an exact sequence $0\la H'(-)\la \Hom(-,F)\la H(-)\la 0$,
and it follows that $H'$ is also a finite cohomological functor.
Applying Lemma~\ref{L74.1} again we have a surjection
$\Hom(-,F')\la H'(-)$, and this assembles to an exact
sequence $\Hom(-,F')\la\Hom(-,F)\la H(-)\la 0$.

Yoneda's lemma tells us that the natural transformation
$\Hom(-,F')\la \Hom(-,F)$ must be of the form
$\Hom(-,\alpha)$ for some $\alpha:F'\la F$, that the natural
transformation $\Hom(-,F)\la H(-)$ corresponds to an element
$y\in H(F)$, and that the vanishing
of the composite $\Hom(-,F')\la\Hom(-,F)\la H(-)$
means that the map $H(\alpha):H(F)\la H(F')$ must take
$y\in H(F)$ to zero. Now complete $\alpha:F'\la F$ to a triangle
$F'\stackrel{\alpha}\la F\stackrel{\beta}\la F''\la$.
The exactness of
$H(F'')\stackrel{H(\beta)}\la H(F)\stackrel{H(\alpha)}\la H(F')$,
coupled with the fact that $H(\alpha)y=0$, 
means that there must be an element $x\in H(F'')$ with
$y=H(\beta)x$.
By Yoneda this means that we obtain a natural transformation
$\Hom(-,F'')\la H(-)$ so that the
diagram below commutes
\[\xymatrix{
  \Hom(-,F')\ar[r]\ar@{=}[d] & \Hom(-,F)\ar[r]\ar@{=}[d]  & \Hom(-,F'')\ar[d]& \\
\Hom(-,F')\ar[r] & \Hom(-,F)\ar[r] & H(-)\ar[r] &  0
}\]
and the rows are exact. It immediately follows that the map
$\Hom(-,F'')\la H(-)$ is a split epimorphism.

Choose a splitting $H(-)\la \Hom(-,F'')$, for example the one
coming from the diagram above. The composite
$\Hom(-,F'')\la H(-)\la\Hom(-,F'')$ is an idempotent
natural endomorphism of a representable functor, therefore
of the form $\Hom(-,e)$ where $e:F''\la F''$ is idempotent.
If $e$ splits then $H$ is representable.
\eprf

\section{Approximating systems}
\label{S27}

In this short section we assemble together a few fairly
routine, elementary facts about
countable direct limits of representable functors. We work in the
generality of $R$--linear functors
between $R$--linear categories, where $R$ is a commutative ring.
What is slightly unusual is that we allow three
interacting categories, 
and hence need a moment's reflection to
make sure that the standard arguments
still work.

\dfn{D27.109}
Let $R$ be a commutative ring, let $\ct$ be an
$R$--linear category, let $\ca,\cb$ be full subcategories of
$\ct$, and let $H:\cb\op\la\Mod R$ be
an $R$--linear functor.
An \emph{$\ca$--approximating system for $H$}
is a sequence in
$\ca$ of  objects and morphisms $E_1\la E_2\la E_3\la\cdots$,
so that
\be
\item
There is a cofinal subsequence of $E_*$ whose objects belong to
$\ca\cap\cb$.
\item
We are given an isomorphism $\colim\,\Hom(-,E_i)\la H(-)$
of functors $\cb\op\la\Mod R$.
\ee
\edfn

\nin
In the proof of
Theorem~\ref{T1.-1}(i) we
will mostly consider the case where $\ca$ is contained
in $\cb$, but in the proof
of Theorem~\ref{T1.-1}(ii)
we will need the more flexible notion.

Since we will freely use approximating systems in our
constructions,
it is comforting to know that they are all the same
up to subsequences. More precisely we have
\lem{L27.113}
Suppose we have an $R$--linear functor $H:\cb\op\la\Mod R$,
and two $\ca$--approximating  systems $E_*$ and $F_*$ for $H$. Then the
systems $E_*$ and $F_*$ are ind-isomorphic. We remind the reader: this means
that 
there exists an $\ca$--approximating system $L_*$ for the functor $H$,
more explicitly 
 $L_1\la L_2\la L_3\la\cdots$, and subsequences $E'_*\subset E_*$, $F'_*\subset F_*$ and $L'_*,L''_*\subset L_*$ with
$E'_*=L'_*$ and $F'_*=L''_*$.
\elem

\prf
Since we may pass to subsequences we will assume that both $E_*$ and $F_*$
belong to $\ca\cap\cb$.
To slightly compress the argument we will extend the sequences by zero;
we will set $E_0=F_0=L_0=0$,
look at the sequences
$E_0\la E_1\la E_2\la\cdots$ and
$F_0\la F_1\la F_2\la\cdots$,
and out of them construct the sequence $L_0\la L_1\la L_2\la\cdots$.

Let $L_0=0$ and $L_1=E_1$, and inductively proceed as follows.
Assume that, for
some $m\geq 1$,
the sequence $L_0\la L_1\la\cdots L_{m-1}\la L_m$ and the map
from the sequence to $H(-)$ have been defined,
in such a way that $L_1\la L_3\la L_5\la\cdots$ is a subsequence of $E_*$
and $L_0\la L_2\la L_4\la L_6\la\cdots$ is a subsequence of $F_*$.
We wish to
extend to $m+1$. There are two cases, $m$ can be odd or even, but up to
interchanging $E_*$ and $F_*$
in the argument below they are the same. We will therefore
assume $m$ odd and leave to the reader the even case.

Then $L_{m-1}$ is 
$F_{I}$ for some $I\geq0$, while
$L_m$ may be assumed equal to $E_J$ for some $J>I$. 
Viewed in the category of functors $(\ca\cap\cb)\op\la\Mod R$,
the map
$\Hom(-,L_m)=\Hom(-,E_J)\la H(-)=\colim\,\Hom(-,F_i)$
is a natural transformation from the representable functor
$\Hom(-,E_J)$ to the colimit. Yoneda tells us that
it corresponds to an element in $H(E_J)=\colim\,\Hom(E_J,F_i)$,
where the colimit is over $i$. We may therefore choose an $I'>J$
and a morphism $L_m=E_J\la F_{I'}$ which delivers the right
element in
the colimit. The composite $F_{I}=L_{m-1}\la L_m\la F_{I'}$
does not have to agree with the map $F_{I}\la F_{I'}$ of the
sequence $F_*$, but they have the same image in
$H(F_I)=\colim\,\Hom(F_{I},F_{i})$
[where the colimit is again over $i$].
That is: after
composing with some $F_{I'}\la F_{I''}$ in the sequence $F_*$ they become
equal.
Set $L_{m+1}=F_{I''}$.
\eprf

\lem{L27.33}
Let $\ct$ be an $R$--linear 
category, let $\ca,\cb$ be full subcategories,
let $H,H':\cb\op\la\Mod R$ be two $R$--linear functors,
and assume we are given for $H$ an $\ca$--approximating
system $E_1\la E_2\la E_3\la\cdots$. Replacing the sequence $E_*$ by
a subsequence belonging to $\ca\cap\cb$, there is a natural
isomorphism $\Hom(H,H')\cong\clim H'(E_m)$ in the category
of functors $\cb\op\la\Mod R$.
\elem

\prf
We have isomorphisms
\begin{eqnarray*}
  \Hom(H,H') &=& \Hom\big(\colim\,\Hom(-,E_m)\,,\,H'(-)\big)\\
  &=& \clim\Hom\big(\Hom(-,E_m)\,,\,H'(-)\big)\\
  &=& \clim H'(E_m)
\end{eqnarray*}
where the last isomorphism is by Yoneda.
\eprf

\cor{C27.5}
Suppose we are given $R$--linear categories $\ca\subset\cb$ and
two $R$--linear functors $H,H':\cb\op\la\Mod R$. If $H$ has
an $\ca$--approximating
system then restriction to the subcategory $\ca\subset\cb$ is a natural
bijection $\Hom(H,H')\la\Hom\big(H|_\ca^{},H'|_\ca^{}\big)$.
\ecor

\prf
Choose an $\ca$--approximating system $E_1\la E_2\la E_3\la\cdots$ for $H$.
Lemma~\ref{L27.33} tells us that both sets are in bijection with
$\clim H'(E_m)$, and the bijection commutes with the restriction map.
\eprf

\lem{L27.592}
Let $\cb$ be an $R$--linear category, let $H,H':\cb\op\la\Mod R$ be two
$R$--linear functors, and let $\ph:H\la H'$ be a natural transformation.
If each of $H,H'$ has a $\cb$--approximating system, let's say
$E_1\la E_2\la E_3\la\cdots$ for $H$ and $E'_1\la E'_2\la E'_3\la\cdots$
for $H'$ then, after replacing $E'_1\la E'_2\la E'_3\la\cdots$
by a subsequence, we can produce a map of sequences $f_*:E_*\la E'_*$
so that $\ph:H\la H'$ is the colimit of the image of
$f_*$ under Yoneda.

Moreover: assume we specify in advance subsequences $\wt E_*\subset E_*$
and $\wt E'_*\subset E'_*$. Then the replacement of $E'_*$ by a subsequence
on which the morphism $f$ is defined, as in the last paragraph, can be
made in such a way that if $E_n$ belongs to $\wt E_*$ then $E'_n$ will belong
to $\wt E'_n$. And we can even handle more than one subsequence: we can
specify in advance a finite set of matching subsequences, pairs
\[
\Big(\{\text{subsequence of $E_*$}\}\quad,\quad
\{\text{subsequence of $E'_*$}\}\Big)\ ,
\]
and, in the replacement of $E'_*$ by a subsequence as in the previous
paragraph, we can guarantee that if $E_n$ belongs to one of the specified
subsequences of $E_*$ then $E'_n$ is chosen to belong to the matching
subsequence of $E'_*$.
\elem

\prf
By Lemma~\ref{L27.33} we have an isomorphism
\[\Hom(H,H')\cong \clim H'(E_i)=\climi\colimj\Hom(E_i, E'_j).
\]
Our element $\ph\in\Hom(H,H')$ must therefore correspond to an inverse
system of elements $\ph_i\in\colimj\Hom(E_i,E'_j)$. We proceed inductively.
\be
\item
  Choose some integer $j_1^{}$ and a preimage in $\Hom(E_1,E'_{j_1^{}})$
  of $\ph_1\in\colimj\Hom(E_1,E_j)$. Call this map $f_1:E_1\la E'_{j_1^{}}$.
  If $E_1$ belongs to one of the prescribed subsequences of $E_*$ then
  choose $j_1^{}$ so that $ E'_{j_1^{}}$ belongs to the matching subsequence
  of $E'_*$.
\item
  Suppose the sequence has been constructed up to an integer $m\geq1$. In
  particular we have a map $E_m\la  E'_{j_m^{}}$, whose image under
  the natural map $\Hom(E_m,E'_{j_m^{}})\la\colimj\Hom(E_m,E'_j)$
  is $\ph_m$.

  We have the element $\ph_{m+1}\in\colimj\Hom(E_{m+1},E'_j)$, we
  can choose a preimage in $\Hom(E_{m+1}, E_J)$ for some integer $J$,
  and we may assume $J>j_{m}^{}$. This gives us a map
  $f':E_{m+1}\la E'_J$. Now the square
  \[\xymatrix{
E_m\ar[r]^{f_m}\ar[d] & E'_{j_m^{}}\ar[d] \\
E_{m+1}\ar[r]^{f'}   & E'_J
  }\]
  need not commute, but the two composites both go, via the map
  $\Hom(E_m,E'_J)\la \colimj\Hom(E_m,E'_j)$, to the same element $\ph_m$.
  Hence replacing $E'_J$ by some $E'_{j_{m+1}^{}}$ with ${j_{m+1}^{}}>J$,
  we may assume the square commutes. And if $E_{m+1}$ belongs to one of
  the prescribed subsequences of $E_*$, choose ${j_{m+1}^{}}>J$ so that
  $E'_{j_{m+1}^{}}$ belongs to the matching subsequence of $E'_*$.
\ee
We have replaced $E'_*$ by a subsequence and produced a map of sequences
$f_*:E_*\la E'_*$.
The reader can check that, if we apply Yoneda to the map of
sequences $f_*:E_*\la E'_*$ and then take colimits, we recover
$\ph:H\la H'$.
\eprf

\rmk{R27.7}
In the remainder of the paper we will
use approximating systems in the following situation.
We will work in some ambient $R$--linear triangulated category $\ct$,
and will assume that $\ct$ has coproducts.
What is special in this case is that, given a functor $H:\cb\op\la\Mod R$
and an $\ca$--approximating
system $E_1\la E_2\la E_3\la\cdots$ for $H$, we can construct in
$\ct$ the homotopy colimit $F=\hoco E_i$.
For $(-)$ in the category $\cb$ we  have a natural map
\[\xymatrix@C+40pt{
H(-)\ar@{=}[r] & \colim\,\Hom(-,E_i)\ar[r] & \Hom(-,F)
}\]
and we will be interested in approximating sequences for which this
map $H(-)\la \Hom(-,F)$ is an isomorphism of functors $\cb\op\la\Mod R$.

In this situation we will say that $E_1\la E_2\la E_3\la\cdots$
is an $\ca$--approximating system for $F$ over $\cb$. 
\ermk

\rmk{R27.99309}
In the proof of Theorem~\ref{T1.-1}(i), the case of interest is where 
$\cb\subset\ct^c$, that is the objects of $\cb$ are all compact.
From
\cite[Lemma~2.8]{Neeman96} we know that, for any compact
object $K\in\ct$ and any sequence of objects of $\ct$,
the natural map is an isomorphism
$H(K)\cong\colim\,\Hom(K,E_i)\la \Hom\big(K,\hoco E_i\big)$.
Thus we're automatically in the situation of Remark~\ref{R27.7};
any sequence $E_*$ in $\ca$ is an $\ca$--approximating system for
$F=\hoco E_i$, over any $\cb\subset\ct^c$.
In Section~\ref{S37} we
will therefore
allow ourselves to occasionally leave unspecified the category
$\cb\subset\ct^c$, and just say that $E_*$ is an $\ca$--approximating
system for $F$.

But here we are careful to specify $\cb$, because
this will come up in the proof of Theorem~\ref{T1.-1}(ii).
\ermk

In the generality of Remark~\ref{R27.7} we note the following
little observation.

\lem{L27.976}
Let $\ct$ be a triangulated category with coproducts, let $\ca,\cb$
be subcategories, and assume $E_1\la E_2\la E_3\la\cdots$ is
a $\ca$--approximating system for $F\in\ct$ over $\cb$.

If $G$ is another object of $\ct$, and
if $\ph:\Hom(-,F)|_\cb\la\Hom(-,G)|_\cb$
is a natural transformation of functors on $\cb$, then there exists in
$\ct$ a
(non-unique) morphism $f:F\la G$ with $\ph=\Hom(-,f)|_\cb$.
\elem

\prf
By Lemma~\ref{L27.33} the natural transformation $\ph:\Hom(-,F)\la \Hom(-,G)$
corresponds to an element in $\clim\Hom(E_i,G)$. Thus for each $i$ we are
given in $\ct$ a
morphism $f_i:E_i\la G$, compatibly with the sequence maps $E_i\la E_{i+1}$.
The compatibility means that the composite
\[\xymatrix@C+20pt{
  \ds\bigoplus_{i=1}^\infty E_i\ar[rr]^-{1-\text{\rm shift}} &&
  \ds\bigoplus_{i=1}^\infty E_i \ar[rr]^-{(f_1^{},f_2,f_3,\ldots)} && G
}\]
must vanish. Hence the map $\bigoplus_{i=1}^\infty E_i\la G$ factors
(non-uniquely) through
$F=\hoco E_i$, which is the third edge in the triangle
\[\xymatrix@C+20pt{
  \ds\bigoplus_{i=1}^\infty E_i\ar[rr]^-{1-\text{\rm shift}} &&
  \ds\bigoplus_{i=1}^\infty E_i \ar[rr]^-{(f_1^{},f_2,f_3,\ldots)} && \hoco E_i\ .
}\]
\eprf

\section{A couple of technical lemmas}
\label{S997}

In Section~\ref{S37} we will prove Theorem~\ref{T1.-1}(i),
and in Section~\ref{SS3} we will prove
Theorem~\ref{T1.-1}(ii). Both proofs
will rely heavily on
a couple of technical lemmas---in this section
we state and prove these lemmas,
in sufficient generality to
cover both applications. Let us therefore set
up a little notation.

\ntn{N997.1}
Throughout this section $R$ will be a commutative
ring, $\ct$ will be an $R$--linear triangulated category with
coproducts, and $\cs\subset\ct$
will be a triangulated subcategory.
The Yoneda functor $\cy:\ct\la\Hom_R^{}\big[\cs\op\,,\,\Mod R\big]$ will
be the map taking $t\in\ct$ to $\Hom(-,t)$, where $\Hom(-,t)$ is
viewed as an
$R$--linear functor $\cs\op\la\Mod R$.

We remind the reader of
Definition~\ref{D27.109}: suppose $\ca$ is a full subcategory of $\ct$
closed under
direct summands, finite coproducts and suspensions, and $H$
is an $\cs$--cohomological functor,
meaning $H:\cs\op\la\Mod R$ is
an $R$--linear cohomological functor. Then a
$\ca$--approximating system for $H$ is a sequence
$E_1\la E_2\la E_3\la\cdots$ in $\ca$, with a subsequence in $\ca\cap\cs$,
and an isomorphism $\colim\,\cy(E_i)\la H(-)$.
\entn

\dfn{D997.203}
Let $\ct$ be an $R$--linear triangulated category with coproducts and let $\cs$
be a triangulated subcategory. If $H$ is a $\cs$--cohomological functor,
we define
$\T H$ by the rule $\T H(s)=H(\Tm s)$.

A \emph{weak triangle} in the category $\Hom_R^{}\big[\cs\op\,,\,\Mod R\big]$
is a sequence  of cohomological
functors $A\stackrel u\la B\stackrel v\la C\stackrel w\la\T A$
such that
any rotation of the following is true: given any triangle
$a\stackrel {u'}\la b\stackrel {v'}\la c\stackrel {w'}\la \T a$
in the category $\cs$ and a 
  commutative diagram
  \[\xymatrix{
    \cy(a)\ar[r]^{\cy(u')}\ar[d]^{f} & \cy(b)\ar[r]^{\cy(v')}\ar[d]^{g} &
    \cy(c)\ar[r]^{\cy(w')} &\cy( \T a) \\
A\ar[r]^{u} & B\ar[r]^{v} &
    C\ar[r]^{w} & \T A 
  }\]
  there is an  extension 
  to a commutative diagram
  \[\xymatrix{
   \cy(a)\ar[r]^{\cy(u')}\ar[d]^{f} & \cy(b)\ar[r]^{\cy(v')}\ar[d]^{g} &
    \cy(c)\ar[r]^{\cy(w')}\ar[d]^{h}  &\cy( \T a)\ar[d]^{\T f} \\
A\ar[r]^{u} & B\ar[r]^{v} &
    C\ar[r]^{w} & \T A 
  }\]
A diagram
$\wh A\stackrel{\wh u}\la \wh B\stackrel{\wh v}\la\wh C\stackrel{\wh w}\la\T \wh A$ in the category $\ct$ is called a weak triangle if the functor
$\cy$ takes it to a weak triangle in $\Hom_R^{}\big[\cs\op\,,\,\Mod R\big]$.
\edfn

\rmk{R997.207}
We remind the reader of Remark~\ref{R27.7}
and Lemma~\ref{L27.976}: if
$A$, $B$ and $C$ have $\ct$--approximating systems
$\mathfrak{A}_*$, $\mathfrak{B}_*$ and $\mathfrak{C}_*$, we may form in $\ct$
the
homotopy colimits $\wh A=\hoco\mathfrak{A}_*$, $\wh B=\hoco\mathfrak{B}_*$
and $\wh C=\hoco\mathfrak{C}_*$. Remark~\ref{R27.7}
tells us that 
there are canonical maps $\alpha:A\la\cy(\wh A)$, $\beta:B\la\cy(\wh B)$
and $\gamma:C\la\cy(\wh C)$. Since our plan is to apply the lemmas in
this section to prove representability theorems, we will mostly
be interested in cases where $\alpha$, $\beta$ and $\gamma$ are
isomorphisms. In this case Lemma~\ref{L27.976}
says that the maps $u$, $v$ and $w$ may be lifted (non-uniquely) to $\ct$;
we may form in $\ct$
a diagram
$\wh A\stackrel{\wh u}\la \wh B\stackrel{\wh v}\la\wh C\stackrel{\wh w}\la\T \wh A$
whose image under $\cy$ is (canonically) isomorphic to 
$A\stackrel u\la B\stackrel v\la C\stackrel w\la\T A$.
\ermk

\lem{L997.362}
Suppose $A\la B\la C\la\T A$ is a weak triangle in
$\Hom_R^{}\big[\cs\op\,,\,\Mod R\big]$. For any
$s\in\cs$ the functor $\Hom\big(\cy(s),-\big)$ takes it to an exact sequence.
\elem

\prf
Given any map $f:\cy(s)\la A$, we can consider the
commutative diagram
  \[\xymatrix{
    \cy(s)\ar@{=}[r]\ar[d]_-{f} & \cy(s)\ar[r]\ar[d]^{uf} &
    0\ar[r] & \cy(\T s) \\
A\ar[r]^{u} & B\ar[r]^{v} &
    C\ar[r]^{w} & \T A 
  }\]
  and the fact that this commutative diagram can be extended gives
  the vanishing of $vuf$.
  
A morphism $g:\cy(s)\la B$ so that $vg=0$ gives a
commutative diagram
  \[\xymatrix{
    \cy(s)\ar@{=}[r] & \cy(s)\ar[r]\ar[d]^{g} &
    0\ar[r]\ar[d] & \cy(\T s) \\
A\ar[r]^{u} & B\ar[r]^{v} &
    C\ar[r]^{w} & \T A 
  }\]
and the existence of an extension gives a morphism $f:\cy(s)\la A$ with
$g=uf$. 
\eprf

\lem{L1.35}
With the conventions of Notation~\ref{N997.1}
and Definition~\ref{D997.203}
suppose we are given:
\be
\item
A morphism
$\alpha: A\la B$ in the category $\Hom_R^{}\big[\cs\op\,,\,\Mod R\big]$.  
\item
A morphism $\alpha_*:\mathfrak{A}_*\la\mathfrak{B}_*$ of sequences in $\ct$,
and an isomorphism in $\Hom_R^{}\big[\cs\op\,,\,\Mod R\big]$ of
$\alpha:A\la B$ with the colimit of
$\cy(\alpha_*):\cy(\mathfrak{A}_*)\la\cy(\mathfrak{B}_*)$.
\item
The sequence $\alpha_*$ is assumed to have a subsequence in $\cs$.
\setcounter{enumiv}{\value{enumi}}
\ee
With just these hypotheses we may complete
$\alpha_*:\mathfrak{A}_*\la\mathfrak{B}_*$ to a sequence
$\mathfrak{A}_*\stackrel{\alpha_*}\la\mathfrak{B}_*\stackrel{\beta_*}\la\mathfrak{C}_*\stackrel{\gamma_*}\la\T\mathfrak{A}_*$
of triangles
in $\ct$,
and the colimit
of
$\cy(\mathfrak{A}_*)\stackrel{\cy(\alpha_*)}\la\cy(\mathfrak{B}_*)\stackrel{\cy(\beta_*)}\la\cy(\mathfrak{C}_*)\stackrel{\cy(\gamma_*)}\la\T\cy(\mathfrak{A}_*)$
is a weak triangle 
$A\stackrel{\alpha}\la B\stackrel{\beta}\la C\stackrel{\gamma}\la \T A$.

Suppose we add the following assumptions:
\be
\setcounter{enumi}{\value{enumiv}}
\item
We are given
two subcategories $\ca\subset\cb\subset\cs$, closed under
finite coproducts, direct summands and
suspensions.
\item
There is a subsequence of $\alpha_*:\mathfrak{A}_*\la\mathfrak{B}_*$
such that the $\mathfrak{A}_i$ belongs to $\ca$ and the
$\mathfrak{B}_i$ belong
to $\cb$. Put $\cc=\Smr(\cb*\ca)$.
\item
Assume furthermore that we are given a $\cc$--cohomological
functor $H$ and a natural
transformation  of $\cb$--cohomological
functors $\ph: B|_\cb\la H|_\cb$. Assume that, on the category
$\ca\subset\cb$, the composite
\[\xymatrix@C+30pt{
A|_\ca\ar[r]^-{\alpha|_\ca} & B|_\ca\ar[r]^-{\ph|_\ca} & H|_\ca
}\]
vanishes.
\item
Assume further that the approximating system 
$\mathfrak{A}_*$ for $A$ is such that
each morphism $\mathfrak{A}_i\la \mathfrak{A}_{i+1}$
is a split monomorphism.
\setcounter{enumiv}{\value{enumi}}
\ee
Then there exists a map $\psi:C|_{\cc}\la H$
so that
$\ph:B|_\cb\la H|_\cb$ is equal to the composite
\[\xymatrix@C+30pt{
B|_\cb\ar[r]^-{\beta|_{\cb}} & C|_\cb\ar[r]^-{\psi|_\cb} & H|_\cb\ .
}\]
\elem

\prf
We are given a morphism of sequences $\alpha_*:\mathfrak{A}_*\la\mathfrak{B}_*$,
meaning for each $m>0$ we have a commutative square
\[\xymatrix@C+30pt{
  \mathfrak{A}_m\ar[r]^{\alpha_m}\ar[d] & \mathfrak{B}_m\ar[d] \\
  \mathfrak{A}_{m+1}\ar[r]^{\alpha_{m+1}} & \mathfrak{B}_{m+1}
}\]
We extend this to a morphism of triangles
\[\xymatrix@C+30pt{
  \mathfrak{A}_m\ar[r]^{\alpha_m}\ar[d] & \mathfrak{B}_m\ar[d]\ar[r]^{\beta_m}
    & \mathfrak{C}_m\ar[r]^{\gamma_m}\ar[d]&\T \mathfrak{A}_m\ar[d]\\
  \mathfrak{A}_{m+1}\ar[r]^{\alpha_{m+1}}& \mathfrak{B}_{m+1}\ar[r]^{\beta_{m+1}}&
  \mathfrak{C}_{m+1}\ar[r]^{\gamma_{m+1}} &\T \mathfrak{A}_{m+1}
}\]
This
produces for us in $\ct$ the sequence of
triangles
$\mathfrak{A}_*\stackrel{\alpha_*}\la\mathfrak{B}_*\stackrel{\beta_*}\la\mathfrak{C}_*\stackrel{\gamma_*}\la\T\mathfrak{A}_*$, with a subsequence in $\cs$,
and it is easy to see that the colimit
of
$\cy(\mathfrak{A}_*)\stackrel{\cy(\alpha_*)}\la\cy(\mathfrak{B}_*)\stackrel{\cy(\beta_*)}\la\cy(\mathfrak{C}_*)\stackrel{\cy(\gamma_*)}\la\T\cy(\mathfrak{A}_*)$
is a weak triangle 
$A\stackrel{\alpha}\la B\stackrel{\beta}\la C\stackrel{\gamma}\la \T A$.

It remains to prove the part with the further assumptions added. Note
that, by passing to a subsequence, we may assume $\mathfrak{A}_i\in\ca$
and $\mathfrak{B}_i\in\cb$,
and hence $\mathfrak{C}_i\in\cb*\ca\subset\cc$.
As $H$ is a cohomological functor on $\cc$ and $\ca\subset\cb\subset\cc$,
we have that, for each integer $m$, the
sequence
$H(\T \mathfrak{A}_m)\stackrel{H(\gamma_m)}\la H(\mathfrak{C}_m)\stackrel{H(\beta_m)}\la H(\mathfrak{B}_m)\stackrel{H(\alpha_m)}\la H(\mathfrak{A}_m)$
must be exact. As $m$ increases this 
gives an inverse system of exact sequences, which we now propose to
analyze. The short exact sequences
$0\la \Ker\big(H(\gamma_m)\big)
\la H(\T \mathfrak{A}_m)\la \Ima\big(H(\gamma_m)\big)\la 0$
give an exact sequence
\[
\xymatrix@C+20pt{
  \clim^1\, H(\T \mathfrak{A}_m) \ar[r]  & \clim^1\,\Ima\big(H(\gamma_m)\big)\ar[r] &
  \clim^2\, \Ker\big(H(\gamma_m)\big)
}
\]
In (vii) we assumed that the maps $\mathfrak{A}_m\la \mathfrak{A}_{m+1}$
are split monomorphisms, hence
the maps $H(\T \mathfrak{A}_{m+1})\la H(\T \mathfrak{A}_m)$ are split
epimorphisms, making the sequence
Mittag-Leffler. Therefore $\clim^1\,H(\T \mathfrak{A}_m)=0$.
We have
$\clim^2\, \Ker\big(H(\gamma_m)\big)=0$ just because we're dealing
with a countable limit. We conclude that
$\clim^1\,\Ima\big(H(\gamma_m)\big)=0$.

Now consider the inverse system of short exact sequences
$0\la \Ima\big(H(\gamma_m)\big)\la
H(\mathfrak{C}_m)\la \Ima\big(H(\beta_m)\big)\la 0$.
Passing to the limit we obtain an exact sequence
\[
\xymatrix@C+20pt{
  \clim H(\mathfrak{C}_m) \ar[r]  & \clim\Ima\big(H(\beta_m)\big)\ar[r] &
  \clim^1\, \Ima\big(H(\gamma_m)\big)
}
\]
We have proved the vanishing of $\clim^1\, \Ima\big(H(\gamma_m)\big)$,
allowing us to
conclude that the map $\clim H(\mathfrak{C}_m) \la \clim\Ima\big(H(\beta_m)\big)$
is an epimorphism. Finally we observe the exact sequences
$0\la \Ima\big(H(\beta_m)\big) \la H(\mathfrak{B}_m)\la H(\mathfrak{A}_m)$
and, since inverse limit is left exact, we deduce the exactness of
\[
\xymatrix@C+20pt{
0\ar[r]  & \clim\Ima\big(H(\beta_m)\big)\ar[r] &
  \clim H(\mathfrak{B}_m)\ar[r] & \clim H(\mathfrak{A}_m)
}
\]
Combining the results we have the exactness of
\[
\xymatrix@C+20pt{
\clim H(\mathfrak{C}_m)\ar[r] &
  \clim H(\mathfrak{B}_m)\ar[r] & \clim H(\mathfrak{A}_m)\ .
}
\]
Now Lemma~\ref{L27.33} tells us that $\ph:B|_\cb\la H|_\cb$
corresponds to an element $f\in\clim H(\mathfrak{B}_m)$, and the
vanishing of the composite
$A|_\ca\la B|_\ca\la H|_\ca$
translates to saying that the image of $f$ under the map
$\clim H(\mathfrak{B}_m)\la\clim H(\mathfrak{A}_m)$ vanishes.
The exactness tells us that $f$ is in the image of the
map $\clim H(\mathfrak{C}_m)\la\clim H(\mathfrak{B}_m)$. This exactly says
that there is a natural transformation $\psi:C|_\cc\la H$
with $\ph=\psi\circ\beta$. We have proved the ``extra assumptions'' part.
\eprf

\lem{L997.601}
With the conventions of Notation~\ref{N997.1}
and Definition~\ref{D997.203}
suppose we are given:
\be
\item
Two full subcategories $\ca\subset\cb$ of the category $\cs$, closed under
finite coproducts, direct summands and
suspensions.
\item
Put $\cc=\Smr(\cb*\ca)$. Assume
we are also given a $\cc$--cohomological functor
$H$.
\item
We are given a weak triangle 
$A\stackrel{\alpha}\la B\stackrel{\beta}\la C\stackrel{\gamma}\la \T A$,
and a natural transformation of $\cc$--cohomological functors
$\psi:C|_\cc\la H$.
\item
The composite $(\psi\beta)|_\cb:B|_\cb\la H|_\cb$ is surjective.
\item
The sequence
\[\xymatrix@C+30pt{
A|_\ca\ar[r]^-{\alpha|_\ca} & B|_\ca\ar[r]^-{(\psi\beta)|_\ca} & H|_\ca\ar[r] &0
}\]
is exact.
\ee
Then the map $\psi:C|_\cc\la H$ 
an epimorphism.
\elem

\prf
We need to show the surjectivity of the map $\psi:C(c)\la H(c)$ for every
$c\in\cc=\Smr(\cb*\ca)$; without loss of generality we may assume
$c\in\cb*\ca$.
Choose a triangle
$a\stackrel{\alpha'}\la b\stackrel{\beta'}\la c\stackrel{\gamma'}\la \T a$
with $b\in\cb$ and $a\in\ca$. Given any element $y\in H(c)$,
the map $H(\beta'):H(c)\la H(b)$ takes $y$ to an element $H(\beta')(y)$
which must be in the image of the surjective map $\psi\beta:B(b)\la H(b)$.
After all $b$ is an object in $\cb$, and the map
$\psi\beta:B(b)\la H(b)$ is an epimorphism on objects $b\in\cb$.
Choose an element $g\in B(b)$ mapping under $\psi\beta$ to $H(\beta')(y)$.
The naturality of $\psi\beta$ means that the square below commutes
\[\xymatrix@C+40pt{
  & B(b)\ar[r]^-{B(\alpha')}\ar[d]^{\psi\beta}
  & B(a) \ar[d]^{\psi\beta} \\
  H(c) \ar[r]^-{H(\beta')} & H(b) \ar[r]^-{H(\alpha')} & H(a)
}\]
If we apply the equal composites in the square to $g\in B(b)$ we discover
that it goes to $H(\alpha')H(\beta')(y)=0$, where the vanishing
is because $\beta'\alpha'=0$.
Therefore the map $B(\alpha')$ takes $g\in B(b)$
to an element in the kernel of $\psi\beta: B(a)\la H(a)$.
As $a$ belongs to $\ca$ the map
$\alpha:A(a)\la B(a)$ surjects onto
this kernel; there is an element $f\in A(a)$ with
$\alpha(f)=B(\alpha')(g)$.

We have produced elements $f\in A(a)$ and $g\in B(b)$, and
Yoneda allows us to view them as
natural transformations $f:\cy(a)\la A$ and $g:\cy(b)\la B$.
The equality $\alpha(f)=B(\alpha')(g)$ transforms into the assertion that
the square below commutes
\[\xymatrix@C+40pt{
  \cy(a)\ar[r]^-{\cy(\alpha')}\ar[d]^f & \cy(b)\ar[r]^-{\cy(\beta')} \ar[d]^g & \cy(c)\ar[r]^-{\cy(\gamma')} & \cy(\T a) \\
  A\ar[r]^-{\alpha} & B\ar[r]^-{\beta}&  C \ar[r]^-\gamma & \T A 
}\]
The top row is the image under Yoneda of a triangle in $\cs$,
while the bottom row is a
weak triangle; hence we may complete to a commutative diagram
\[\xymatrix@C+40pt{
  \cy(a)\ar[r]^-{\cy(\alpha')}\ar[d]^f & \cy(b)\ar[r]^-{\cy(\beta')} \ar[d]^g & \cy(c)\ar[r]^-{\cy(\gamma')}\ar[d]^h & \cy(\T a)\ar[d]^{\T f} \\
  A\ar[r]^-{\alpha} & B\ar[r]^-{\beta}&  C \ar[r]^-\gamma & \T A 
}\]
We have produced a morphism $h:\cy(c)\la C$, which we may view as an
element $h\in C(c)$. And the commutativity of the
middle square translates, under Yoneda, to the statement that
$C(\beta'):C(c)\la C(b)$ takes $h\in C(c)$ to $\beta(g)\in C(b)$.
Applying
$\psi:C\la H$ we obtain the second equality below
\begin{eqnarray*}
  H(\beta')(y) &=& \psi\beta(g)\\
  &=&\psi C(\beta')(h) \\
  &=&H(\beta')\psi(h)\ .
\end{eqnarray*}
The first equality is by construction of $g\in B(b)$,
and the third is
the naturality of $\psi$. Therefore the map
$H(\beta'):H(c)\la H(b)$ annihilates
$y-\psi(h)$. Because $H$ is cohomological there is
an element $x\in H(\T a)$ with $H(\gamma')(x)=y-\psi(h)$.
But $a\in\ca\subset\cb$ and we may choose
a $\theta\in B(\T a)$ with $\psi\beta(\theta)=x$. We
have
\begin{eqnarray*}
  y-\psi(h)
  &=& H(\gamma')\psi\beta(\theta)\\
  &=&\psi C(\gamma')\beta(\theta)
\end{eqnarray*}
where the first equality is the construction of $\theta$,
and the second is
the naturality of $\psi:C\la H$.
These equalities combine to the formula
$y=\psi\big[h+C(\gamma')\beta(\theta)\big]$, which
exihibits $y\in H(c)$ as lying in the image
of $\psi:C(c)\la H(c)$.
\eprf

The next lemma will not be 
needed until the proof of Theorem~\ref{TT3.13}.

\lem{L997.223}
With the conventions of Notation~\ref{N997.1}
and Definition~\ref{D997.203}
suppose we are given:
\be
\item
Two full subcategories $\ca\subset\cb$ of the category $\cs$, closed under
finite coproducts, direct summands and
suspensions.
\item
Put $\cc=\Smr(\cb*\ca)$.
\item
In the category $\Hom_R^{}\big[\cs\op\,,\,\Mod R\big]$ we are given a diagram
of cohomological functors
\[\xymatrix@C+20pt{
                  & \wt B\ar[d]^-\delta  \\
A\ar[r]^-{\alpha} & B\ar[r]^-{\beta} & C\ar[r]^-{\gamma}\ar[d]^-\psi & \T A\\
  & & H
}\]
where the middle row is a weak triangle.
\item
The composite $(\psi\beta)|_\cb^{}: B|_\cb^{}\la H|_\cb^{}$ is surjective.
\item
The kernel of the morphism $(\psi\beta\delta)|_\cb:\wt B|_\cb\la H|_\cb$
is annihilated by $\delta:\wt B\la B$.
\item
The sequence
\[\xymatrix@C+30pt{
A|_\ca\ar[r]^-{\alpha|_\ca} & B|_\ca\ar[r]^-{(\psi\beta)|_\ca} & H|_\ca\ar[r] &0
}\]
is exact.
\setcounter{enumiv}{\value{enumi}}
\ee
Then the map $(\beta\delta):\wt B\la C$ annihilates the kernel of
$(\psi\beta\delta)|_\cc:\wt B|_\cc\la H|_\cc$.
\elem

\prf
We need to show that, if $c\in\cc=\Smr(\cb*\ca)$ and $y\in\wt B(c)$
is annihilated by the map $\psi\beta\delta:\wt B(c)\la H(c)$, then
$y$ is already annihilated by the shorter map $\beta\delta:\wt B(c)\la C(c)$.
Note that without loss of generality we may assume $c\in\cb*\ca$. Choose
therefore a triangle
$a\stackrel{\alpha'}\la b\stackrel{\beta'}\la c\stackrel{\gamma'}\la \T a$
with $b\in\cb$ and $a\in\ca$, and consider the commutative diagram with
exact rows
\[\xymatrix@C+20pt{
 & & \wt B(c)\ar[d]^-{\delta}\ar[r]^-{\wt B(\beta')} &
\wt B(b)\ar[d]^-{\delta} \\
B(\T b)\ar[d]^-{\psi\beta}\ar[r]^-{B(\T\alpha')}  &
B(\T a)\ar[d]^-{\psi\beta}\ar[r]^-{B(\gamma')}  &
 B(c)\ar[d]^-{\psi\beta}\ar[r]^-{B(\beta')} &
B(b)\ar[d]^-{\psi\beta} \\
H(\T b)\ar[r]^-{H(\T\alpha')} &H(\T a)\ar[r]^-{H(\gamma')} &  H(c)\ar[r]^-{H(\beta')} & H(b)
}\]
We are given an element $y\in\wt B(c)$ such that the vertical composite
in the third column annihilates it. Therefore $\wt B(\beta')(y)$ is an
element of $\wt B(b)$ annihilated by the vertical composite in
the fourth column. By assumption (v) the element $\wt B(\beta')(y)$ is already
killed by $\delta:\wt B(b)\la B(b)$, and we conclude that the equal
composites in the top-right square annihilate $y$. Hence
the map $\delta:\wt B(c)\la B(c)$ must take $y\in\wt B(c)$ to an element
in the image of $B(\gamma')$, and we may therefore
\be
\setcounter{enumi}{\value{enumiv}}
\item
Choose an $x\in B(\T a)$ with
$B(\gamma')(x)=\delta(y)$.
\setcounter{enumiv}{\value{enumi}}
\ee
Recall that the vertical composite in the third column kills $y$,
hence the equal composites in the middle square at the bottom must annihilate
$x$. Therefore $\psi\beta(x)\in H(\T a)$ lies in the kernel of $H(\gamma')$,
which is the image of $H(\T\alpha'):H(\T b)\la H(\T a)$. Now
by (iv)
the map $\psi\beta:B(\T b)\la H(\T b)$ is surjective,
and we may
lift further to $B(\T b)$; we can choose an element $w\in B(\T b)$ whose
image under the equal composites in the bottom-left square are equal to
$\psi\beta(x)$. Therefore we have that,
with $x\in B(\T a)$ as in (vii) and $w\in B(\T b)$ as above, the element
$x-B(\T \alpha')(w)$ is annihilated by $\psi\beta:B(\T a)\la H(\T a)$.
Now (vi) tells us that 
\be
\setcounter{enumi}{\value{enumiv}}
\item
We may choose an element $v\in A(\T a)$ whose
image under $\alpha:A(\T a)\la B(\T a)$ is equal to $x-B(\T \alpha')(w)$.
\setcounter{enumiv}{\value{enumi}}
\ee
To complete the proof consider the commutative diagram with
vanishing horizontal and vertical composites 
\[\xymatrix@C+20pt{
 & A(\T a)\ar[d]^-{\alpha}\ar[r]^-{A(\gamma')} &
A(c)\ar[d]^-{\alpha} \\
B(\T b)\ar[r]^-{B(\T\alpha')}  &
B(\T a)\ar[d]^-{\beta}\ar[r]^-{B(\gamma')}  &
 B(c)\ar[d]^-{\beta} \\
 &C(\T a)\ar[r]^-{C(\gamma')} &  C(c)
}\]
The horizontal map in the top row
takes the element $v\in A(\T a)$ constructed in (viii) to
$A(\gamma')(v)$, which must be annihilated by the vertical composite
in the third column. By the commutativity of the top square coupled
with (viii), this means that $x-B(\T \alpha')(w)$ is an element
of $B(\T a)$ annihilated by the equal composites in the bottom square.
In particular the horizontal map $B(\gamma')$ takes
$x-B(\T \alpha')(w)$ to an element of the kernel
of $\beta:B(c)\la C(c)$. But the map $B(\gamma')$ annihilates
$B(\T \alpha')(w)$, and by (vii) it takes $x$ to $\delta(y)$. We conclude
that $\beta\delta(y)=0$.
\eprf

\section{The proof of Theorem~{\protect{\ref{T1.-1}}}($\mathrm{i}$)}
\label{S37}

It's time to prove Theorem~\ref{T1.-1}(i); we should
focus the general lemmas of Section~\ref{S997}
on the situation at hand. Thus 
in this section we make the following global assumptions:

\ntn{N37.9993}
We specialize the conventions of Notation~\ref{N997.1}
by setting $\cs=\ct^c$, that
is $\cs$ is the subcategory of compact objects in $\ct$.
Thus in this section the functor $\cy$ of Notation~\ref{N997.1}
specializes to
$\cy:\ct\la\Hom_R^{}\big([\ct^c]\op\,,\,\Mod R\big)$,
which takes
an object $t\in\ct$ to $\cy(t)=\Hom(-,t)|_{\ct^c}$.

In the generality of Notation~\ref{N997.1}
we considered $\cs$--cohomological functors $H$ and 
$\cb$--approximating systems $E_1\la E_2\la E_3\la\cdots$.
Because we are now in the special case where $\cs=\ct^c$
Remark~\ref{R27.99309} applies: if $F=\hoco E_i$ then the
natural map $H\la\cy(F)$ must be an isomorphism.

We will furthermore assume that we have chosen in $\ct$
a single compact generator $G$.
We will suppose given
a \tstr\ $(\ct^{\leq0},\ct^{\geq0})$ in the preferred
equivalence class, with $\ct^{\geq0}$ closed under coproducts.
For
example we could let $(\ct^{\leq0},\ct^{\geq0})$ be equal to
$\big(\ct_G^{\leq0},\ct_G^{\geq0}\big)$, see Remark~\ref{R0.7}.
Let $\ca$ be the heart of the \tstr, and
 $\ch:\ct\la\ca$ the homological functor of Reminder~\ref{R-1.1}.
We will assume given an integer $A>0$ with
$\Hom\big(\T^{-A}G,\ct^{\leq0}\big)=0$. The existence of
such an $A$ is equivalent to the hypothesis that $\Hom(G,\T^iG)=0$
for $i\gg0$. And what is important for us
is that this guarantees that the category $\ct^-_c$ is a thick
subcategory of $\ct$.

Finally and most importantly: as in Notation~\ref{N997.1} the
triangulated category $\ct$ is assumed to be $R$--linear for some
commutative ring $R$. But from now on we add the assumption that
the ring $R$ is noetherian and, with $G$ as in the last paragraph,
the $R$--module $\Hom(G,\T^iG)$ is finite for every $i\in\zz$.
Combining this paragraph with the last: the functor $\cy(G)$ is $G$--locally
finite.
\entn

Under some additional approximability
assumptions, Theorem~\ref{T1.-1} 
describes the essential image of the functor $\cy$ taking
$F\in\ct^-_c$ to
$\cy(F)=\Hom(-,F)|_{\ct^c}$, and tells us that the functor $\cy$ is full.
To show that $\cy(F)$ lies in the expected image
one doesn't need any hypotheses beyond the ones above, we prove

\lem{L1.1.3}
With the assumptions of Notation~\ref{N37.9993},
for any $F\in\ct^-_c$ the functor $\cy(F):[\ct^c]\op\la\Mod R$
is a locally finite $\ct^c$--cohomological functor.
\elem

\prf
We are given that the functor $\cy(G)$ is $G$--locally
finite.
In particular: for $i\ll0$ we have $\Hom(\T^iG,G)=0$.
With $(\ct^{\leq0},\ct^{\geq0})$ as in
Notation~\ref{N37.9993}, that is our fixed \tstr\ in the
preferred equivalence class, Remark~\ref{R0.14.7} coupled with
Lemma~\ref{L17.17} tell us that, for any object $K\in\ct^c$,
there is an integer $B>0$ so that
$\Hom\big(\T^{-B}K,\ct^{\leq0}\big)=0$. Remark~\ref{R0.13.99}
gives the inclusion in $F\in\ct^-_c\subset\ct^-$, hence we may
choose an integer $A>0$ with $\T^AF\in\ct^{\leq0}$. We deduce
that $\Hom(\T^{i}K,F)=0$ for all $i\leq-A-B$.

The fact that the functor $\cy(G)$ is $G$--locally
finite also means that, for every integer $i\in\zz$,
the $R$--module $\Hom(\T^iG,G)$ is finite. The full subcategory
$\cl\subset\ct^c$ defined by
\[
\cl\eq\{L\in\ct^c\mid\Hom(\T^iG,L)\text{ is a finite $R$--module for all }i\in\zz\}
\]
is thick and contains $G$, hence $\ct^c=\gen G{}\subset\cl$. Now take
any $L\in\ct^c$ and define the full subcategory $\ck(L)\subset\ct^c$ by
\[
\ck(L)\eq\{K\in\ct^c\mid\Hom(\T^iK,L)\text{ is a finite $R$--module for all }i\in\zz\}
\]
Then $\ck(L)$ is thick and contains $G$, hence $\ct^c=\gen G{}\subset\ck(L)$.
We conclude that $\Hom(\T^iK,L)$ is a finite $R$--module for all $K,L\in\ct^c$
and all $i\in\zz$.

Now fix the integer $i$, the object $K\in\ct^c$ and the object $F\in\ct^-_c$,
and we want to prove that
$\Hom(\T^iK,F)$ is a finite $R$--module. The first paragraph of the proof
produced an integer $B>0$ with $\Hom\big(\T^{-B}K,\ct^{\leq0}\big)=0$, and
since $F$ belongs to $\ct^-_c$
there exists a triangle $L\la F\la D$ with $L\in\ct^c$ and
$D\in\ct^{\leq-i-B-1}$.
In the exact sequence
\[\xymatrix{
\Hom(\T^iK,\Tm D)\ar[r] &
\Hom(\T^iK,L)\ar[r] &
\Hom(\T^iK,F)\ar[r] &
\Hom(\T^iK,D)
}\]
we have that
$\Hom(\T^iK,D)=0
=\Hom(\T^iK,\Tm D)$,
and $\Hom(\T^iK,F)\cong \Hom(\T^iK,L)$ must be a finite
$R$--module by the second paragraph of the proof.
\eprf

\nin
In
Corollary~\ref{C17.29} we learned that, under some
approximability hypotheses, objects in $\ct^-_c$
can be well approximated by sequences with special properties.
We don't need all these properties yet; for the
next few lemmas we formulate
what we will use.

\dfn{D1.29}
Adopting the conventions of Notation~\ref{N37.9993},
a \emph{strong $\gen Gn$--approximating system} is
a sequence of objects and morphisms $E_1\la E_2\la E_3\la\cdots$
\be
\item
  Each $E_m$ belongs to $\gen Gn$.
\item
The map $\ch^i(E_m)\la\ch^i(E_{m+1})$ is an isomorphism whenever $i\geq-m$.
\setcounter{enumiv}{\value{enumi}}
\ee
In this definition we also allow $n=\infty$, we
simply declare $\gen G\infty=\gen G{}=\ct^c$.

Suppose we are also given an object $F\in\ct$, together with
\be
\setcounter{enumi}{\value{enumiv}}
\item
A map of
the approximating system $E_*$ to $F$.
\item
The map in (iii) is such that
$\ch^i(E_m)\la\ch^i(F)$ is an isomorphism whenever $i\geq-m$.
\ee 
Then we declare $E_*$ to be a
\emph{strong $\gen Gn$--approximating system for $F$.}
\edfn

\rmk{R1.31}
Although Definition~\ref{D1.29} was phrased
in terms of the particular choice
of \tstr\ 
$(\ct^{\leq0},\ct^{\geq0})$, made in Notation~\ref{N37.9993},
it is robust---up to passing to subsequences
a
strong $\gen Gn$--approximating system for $F$ will work for
any equivalent \tstr.
\ermk

\lem{L1.1010101}
With the conventions of Definition~\ref{D1.29} we have
\be
\item
Given an object $F\in\ct^-$ and a strong $\gen Gn$--approximating
system $E_*$ for $F$, then the (non-canonical) map $\hoco E_i\la F$
is an isomorphism.
\item
Any object $F\in\ct^-_c$ has a strong $\ct^c$--approximating
system.
\item
Any strong $\gen Gn$--approximating system
$E_1\la E_2\la E_3\la\cdots$ is a
strong $\gen Gn$--approximating system of the
homotopy colimit $F=\hoco E_i$. Moreover
$F$ belongs to $\ct^-_c$.
\ee 
\elem

\prf
We begin by proving (iii).
Suppose $E_1\la E_2\la E_3\la\cdots$ is a strong $\gen Gn$--approximating
system, and let $F=\hoco E_i$. The objects $E_i$ all belong to
$\gen Gn\subset\ct^c\subset\ct^-$. Choose an integer $\ell>0$ with
$E_1\in\ct^{\leq \ell}$. The fact that $\ch^i(E_1)\la\ch^i(E_m)$ is an isomorphism
for all $i\geq-1$ means that $\ch^i(E_m)=0$ for all $i>\ell$ and all $m$,
and Lemma~\ref{L-1.3} gives that the
$E_m$ all lie in $\ct^{\leq \ell}$. Hence the homotopy
colimit $F$ also belongs to $\ct^{\leq \ell}$.

Now Remark~\ref{R-1.9} tells us that the map $\colim\,\ch^i(E_m)\la\ch^i(F)$
is an isomorphism for every $i\in\zz$. By the previous paragraph
the triangle
$E_m\la F\la D_m$ lies in $\ct^-$, and
as $\ch^i(E_m)\la\ch^i(F)$ is an isomorphism for $i\geq-m$
we deduce that $\ch^i(D_m)=0$ for all $i\geq-m$.
Lemma~\ref{L-1.3} guarantees that $D_m\in\ct^{\leq-m-1}$, and as
$E_m\in\gen Gn\subset\ct^c$ and $m>0$ is arbitrary
we have that $F$ satisfies the criterion for belonging to $\ct^-_c$.

Next we prove (i). By (iii) the map $\hoco E_i\la F$ is a morphism from
$\hoco E_i\in\ct^-_c$ to $F\in\ct^-$, hence it is a morphism
in $\ct^-$. The hypothesis of (i), coupled with Remark~\ref{R-1.9},
tell us that $\ch^i\big(\hoco E_i\big)\la \ch^i(F)$ is an isomorphism
for every $i\in\zz$. By Lemma~\ref{L-1.5}(ii) the map
$\hoco E_i\la F$ is an isomorphism.

It remains to prove (ii). Note that, if we assume more approximability
hypotheses on $\ct$, then (ii) is immediate from Corollary~\ref{C17.29}.
But let us see that we don't yet need any strong assumptions.

Take any $F\in\ct^-_c$. There exists a triangle $E_1\la F\la D_1$ with
$E_1\in\ct^c$ and $D_1\in\ct^{\leq-3}$. When $i\geq-1$ exact sequence
$\ch^{i-1}(D_1)\la\ch^i(E_1)\la\ch^i(F)\la\ch^i(D_1)$ has
$\ch^{i-1}(D_1)=0=\ch^i(D_1)$, starting the construction of $E_*$.

Suppose now that we have constructed the sequence up to an integer
$n>0$, that is we have a map $f_m:E_m\la F$, with $E_m\in\ct^c$, and so
that $\ch^i(f_m)$ is an isomorphism for all $i\geq-m$. Lemma~\ref{L17.17}
allows us to choose an integer $N>0$ so that
$\Hom\big(E_m,\ct^{\leq-N}\big)=0$. Because $F$ belongs to $\ct^-_c$ we may
choose a triangle $E_{m+1}\la F\la D_{m+1}$ with $E_{m+1}\in\ct^c$
and $D_{m+1}\in\ct^{\leq-N-m-3}$. As in the paragraph above we
show that the map $\ch^i(E_{m+1})\la\ch^i(F)$ is an isomorphism for all
$i\geq-m-1$. And since the composite $E_m\stackrel{f_n}\la F\la D_{m+1}$
vanishes, the map $f_n$ must factor as $E_m\la E_{m+1}\la F$.
\eprf

\rmk{R37.20202025}
Lemma~\ref{L1.1010101}(i) and Remark~\ref{R27.99309} combine to tell
us that a strong $\gen Gn$--approximating system for $F$,
in the sense of Definition~\ref{D1.29},
is in fact an approximating system for $F$
as defined in Remark~\ref{R27.7}. Our terminology isn't misleading.
\ermk

\rmk{R37.2020202}
Let us now specialize Lemma~\ref{L1.35} to the framework of this section.
Assume we are given
\be
\item
A morphism $\wh\alpha:\wh A\la \wh B$ in the category $\ct^-_c$.
\item
Two integers $n'$ and $n$, as well as a strong $\gen G{n'}$--approximating
system $\mathfrak{A}_*$ for $\wh A$ and a
strong $\gen G{n}$--approximating
system $\mathfrak{B}_*$ for $\wh B$.
\setcounter{enumiv}{\value{enumi}}
\ee
Lemma~\ref{L27.592} allows us to choose a subsequence of
$\mathfrak{B'}_*\subset\mathfrak{B}_*$ and a map of sequences
$\alpha_*:\mathfrak{A}_*\la\mathfrak{B'}_*$ compatible with
$\wh\alpha:\wh A\la \wh B$.
A subsequence of a strong $\gen G{n}$--approximating sequence is
clearly a strong $\gen G{n}$--approximating sequence,
hence $\mathfrak{B'}_*$ is a strong $\gen G{n}$--approximating sequence for
$\wh B$.
Now 
as in Lemma~\ref{L1.35} we extend
$\alpha_*:\mathfrak{A}_*\la\mathfrak{B'}_*$ to a sequence of
triangles, in particular
for each $m>0$ this gives a morphism of triangles
\[\xymatrix@C+20pt{
  \mathfrak{A}_m\ar[r]^{\alpha_m}\ar[d] & \mathfrak{B'}_m\ar[d]\ar[r]^{\beta_m}
    & \mathfrak{C}_m\ar[r]^{\gamma_m}\ar[d]&\T \mathfrak{A}_m\ar[d]
\ar[r]^-{\T\alpha_m} &\T \mathfrak{B'}_m\ar[d]\\
  \mathfrak{A}_{m+1}\ar[r]^{\alpha_{m+1}}& \mathfrak{B'}_{m+1}\ar[r]^{\beta_{m+1}}&
  \mathfrak{C}_{m+1}\ar[r]^{\gamma_{m+1}} &\T \mathfrak{A}_{m+1}
\ar[r]^-{\T\alpha_{m+1}} &\T \mathfrak{B'}_{m+1}
}\]
Applying the functor $\ch^i$ with $i\geq-m$ yields a
commutative diagram in the heart of $\ct$ where the rows are exact, and where
the vertical maps away from the middle are isomorphisms.
By the 5-lemma the middle vertical map, i.e.~the map
$\ch^i(\mathfrak{C}_{m})\la\ch^i(\mathfrak{C}_{m+1})$, must also be an
isomorphism when $i\geq-m$.
We conclude that $\mathfrak{C}_*$ is a strong
$\gen G{n'+n}$--approximating system. Put $\wh C=\hoco\mathfrak{C}_*$.
By Lemma~\ref{L1.1010101}(iii) the object $\wh C$ belongs to $\ct^-_c$
and $\mathfrak{C}_*$ is a strong
$\gen G{n'+n}$--approximating system for $\wh C$,
while Remark~\ref{R997.207} guarantees that the weak triangle
$A\stackrel u\la B\stackrel v\la C\stackrel w\la\T A$
of Lemma~\ref{L1.35} is isomorphic
to the image under $\cy$ of a weak triangle
$\wh A\stackrel{\wh u}\la \wh B\stackrel{\wh v}\la\wh C\stackrel{\wh w}\la\T \wh A$
in the category $\ct^-_c$.

Furthermore: the homological functor $\ch$ takes each
of the triangles
$\mathfrak{A}_m\stackrel{\alpha_m}\la\mathfrak{B'}_m\stackrel{\beta_m}\la\mathfrak{C}_m\stackrel{\gamma_m}\la\T\mathfrak{A}_m$
to a long exact sequence, and by Remark~\ref{R-1.9} the
(eventually stable) colimit is
$\ch$ of the weak
triangle
$\wh A\stackrel{\wh u}\la \wh B\stackrel{\wh v}\la\wh C\stackrel{\wh w}\la\T \wh A$.
Hence $\ch$ takes the weak triangle
$\wh A\stackrel{\wh u}\la \wh B\stackrel{\wh v}\la\wh C\stackrel{\wh w}\la\T \wh A$
to a long exact sequence.
\ermk

\lem{L37.601}
Let the conventions be as in Notation~\ref{N37.9993}. Assume
$H$ is a locally finite $\gen Gn$--cohomological functor.
Then there exists an object
$F\in\ct^-_c$ and an epimorphism of $\gen Gn$--cohomological
functors $\ph:\cy(F)|_{\gen Gn}\la H$. Furthermore the object
$F$ may be chosen to have a strong $\gen Gn$--approximating
system.

We will in fact prove a refinement of the above. Since $H$ is
assumed to be a locally finite $\gen Gn$--cohomological functor its
restriction to $\gen Gm$, for
any integer $m<n$, is a locally finite 
$\gen Gm$--cohomological
functor. Hence for any $m<n$ the first paragraph
delivers an object $F_m\in\ct^-_c$, with a
strong $\gen Gm$--approximating system, and a surjective natural
transformation $\ph_m:\cy(F_m)|_{\gen Gm}\la H|_{\gen Gm}$.
We will actually construct these $F_m$'s compatibly. We will produce
in $\ct^-_c$ a sequence
$F_1\la F_2\la \cdots\la F_{n-1}\la F_n$,
with compatible maps $\ph_m:\cy(F_m)|_{\gen Gn}\la H$, so that
\be
\item
  For each $m>0$ the object $F_m$ has a strong $\gen Gm$--approximating
  system, and
  the map $\ph_m|_{\gen Gm}:\cy(F_m)|_{\gen Gm}\la H|_{\gen Gm}$ is an epimorphism.
\item
  The sequence is such that the kernel
  of the map $(\ph_m)|_{\gen G1}:\cy(F_m)|_{\gen G1}\la H|_{\gen G1}$  
  is annihilated by the map $\cy(F_m)|_{\gen G1}\la\cy(F_{m+1})|_{\gen G1}$.
\setcounter{enumiv}{\value{enumi}}
\ee
\elem

\prf
The proof is by induction on $n$. In the case $n=1$ we prove the 
refinement that allows the induction to proceed
\be
\setcounter{enumi}{\value{enumiv}}
\item
  Suppose $H$ is a locally finite $\gen G1$--cohomological functor.
  Then we may construct an
  object $F\in\ct^-_c$ and an epimorphism $\cy(F)|_{\gen G1}\la H$.
  Furthermore
  the object $F$ can be chosen to have a strong $\gen G1$--approximating
  system $E_1\la E_2\la E_3\la\cdots$ \emph{in which every morphism
  $E_i\la E_{i+1}$ is a split monomorphism.}
\setcounter{enumiv}{\value{enumi}}
\ee
The proof of (iii) is easy: we have that $H(\T^i G)$ is a finite $R$--module
for every
$i\in\zz$, and vanishes if $i\ll0$. For each $i$ with $H(\T^iG)\neq0$
choose a
finite number of generators $\{f_{ij},\,j\in J_i\}$
for the $R$--module $H(\T^iG)$. By
Yoneda
every $f_{ij}\in H(\T^iG)$ corresponds
to a morphism $\ph_{ij}^{}:\cy(\T^i G)|_{\gen G1}\la H$. Let $F$ be
defined by
\[
F =\coprod_{i\in\zz}\bigoplus_{j\in J_i}\T^i G
\]
and let the morphism $\ph:\cy(F)|_{\gen G1}\la H$
be given by
\[\xymatrix@C+20pt{
  \cy(F)|_{\gen G1}\ar@{=}[r]&
  \ds\bigoplus_{i\in\zz}\bigoplus_{j\in J_i}\cy(\T^iG)|_{\gen G1}
  \ar[rr]^-{(\ph_{ij})} && H
}\]
where $(\ph_{ij})$ stands for the row matrix with entries $\ph_{ij}$; on
the $i,j$ summand the map is $\ph_{ij}$.
Finally: because the \tstr\ is in the preferred
equivalence class there is an integer $B>0$ with
$\T^BG\in\ct^{\leq0}$. For $m>0$ we define
\[
E_m=\bigoplus_{i\leq m+B}\,\bigoplus_{j\in J_i}\T^i G
\]
The sum is finite by hypothesis, making
$E_m$ an object of $\gen G1$. The obvious map $E_m\la E_{m+1}$
is a split monomorphism, and in the decomposition $F\cong E_m\oplus \wt F$ we
have that $\wt F$, being the coproduct of $\T^iG$ for $i\geq m+B+1$,
belongs to $\ct^{\leq-m-1}$. Therefore the map $\ch^i(E_m)\la\ch^i(F)$
is an isomorphism when $i\geq-m$,
making $E_*$ is a strong $\gen G1$--approximating
system for $F$.

Now for the induction step. Suppose $n\geq 1$ is an integer and we know the
Lemma for all integers $\leq n$. We wish to show it holds for $n+1$. Let
$H$ be a locally finite $\gen G{n+1}$--cohomological functor.
Then the restriction of $H$ to $\gen Gn$ is a locally finite
$\gen Gn$--cohomological functor, and we may apply the induction hypothesis
to produce in $\ct^-_c$ a sequence
$F_1\la F_2\la \cdots\la F_{n-1}\la F_n$,
with compatible
surjections $\ph_m:\cy(F_m)|_{\gen Gm}\la H|_{\gen Gm}$.
In particular
the map $\ph_n:\cy(F_n)|_{\gen Gn}\la H|_{\gen Gn}$
is an epimorphism.
Complete the natural transformation $\ph_n$ to a short exact sequence
\[\xymatrix@C+20pt{
  0\ar[r] & H'\ar[r] &\cy(F_n)|_{\gen Gn}\ar[r]^-{\ph_n} & H|_{\gen Gn}
  \ar[r] & 0
}\]
of functors on $\gen Gn$. Since $H|_{\gen Gn}$ and $\cy(F_n)|_{\gen Gn}$
are
locally finite $\gen Gn$--cohomological functors so
is $H'$, and induction
applies. For any $m\leq n$ we may choose a surjection
$\ph'_m:\cy(F')|_{\gen Gm}\la H'|_{\gen Gm}$,
as in (i) and (ii). We wish to consider
the special case $m=1$, where we can assume our $F'$ is as in (iii).
That is we choose an object $F'\in\ct^-_c$, which admits a strong
$\gen G1$--approximating system and a surjection
$\ph':\cy(F')|_{\gen G1}\la H'|_{\gen G1}$.
 And
we may further assume that our approximating system
$E'_1\la E'_2\la E'_3\la\cdots$ for $F'$ is such that every morphism
$E'_i\la E'_{i+1}$ is a split monomorphism.

We have natural transformations
$\cy(F')|_{\gen G1}\la H'|_{\gen G1}\la\cy(F_n)|_{\gen G1}$,
and as $F'$ admits a
$\gen G1$--approximating system Lemma~\ref{L27.976} tells
us that the composite
is equal to $\cy(\alpha_n)|_{\gen G1}$ for
some morphism $\alpha_n:F'\la F_n$
in the category $\ct$. And now Lemma~\ref{L1.35} applies;
see Remark~\ref{R37.2020202} for an elaboration of how it specializes
to the current context.
We learn that
\be
\setcounter{enumi}{\value{enumiv}}
\item
  There exists a weak triangle
  $F'\stackrel{\alpha_n}\la F_n\stackrel{\beta_n}\la F_{n+1}$ in
  the category $\ct^-_c$,
with $F_{n+1}$ admitting a strong $\gen G{n+1}$--approximating
system.
\item
  There is natural transformation
  $\ph_{n+1}:\cy(F_{n+1})|_{\gen G{m+1}}\la H$,
  such that $\ph_n:\cy(F_n)|_{\gen G{n}}\la H|_{\gen G{n}}$ is equal to
  the composite
\[\xymatrix@C+40pt{
  \cy(F_n)|_{\gen Gn}\ar[r]^-{\cy(\beta_n)|_{\gen Gn}} & \cy(F_{n+1})|_{\gen Gn} \ar[r]^-{\ph_{n+1}|_{\gen Gn}} &
  H|_{\gen Gn}\ .
}\]
\setcounter{enumiv}{\value{enumi}}
\ee
Comparing the two exact sequences
\[\xymatrix@C+40pt@R-15pt{
  \cy(F')|_{\gen G1}\ar[r]^-{\cy(\alpha_n)|_{\gen G1}} & \cy(F_{n})|_{\gen G1} \ar[r]^-{\ph_{n}|_{\gen G1}} &
  H|_{\gen G1} \\
 \cy(F')|_{\gen G1}\ar[r]^-{\cy(\alpha_n)|_{\gen G1}} & \cy(F_{n})|_{\gen G1} \ar[r]^-{\cy(\beta_{n})|_{\gen G1}} &
  \cy(F_{n+1})|_{\gen G1}
}\]
we conclude that $\ph_{n}|_{\gen G1}$ and $\cy(\beta_{n})|_{\gen G1}$ have
the same kernel, that is
the map $\beta_n:F_n\la F_{n+1}$ satisfies (ii).

To finish the proof of
(i) it remains to show that $\ph_{n+1}:\cy(F_{n+1})|_{\gen G{n+1}}\la H$ is
an epimorphism, but this is now immediate from
Lemma~\ref{L997.601}.
\eprf

\rmk{R37.1003}
Let the conventions be as in Notation~\ref{N37.9993} and assume
$H$ is a locally finite $\ct^c$--cohomological functor.
For every integer
$n>0$ the restriction of $H$ to $\gen Gn$ is a locally finite
$\gen Gn$--homological functor,
and Lemma~\ref{L37.601} permits us to construct a
sequence $F_1\stackrel{\beta_1^{}}\la F_2
\stackrel{\beta_2^{}}\la F_3\stackrel{\beta_3^{}}\la\cdots$ in the
category $\ct^-_c$, together with compatible
epimorphisms $\ph_n:\cy(F_n)|_{\gen Gn}\la H|_{\gen Gn}$.
Since $F_n$ is constructed to have a $\gen Gn$--approximating
system, Corollary~\ref{C27.5} says that each $\ph_n$ lifts
uniquely to a natural transformation defined on
all of $\ct^c$, which we denote 
$\ph_n:\cy(F_n)\la H$. Note
the abuse of notation, where $\ph_n$ is allowed to stand
both for the natural transformation of functors on $\ct^c$
and for its restriction to $\gen Gn$. The triangle
\[\xymatrix@C+40pt@R-20pt{
    & \cy(F_{n+1}) \ar[dd]^{\ph_{n+1}}\\
  \cy(F_{n}) \ar[rd]_-{\ph_{n}}
  \ar[ur]^-{\cy(\beta_{n}^{})} & \\
  & H
}\]
commutes when restricted to the subcategory $\gen Gn\subset\gen G{}$,
and the fact that $F_n$ has a $\gen Gn$--approximating system
coupled with the uniqueness assertion of Corollary~\ref{C27.5}
tells us that the triangle commutes on the nose, in
the category $\Hom\big([\ct^c]\op,\Mod R\big)$.
\ermk

\pro{P37.1005}
Let the conventions be as in Notation~\ref{N37.9993} and assume
$H$ is a locally finite $\ct^c$--cohomological functor.
Then there exists an object $F\in\ct$ and an isomorphism
$\ph:\cy(F)\la H$.

Now let $\cg=\oplus_{C\in\ct^c}^{}C$; for those worried about set theoretic
issues this means that $\cg$ is the coproduct, over the isomorphism
classes of objects in $\ct^c$, of a representative in each
isomorphism class. Then $F$ may be chosen to lie in $\ogenun \cg4$.
\epro

\prf
In Remark~\ref{R37.1003} we noted
that 
Lemma~\ref{L37.601} constructs for us a
sequence
$F_1\stackrel{\beta_1^{}}\la F_2
\stackrel{\beta_2^{}}\la F_3\stackrel{\beta_3^{}}\la\cdots$
in the
category $\ct$, together with compatible
maps $\ph_n:\cy(F_n)\la H$;
there is an induced map
$\colim\,\cy(F_n)\la H$.
If we define $F$ to be  $F=\hoco F_n$ then we have an
object $F\in\ct$, and 
\cite[Lemma~2.8]{Neeman96} tells
us that the natural map
$\colim\,\cy(F_n)\la\cy(F)$
is an isomorphism. We have constructed a map
$\ph:\cy(F)\la H$ and will prove that
$\ph$ is an isomorphism.

Let us consider the restriction of the natural transformation
$\ph$ to the subcategory $\gen G1\subset\ct^c$.
The natural transformation $\ph|_{\gen G1}:\cy(F)|_{\gen G1}\la H|_{\gen G1}$
is the map to $H|_{\gen G1}$ from the colimit of the sequence
\[\xymatrix@C+40pt{
  \cy(F_1)|_{\gen G1}\ar[r]^-{\cy(\beta_1^{})|_{\gen G1}}
  &\cy(F_2)|_{\gen G1}\ar[r]^-{\cy(\beta_2^{})|_{\gen G1}}
 &\cdots
}\]
and Lemma~\ref{L37.601}(ii) says the sequence is
such that each map 
$\cy(\beta_n^{})|_{\gen G1}$ factors
as $\cy(F_n)|_{\gen G1}\la H|_{\gen G1}\la \cy(F_{n+1})|_{\gen G1}$.
Hence the colimit agrees with the colimit
of the ind-isomorphic
constant sequence $H|_{\gen G1}\la H|_{\gen G1}\la H|_{\gen G1}\la $,
and this proves that the restriction of $\ph:\cy(F)\la H$
to the category $\gen G1$ is an isomorphism. Concretely: for every
$i\in\zz$ the map $\ph:\Hom(\T^iG,F)\la H(\T^iG)$
is an isomorphism.
The full subcategory $\ck\subset\ct^c$ defined by
\[
\ck\eq\left\{
K\in\ct^c\,\left|\begin{array}{c}
\forall i\in\zz\text{\rm\ the map}\\
\ph:\Hom(\T^iK,F)\la H(\T^iK)\\
\text{\rm is an isomorphism}
\end{array}
\right.\right\}
\]
is thick and contains $G$, hence $\ck\subset\ct^c=\gen G{}\subset\ck$.

It remains to prove that the $F$ we constructed belongs to $\ogenun \cg4$.
We begin with the observation that each $F_n$ in the sequence
$F_1\stackrel{\beta_1^{}}\la F_2
\stackrel{\beta_2^{}}\la F_3\stackrel{\beta_3^{}}\la\cdots$
of Remark~\ref{R37.1003}  has a strong $\gen Gn$--approximating system.
This means that $F_n=\hoco E_i^n$, with each $E_i^n\in\gen Gn$.
The triangle
\[\xymatrix@C+30pt{
  \ds\bigoplus_{i=1}^\infty E_i^n \ar[r] &
  \ds\bigoplus_{i=1}^\infty E_i^n \ar[r] &
  F_n\ar[r] &
  \ds\bigoplus_{i=1}^\infty\T E_i^n  
}\]
tells us that $F_n$ must belong to
$\ogenun\cg1*\ogenun\cg1\subset\ogenun \cg 2$.
But now $F=\hoco F_n$, and the triangle
\[\xymatrix@C+30pt{
  \ds\bigoplus_{n=1}^\infty F_n \ar[r] &
  \ds\bigoplus_{n=1}^\infty F_n \ar[r] &
  F\ar[r] &
  \ds\bigoplus_{n=1}^\infty\T F_n  
}\]
gives that $F$ belongs to $\ogenun\cg2*\ogenun\cg2\subset\ogenun\cg4$.
\eprf

\ntn{N37.10245}
This is as far as we get with the
assumptions of Notation~\ref{N37.9993}.
From now on we will assume further
that $\ct$ is weakly approximable.
\entn

\lem{L37.1025}
Let the conventions be as in Notation~\ref{N37.10245}. 
Assume
$H$ is a locally finite $\ct^c$--cohomological functor.
There exists an integer $\wt A>0$ and, for any
$n>0$, an object $F_n\in\ct^-_c\cap\ct^{\leq\wt A}$
as well as a natural transformation
$\ph_n:\cy(F_n)\la H$
which is surjective when restricted to $\gen Gn$.
\elem

\prf
Remark~\ref{R37.1003} produced for us an object $F_n\in\ct^-_c$
and
a natural transformation $\ph_n:\cy(F_n)\la H$,
so that the restriction to $\gen Gn$ of $\ph_n$ is surjective.
The new assertion is that we may choose $F_n$ to lie
in $\ct^{\leq\wt A}$ for some $\wt A>0$ independent of $n$.

Because $H$ is locally finite there exists an integer
$B'>0$ with $H(\T^{-i}G)=0$ for all $i\geq B'$. Since 
$H$ is cohomological it follows that $H(E)=0$ for all
$E\in\genuf G{}{B'}$. And Corollary~\ref{C17.29} gives us an
integer $B>0$ so that, for every integer $m>0$, every object
$F\in\ct^-_c\cap\ct^{\leq0}$ admits a triangle
$E_m\la F\la D_m$ with $D_m\in\ct^{\leq-m}$ and
$E_m\in\genu G{}{1-m-B,B}\subset\genuf G{}{1-m-B}$. I assert that
$\wt A=B+B'$ works.

Let us begin with the $F_n$ provided by
Remark~\ref{R37.1003}. Because it belongs to $\ct^-_c\subset\ct^-$
there exists an integer $\ell>B+B'$ with $F_n\in\ct^{\leq\ell}$.
Applying Corollary~\ref{C17.29} to the object
$\T^\ell F_n\in\ct^-_c\cap\ct^{\leq0}$, with $m=\ell-B-B'$,
we learn that there exists a triangle $E\stackrel\alpha\la F_n\la D$ with
$D\in\ct^{\leq \ell-m}=\ct^{\leq B+B'}=\ct^{\leq\wt A}$ and
$E\in\genuf G{}{\ell+1-m-B}=\genuf G{}{1+B'}$.
In particular $E\in\ct^c$, and $H(E)=0$ by the choice of
$B'$. Hence
the composite 
\[\xymatrix@C+60pt{
  \cy(E)\ar[r]^-{\cy(\alpha)}  
  &\cy(F_n)\ar[r]^-{\ph} 
 & H
}\]
vanishes.

Now we apply Lemma~\ref{L1.35} with $\ca=\cb=\ct^c$. The object
$F_n\in\ct^-_c$, of Lemma~\ref{L37.601} and Remark~\ref{R37.1003},
comes with a strong $\gen Gn$--approximating system, which
is certainly a strong $\ct^c$--approximating system.
The object $E\in\ct^c$ comes with the trivial strong
$\ct^c$--approximating system
$E\stackrel\id\la E\stackrel\id\la E\stackrel\id\la \cdots$. In this
system the connecting maps are all identities, which are
split monomorphisms. The hypotheses of Lemma~\ref{L1.35}
and Remark~\ref{R37.2020202} hold and
the Lemma produces for us, in $\ct^-_c$, a weak triangle 
$E\stackrel\alpha\la F_n\stackrel\beta\la \wt D\la \T E$
and a factorization of $\ph:\cy(F_n)\la H$
as a composite 
\[\xymatrix@C+60pt{
  \cy(F_n)\ar[r]^{\cy(\beta)}
  &\cy(\wt D)\ar[r]^-{\psi} 
 &H\ .
}\]
The surjectivity of the restriction to $\gen Gn$ of $\ph$ implies the
surjectivity of the restriction to $\gen Gn$ of $\psi$.

It remains to show that $\wt D\in\ct^-_c$ belongs to
$\ct^-_c\cap\ct^{\leq \wt A}$. We know that, in the triangle
$E\la F_n\la D\la \T E$, the object $D$ belongs to $\ct^{\leq\wt A}$.
The long exact
sequence $\ch^{i-1}(D)\la\ch^i(E)\la\ch^i(F_n)\la\ch^i(D)$
tells us that $\ch^i(\alpha):\ch^i(E)\la\ch^i(F)$
is surjective if $i=1+\wt A$ and is an isomorphism
when $i>1+\wt A$. The long exact sequence
$\ch^i(E)\la\ch^i(F_n)\la\ch^i(\wt D)\la\ch^{i+1}(E)\la\ch^{i+1}(F_n)$
says that $\ch^i(\wt D)=0$ if $i\geq 1+\wt A$. By Lemma~\ref{L-1.3}
we conclude that $\wt D\in\ct^{\leq\wt A}$.
\eprf

\ntn{N37.1024}
This is as far as we get with the
assumptions of Notation~\ref{N37.10245},
from now on we assume further
that $\ct$ is approximable---weak approximability will no
longer be enough.
\entn

\lem{L37.1036}
Let the conventions be as in Notation~\ref{N37.1024}, 
and assume
$H$ is a locally finite $\ct^c$--cohomological functor.
Choose an integer $B>0$ as in
Lemma~\ref{L17.25}.

Suppose $\wt F,F'$ are objects in $\ct^-_c\cap\ct^{\leq0}$,
$E$ is an object in $\ct^c\cap\ct^{\leq0}$, and we have a morphism
$\alpha:E\la \wt F$ in $\ct^-_c$.
Assume we are given an integer $m>0$, as
well as natural transformations
$\wt\ph:\cy(\wt F)\la H$ and
$\ph':\cy(F')\la H$,
so that $\wt\ph$ restricts to an epimorphism on $\gen G{mB}$
and $\ph'$ restricts to an epimorphism on $\gen G{(m+1)B}$.

Then there exists in $\ct^-_c\cap\ct^{\leq0}$ a commutative diagram
\[\xymatrix{
E\ar[rr]^{\e}\ar[dr]_{\alpha}&  &\ar[dl]^{\gamma} E'\ar[dr]^{\alpha'} & & F'\ar[dl]_-{\beta} \\
& \wt F &                 & \wt F'
}\]
and there exists a natural transformation
$\wt\ph':\cy(\wt F')\la H$
so that
\be
\item
  The object $E'$ belongs to $\ct^c\subset\ct^-_c$.
\item
  The maps $\ch^i(\alpha')$ and $\ch^i(\gamma)$ are
  isomorphisms for all $i\geq-m+2$. 
\item
  The triangle 
\[\xymatrix@C+50pt@R-20pt{
    & \cy(\wt F') \ar[dd]^{\wt\ph'}\\
  \cy(F') \ar[rd]_-{\ph'}
  \ar[ur]^-{\cy(\beta)} & \\
  & H
}\]
  commutes; the surjectivity of the restriction to $\gen G{(m+1)B}$
  of $\ph'$ therefore implies the surjectivity of the
  restriction to $\gen G{(m+1)B}$
  of $\wt\ph'$.
\item
  The square below commutes
  \[\xymatrix@C+50pt{
    \cy(E') \ar[r]^-{\cy(\gamma)}
    \ar[d]^-{\cy(\alpha')}& \cy(\wt F) \ar[d]^{\wt\ph}\\
  \cy(\wt F') \ar[r]_-{\wt\ph'}  & H
}\]
\setcounter{enumiv}{\value{enumi}}
\ee
\elem

\prf
Corollary~\ref{C17.29}, applied to the object $F'\in\ct^-_c\cap\ct^{\leq0}$
and the integer $m$, permits us to construct in $\ct^-_c$ a triangle
$E'_m\stackrel a\la F'\stackrel b\la D_m$ with $D_m\in\ct^{\leq-m}$ and
$E'_m\in\genu G{mB}{1-m-B,B}\subset\gen G{mB}$.
Because $E'_m$ belongs to $\ct^c$ the
composite 
\[\xymatrix@C+40pt{
\cy(E'_m)\ar[r]^-{\cy(a)} &
  \cy(F')\ar[r]^-{\ph'} & H
}\]
is a natural transformation from a representable functor
on $\ct^c$ to $H$, and corresponds to an element $x\in H(E'_m)$.
As $E'_m$ belongs to $\gen G{mB}$ 
the morphism $\wt\ph:\Hom(E'_m,\wt F)\la H(E'_m)$ is
surjective; there is a map $f:E'_m\la \wt F$ with
$\wt\ph(f)=x$. Yoneda translates this to mean that the 
square below commutes
  \[\xymatrix@C+50pt{
    \cy(E'_m) \ar[r]^-{\cy(f)}
    \ar[d]^-{\cy(a)}& \cy(\wt F) \ar[d]^{\wt\ph}\\
  \cy(F') \ar[r]_-{\ph'}  & H
}\]

We now have a morphism $(\alpha,f):E\oplus E'_m\la \wt F$.
The object $E\oplus E'_m$ belongs to $\ct^c$ while
the object $\wt F\in\ct^-_c$ has a strong $\ct^c$--approximating
system---see Lemma~\ref{L1.1010101}(ii). There exists a sequence
$\wt E_1\la \wt E_2\la\wt E_3\la\cdots$ in $\ct^c$,
and a map of $\wt E_*$ to $\wt F$, and so that
$\ch^i(\wt E_j)\la\ch^i(\wt F)$ is an isomorphism whenever
$i\geq-j$.
The morphism $E\oplus E'_m\la \wt F$ must factor
through some $\wt E_j$---see Lemma~\ref{L1.1010101}(i) and
\cite[Lemma~2.8]{Neeman96}.
Choose such an $\wt E_j$ with $j\geq m$, declare $E'=\wt E_j$, and
let $\e:E\la E'$, $g:E'_m\la E'$ and $\gamma:E'\la \wt F$ be
the obvious maps. Then $\gamma\e=\alpha$ and $\gamma g=f$.

By construction the map $\ch^\ell(\gamma):\ch^\ell(E')\la\ch^\ell(\wt F)$ is an
isomorphism for all $\ell\geq-m$.
Recalling 
the triangle $E'_m\stackrel a\la F'\la D'_m$ with 
$D'_m\in\ct^{\leq-m}$, the exact sequence
\[\xymatrix@R-20pt{
\ch^{\ell-1}(D'_m)\ar[r]  & \ch^\ell(E'_m)\ar[r]^{\ch^\ell(a)} &
\ch^\ell(F')\ar[r] & \ch^\ell(D'_m)
}\]
teaches us that
\be
\setcounter{enumi}{\value{enumiv}}
\item
The maps $\ch^\ell(\gamma)$ and $\ch^\ell(a)$
are isomorphisms for all $\ell\geq-m+2$.
\setcounter{enumiv}{\value{enumi}}
\ee
Since $m\geq1$ by
assumption, we learn in particular that
if $\ell\geq1$ then $\ch^\ell(E')\cong\ch^\ell(\wt F)=0$
and $\ch^\ell(E'_m)\cong\ch^\ell(F')=0$. As both
$E'$ and $E'_m$ are objects of $\ct^c\subset\ct^-$,
Lemma~\ref{L-1.3} informs us that
\be
\setcounter{enumi}{\value{enumiv}}
\item
$E',E'_m$ both
lie in $\ct^{\leq0}$.
\setcounter{enumiv}{\value{enumi}}
\ee

Now consider the commutative square
  \[\xymatrix@C+50pt{
    \cy(E'_m) \ar[r]^-{\cy(g)}
    \ar[d]^-{\cy(a)}& \cy(E') \ar[d]^{\rho=\wt\ph\circ\cy(\gamma)}\\
  \cy(F') \ar[r]_-{\ph'}  & H
}\]
In other words: we have in $\ct^-_c\cup\ct^{\leq0}$ a morphism 
$\s=\left(\begin{array}{r}
-g\\ a
\end{array}\right)\colon E'_m \la E'\oplus F'$,
as well as a natural transformation $(\rho,\ph'):\cy(E'\oplus
F')\la H$, and the composite
 \[\xymatrix@C+50pt{
    \cy(E'_m) \ar[r]^-{\cy(\s)}
    & \cy(E'\oplus F') \ar[r]^-{(\rho,\ph')} &
  H
}\]
 vanishes. The object $E'\oplus F'\in\ct^-_c$ has a strong
 $\ct^c$--approximating system by Lemma~\ref{L1.1010101}(ii),
 and the object $E'_m\in\ct^c$ has the trivial strong
 $\ct^c$--approximating system
$E'_m\stackrel\id\la E'_m\stackrel\id\la E'_m\stackrel\id\la\cdots$.
We may apply Lemma~\ref{L1.35} as specialized in Remark~\ref{R37.2020202},
with $\ca=\cb=\ct^c$, to
deduce that $E'_m\stackrel\s\la E'\oplus F'$ may be completed in $\ct^-_c$
to a 
weak triangle
\[\xymatrix@C+20pt{
  E'_m \ar[r]^-{\s} & E'\oplus F' \ar[r]^-{(\alpha',\beta)} &\wt F'
  \ar[r]^-{\tau} &\T E'_m
}\]
in such a way that the morphism $(\rho,\ph'):\Hom(-,E'\oplus
F')|_{\ct^c}\la H(-)$ factors as
$\wt\ph'\circ\Hom\big(-,(\alpha',\beta)\big)$ for some natural
transformation
$\wt\ph':\Hom(-,\wt F')|_{\ct^c}\la H(-)$.
We have constructed in $\ct^-_c$ a commutative square
\[\xymatrix@C+20pt{
  E'_m \ar[r]^-{g}\ar[d]^a & E'\ar[d]^{\alpha'} \\
F'\ar[r]^\beta & \wt F'
}\] 
and a natural transformation
$\wt\ph':\Hom(-,\wt F')|_{\ct^c}\la H(-)$
so that the diagram below commutes
\[\xymatrix@C+50pt@R+10pt{
    \cy(E'_m) \ar[r]^-{\cy(g)}
    \ar[d]^-{\cy(a)}& \cy(E') \ar[r]^{\cy(\gamma)}
    \ar[d]^-{\cy(\alpha')}
    & \cy(\wt F) \ar[dd]^{\wt\ph}\\
    \cy(F') \ar@/_2pc/[rrd]_-{\ph'}\ar[r]^-{\cy(\beta)}
    & \cy(\wt F') \ar[rd]^-{\wt\ph'}  &  \\
 & & H
}\]
This finishes our construction of the diagram 
\[\xymatrix{
E\ar[rr]^{\e}\ar[dr]_{\alpha}&  &\ar[dl]^{\gamma} E'\ar[dr]^{\alpha'} & & F'\ar[dl]_-{\beta} \\
& \wt F &                 & \wt F'
}\]
and the natural transformation
$\wt\ph':\cy(\wt F')\la H$. The assertions (i), (iii) and
(iv) of the Lemma have already been proved, as well as half of (ii). It remains
to show that $\wt F'\in\ct^-_c$ also belongs to $\ct^{\leq0}$, and that
the morphisms $\ch^\ell(\alpha')$ are isomorphisms for $\ell\geq-m+2$.

To prove this we recall that the commutative square
\[\xymatrix@C+20pt{
  \ch^\ell(E'_m) \ar[r]^-{\ch^\ell(g)}\ar[d]_-{\ch^\ell(a)} & \ch^\ell(E')\ar[d]^{\ch^\ell(\alpha')} \\
\ch^\ell(F')\ar[r]^-{\ch^\ell(\beta)} & \ch^\ell(\wt F')
}\] 
comes from applying $\ch^\ell$ to the weak triangle 
\[\xymatrix@C+20pt{
  E'_m \ar[r]^-{\s} & E'\oplus F' \ar[r]^-{(\alpha',\beta)} &\wt F'
  \ar[r]^-{\tau} &\T E'_m
}\]
The morphism
$\ch^\ell(\tau):\ch^\ell(\wt F')\la \ch^{\ell+1}(E'_m)$
fits in a long exact sequence. In (v) we proved that
$\ch^\ell(a):\ch^\ell(E'_m)\la\ch^\ell(F')$
is an isomorphism for $\ell\geq-m+2$, which makes the map
$\ch^\ell(\s):\ch^\ell(E'_m) \la \ch^\ell(E'\oplus F')$ a split
monomorphism for all $\ell\geq-m+2$. The exactness
tells us first that $\ch^\ell(\tau)=0$ for all
$\ell\geq-m+1$, and then that the
square
\[\xymatrix@C+20pt{
  \ch^\ell(E'_m) \ar[r]^-{\ch^\ell(g)}\ar[d]_-{\ch^\ell(a)} & \ch^\ell(E')\ar[d]^{\ch^\ell(\alpha')} \\
\ch^\ell(F')\ar[r]^-{\ch^\ell(\beta)} & \ch^\ell(\wt F')
}\] 
is bicartesian
for all $\ell\geq-m+2$. As long as $\ell\geq-m+2$,
the fact that $\ch^\ell(a)$
is an isomorphism forces $\ch^\ell(\alpha')$ to also be.

In particular: since $m\geq1$ we have that $-m+2\leq1$, and deduce that
$\ch^\ell(\wt F')\cong\ch^\ell(E')=0$ for all $\ell\geq1$.
Lemma~\ref{L-1.3} now gives that $\wt F'$ belongs to $\ct^{\leq0}$.
\eprf

\lem{L37.1046}
Let the conventions be as in Notation~\ref{N37.1024}, 
and assume
$H$ is a locally finite $\ct^c$--cohomological functor.
There exists an object $F\in\ct^-_c$ and an
epimorphism $\Hom(-,F)|_{\ct^c}\la H(-)$.
\elem

\prf
By Lemma~\ref{L37.1025} we may assume given:
an integer $\wt A>0$, objects $F_n\in\ct^-_c\cap\ct^{\leq\wt A}$ and
natural transformations $\ph_n:\cy(F_n)\la H$
which restrict to epimorphisms on $\gen Gn$. Replacing the functor
$H$ by $H(\T^{\wt A}-)$ we may assume $\wt A=0$. Let $B>0$ be
the integer whose existence is given by
Lemma~\ref{L17.25}.
Next we need to make our construction.

We will proceed inductively, using Lemma~\ref{L37.1036}, to construct
in $\ct^-_c\cap\ct^{\leq0}$ a sequence
\[\xymatrix@C+15pt{
E_1\ar[rr]^{\e_2}\ar[dr]_{\alpha_1}&  &\ar[dl]^{\gamma_2} E_2\ar[dr]_{\alpha_2}\ar[rr]^{\e_3} & & E_3\ar[dl]^-{\gamma_3}\ar[dr]_{\alpha_3} & \cdots \\
& \wt F_1 &                 & \wt F_2 & & \wt F_3
}\]
as well as natural transformations $\wt\ph_i:\cy(\wt F_i)\la H$,
satisfying the following
\be
\item
  The restriction of $\wt\ph_i$ to the subcategory $\gen G{(i+3)B}$
  is surjective.
\item
  For each $i>0$ the square
\[\xymatrix@C+40pt{
\cy(E_{i+1}) \ar[r]^-{\cy(\alpha_{i+1})}\ar[d]_-{\cy(\gamma_{i+1})} & \cy(\wt F_{i+1}) \ar[d]^{\wt\ph_{i+1}}\\
\cy(\wt F_i) \ar[r]^{\wt\ph_i} & H
}\]
commutes.
\item
  The morphisms $\ch^\ell(\gamma_i)$ and $\ch^\ell(\alpha_{i})$ are isomorphisms
  whenever $\ell\geq-i$.
\setcounter{enumiv}{\value{enumi}}
\ee
To start the induction we
declare $\wt F_1=F_{4B}$ and $\wt\ph_1=\ph_{4B}$. Choose a triangle
$E_1\stackrel{\alpha_1}\la\wt F_1\la D_1$
with $D_1\in\ct^{\leq-3}$; we immediately have that $\ch^\ell(\alpha_1)$
is an isomorphism for $\ell\geq-1$. In particular for $\ell\geq1$ we have
$\ch^\ell(E_1)\cong\ch^\ell(F_{4B})=0$; hence $E_1\in\ct^c\cap\ct^{\leq0}$.

Suppose our induction has proceeded as far as $n$.
In particular: we have produced in $\ct^-_c\cap\ct^{\leq0}$ a
morphism $\alpha_n:E_n\la \wt F_n$, with $E_n\in\ct^c$, as well
as a natural transformation $\wt\ph_n:\cy(\wt F_n)\la H$
which is surjective when restricted to $\gen G{(n+3)B}$.
And we have done it in such a way that $\ch^\ell(\alpha_n)$ is
an isomorphism for $\ell\geq-n$. We wish to go on to $n+1$.

Now the first paragraph of the proof gives us an object
$F_{(n+4)B}\in\ct^-_c\cap\ct^{\leq0}$, as well as a natural transformation
$\ph_{(n+4)B}^{}:\cy(F_{(n+4)B})\la H$ whose
restriction to $\gen G{(n+4)B}$ is surjective. Lemma~\ref{L37.1036},
with $m=n+3$,
allows us to construct in $\ct^-_c\cap\ct^{\leq0}$ the diagram
\[\xymatrix{
E_n\ar[rr]^{\e_{n+1}^{}}\ar[dr]_{\alpha_{n}}&  &\ar[dl]^{\gamma_{n+1}} E_{n+1}\ar[dr]^{\alpha_{n+1}} & & F_{(n+4)B}\ar[dl]_-{\beta} \\
& \wt F_n &                 & \wt F_{n+1}
}\]
as well as the natural transformation
$\wt\ph_{n+1}:\cy(\wt F_{n+1})\la H$
satisfying lots of properties:
Lemma~\ref{L37.1036}(i) tells us that $E_{n+1}$ belongs to $\ct^c$.
Lemma~\ref{L37.1036}(ii) gives that $\ch^\ell(\alpha_{n+1})$
and $\ch^\ell(\gamma_{n+1})$ are isomorphisms if $\ell\geq-n-3+2=-n-1$.
Lemma~\ref{L37.1036}(iii) says that $\wt\ph_{n+1}$ is an
epimorphism when restricted to $\gen G{(n+4)B}$,
and Lemma~\ref{L37.1036}(iv) says that the square in (ii) above
commutes. This finishes the induction.

It remains to see how to deduce the Lemma. We have produced a sequence
$E_1\stackrel{\e_2}\la E_3\stackrel{\e_3}\la E_3\stackrel{\e_4}\la \cdots$
and, for each $i>0$
\be
\setcounter{enumi}{\value{enumiv}}
\item
We define $\psi_i:\cy(E_i)\la H$ to be the composite
\[\xymatrix@C+40pt{
\cy(E_{i}) \ar[r]^-{\cy(\alpha_{i})} &
\cy(\wt F_i) \ar[r]^{\wt\ph_i} & H
}\]
\setcounter{enumiv}{\value{enumi}}
\ee
With these definitions we will prove
\be
\setcounter{enumi}{\value{enumiv}}
\item
The following triangles commute
\[\xymatrix@C+50pt@R-20pt{
    & \cy(E_{i+1}) \ar[dd]^{\psi_{i+1}}\\
  \cy(E_i) \ar[rd]_-{\psi_i}
  \ar[ur]^-{\cy(\e_{i+1})} & \\
  & H
}\]
\setcounter{enumiv}{\value{enumi}}
\ee
This permits us to make the next definitions
\be
\setcounter{enumi}{\value{enumiv}}
\item
We put $F=\hoco E_i$, and let $\psi:\cy(F)\la H$ be
\[\xymatrix@C+20pt{
  \cy(F) \ar@{=}[r] &\colim\,\cy(E_i)\ar[r] & H
}\]
that is the colimit of the maps $\psi_i$.
\setcounter{enumiv}{\value{enumi}}
\ee
And with the definitions made, we will prove
\be
\setcounter{enumi}{\value{enumiv}}
\item
The object $F$ belongs to $\ct^-_c$.
\item
  The map $\psi:\cy(F)|_{\ct^c}\la H$ is an epimorphism.
\setcounter{enumiv}{\value{enumi}}
\ee
Together, (vii) and (viii) contain the assertion of the Lemma. All that remains
is to prove (v), (vii) and (viii).

To prove (v) the reader should consider the diagram
\[\xymatrix@C+40pt@R-15pt{
  &\cy(E_{i+1}) \ar[r]^-{\cy(\alpha_{i+1})}\ar[dd]^-{\cy(\gamma_{i+1})} & \cy(\wt F_{i+1}) \ar[dd]^{\wt\ph_{i+1}}\\
  \cy(E_{i}) \ar[ur]^-{\cy(\e_{i+1})}
  \ar[dr]_-{\cy(\alpha_{i})} & & \\
&\cy(\wt F_i) \ar[r]^{\wt\ph_i} & H
}\]
We wish to prove the commutativity
of the perimeter.
The square commutes by (ii), and the triangle by applying the
functor $\cy$ to the commutative triangle
\[\xymatrix@C+15pt{
E_i\ar[rr]^{\e_{i+1}}\ar[dr]_{\alpha_i}&  &\ar[dl]^{\gamma_{i+1}} E_{i+1}\\
& \wt F_i &                
}\]

Also: we may apply $\ch^\ell$ to the commutative triangle.
 From (iii) we know that $\ch^\ell(\alpha_i)$
is an isomorphism for $\ell\geq-i$ and also that 
$\ch^\ell(\gamma_{i+1})$ is an isomorphism for $\ell\geq-i-1$. Since
$\ch^\ell(\alpha_i)=\ch^\ell(\gamma_{i+1})\ch^\ell(\e_{i+1})$
we learn that $\ch^\ell(\e_{i+1})$ is an isomorphism when $\ell\geq-i$.
Therefore $E_*$ is a strong $\ct^c$--approximating system as
in Definition~\ref{D1.29}, and Lemma~\ref{L1.1010101} informs
us that $F=\hoco E_i$ belongs to $\ct^-_c$---that is we have proved (vii).
It remains only to prove (viii).

Suppose therefore that $C$ is an object $\ct^c$. There exists
an $n>0$ with $C\in\gen G{(n+3)B}$, and Lemma~\ref{L17.17} says that we may
also choose $n$ so that $\Hom(C,\ct^{-n-1})=0$. Because
$\wt\ph_n:\cy(\wt F_n)\la H$ is surjective on
$\gen G{(n+3)B}$ we have that the map
$\wt\ph_n:\Hom(C,\wt F_n)\la H(C)$ is surjective.
In the triangle $E_n\stackrel{\alpha_n}\la \wt F_n\la \wt D_n$ we have
that $\ch^\ell(\alpha_n)$ is an isomorphism for $\ell\geq-n$, therefore
$\ch^\ell(\wt D_n)=0$ for $\ell\geq-n$, therefore $\wt D_n\in\ct^{\leq-n-1}$. 
In the exact sequence
\[\xymatrix@C+15pt{
\Hom(C,E_n)\ar[rr]^-{\Hom(C,\alpha_n)}& &\Hom(C,\wt F_n)\ar[r] &               
\Hom(C,\wt D_n)
}\]
we have $\Hom(C,\wt D_n)=0$ since
$\wt D_n\in\ct^{\leq-n-1}$. The
map $\psi_n:\Hom(C,E_n)\la H(C)$ of (iv) is the composite of the
two epimorphisms
\[\xymatrix@C+15pt{
  \Hom(C,E_n)\ar[rr]^-{\Hom(C,\alpha_n)}&&
  \Hom(C,\wt F_n)\ar[r]^-{\wt\ph_n} &          
H(C)
}\]
But it factors through $\psi:\Hom(C,F)\la H(C)$, which must
therefore be epi.
\eprf

\rmd{R37.1050}
Let $\ct$ be a triangulated category with coproducts.
A morphism $f:D\la E$ is called
\emph{phantom} if, for every compact object $C\in\ct$, the induced map
$\Hom(C,f):\Hom(C,D)\la\Hom(C,E)$ vanishes. The phantom maps
form an ideal: if $f,f':D\la E$ are phantom then so is $f+f'$, and if
$D'\stackrel e\la D\stackrel f\la E\stackrel g\la E'$ are composable
morphisms with $f$ phantom, then $gfe:D'\la E'$ is also
phantom.
\ermd

\cor{C37.1056}
Let the conventions be as in Notation~\ref{N37.1024}.
Let 
$F'\in\ct$ be an object
such that the functor $H=\cy(F')$
is a locally finite $\ct^c$--cohomological functor.
There exists an object $F\in\ct^-_c$ and a triangle
$F\stackrel f\la F'\stackrel g\la D$ with
$g$ phantom.
\ecor

\prf
Lemma~\ref{L37.1046} gives us an object $F\in\ct^-_c$ and an
epimorphism $\ph:\cy(F)\la H=\cy(F')$.
Since $F$ belongs to $\ct^-_c$ Lemma~\ref{L1.1010101}(ii) produces
for $F$ a (strong) $\ct^c$--approximating system. Lemma~\ref{L27.976}
allows us to realize the natural transformation
$\ph$ as $\cy(f):\cy(F)\la \cy(F')$
for some (non-unique) $f:F\la F'$. Complete $f$ to a triangle
$F\stackrel f\la F'\stackrel g\la D$.
For every object
$C\in\ct$ we have an exact sequence
\[\xymatrix@C+15pt{
  \Hom(C,F)\ar[rr]^-{\Hom(C,f)}&&
  \Hom(C,F')\ar[rr]^-{\Hom(C,g)} &&          
\Hom(C,D)
}\]
and if $C$ is compact the map $\Hom(C,f)=\cy(f)(C)$ is surjective. It follows
that $\Hom(C,g)$ is the zero map.
\eprf

Now that it's time to state the main theorem we include all the hypotheses
explicitly.

\thm{T37.1067}
Let $R$ be a noetherian, commutative ring.
Let $\ct$ be an $R$--linear triangulated category with coproducts,
and suppose it has
a compact generator $G$ 
such that $\Hom(-,G)$ is a $G$--locally finite cohomological functor.
Assume further that $\ct$ is approximable.

Let $\ct^-_c\subset\ct$ be the category of Definition~\ref{D0.13}, where
the \tstr\ with respect to which we define it is in the preferred
equivalence class. Then the functor 
$\cy:\ct^-_c\la \Hom\big((\ct^c)\op,\Mod R\big)$, taking
$F\in\ct^-_c$ to $\cy(F)=\Hom(-,F)|_{\ct^c}$, satisfies
\be
\item
  The objects in the essential image of $\cy$
  are the locally finite $\ct^c$--cohomological functors.
\item
  The functor $\cy$ is full.
\setcounter{enumiv}{\value{enumi}}
\ee
\ethm

\prf
The fact that, for any object $F\in\ct^-_c$, the functor
$\cy(F)$
is a locally finite $\ct^c$--cohomological functor was proved
in Lemma~\ref{L1.1.3}. In Lemma~\ref{L1.1010101}(ii) we saw that any
$F\in\ct^-_c$ admits a $\ct^c$--approximating system, and Lemma~\ref{L27.976}
guarantees that any natural transformation
$\ph:\cy(F)\la \cy(F')$ can be expressed as
$\ph=\cy(f)$ for some $f:F\la F'$; that is the functor is
full. It remains to show that any
locally finite $\ct^c$--cohomological functor $H:(\ct^c)\op\la\Mod R$
is in the essential image; we must show it isomorphic to
$\cy(F)$ for some $F\in\ct^-_c$.

Proposition~\ref{P37.1005} produced a candidate $F$; we have
an $F\in\ogenun\cg4$ and an isomorphism $H\cong \cy(F)$.
We wish to show that $F$ belongs to $\ct^-_c$. We proceed by induction
to prove
\be
\setcounter{enumi}{\value{enumiv}}
\item
  Let $\ci$ be the ideal of phantom maps. For
  each integer $n>0$ there exists a triangle
  $F_n\la F\stackrel{\beta_n}\la D_n$ with $F_n\in\ct^-_c$ and
  $\beta_n\in\ci^n$.
\setcounter{enumiv}{\value{enumi}}
\ee
We prove (iii) by induction on $n$. The case $n=1$ is given by
Corollary~\ref{C37.1056}. Now for the inductive step:
assume that, for some $n\geq1$, we are given
a triangle $F_n\la F\stackrel{\beta_n}\la D_n$ with $F_n\in\ct^-_c$ and
$\beta_n\in\ci^n$.
We know that both $\cy(F_n)$
and $\cy(F)$ are locally finite $\ct^c$--cohomological
functors, and the exact sequence
\[\xymatrix{
  0\ar[r] & \cy(\Tm D_n) \ar[r] & \cy(F_n)\ar[r] &
  \cy(F)\ar[r] & 0
}\]
says that so is $\cy(\Tm D_n)$. Corollary~\ref{C37.1056}
permits us to construct a triangle $F'\la D_n\stackrel\gamma\la D_{n+1}$
with $F'\in\ct^-_c$ and $\gamma\in\ci$. Let
$\beta_{n+1}:F\la D_{n+1}$ be the
composite $F\stackrel{\beta_n}\la D_{n} \stackrel\gamma\la D_{n+1}$.
Since $\beta_n\in\ci^n$ and $\gamma\in\ci$ we deduce
that $\beta_{n+1}\in\ci^{\n+1}$.
If we complete $\beta_{n+1}$ to a
triangle $F_{n+1}\la F\stackrel{\beta_{n+1}}\la D_{n+1}$,
the octahedral axiom allows us to find a triangle $F_n\la F_{n+1}\la F'$.
Since $F_n$ and $F'$ both lie in $\ct^-_c$ so does $F_{n+1}$. This
completes the proof of (iii).

Now consider the triangle $F_4\la F\la D_4$.
The morphism $F\la D_4$ is in $\ci^4$, but $F$ belongs to $\ogenun\cg4$.
One easily shows that
$\big(\ogenun\cg1,\ci\big)$ is a projective
class as in Christensen~\cite[Definition~2.2]{Christensen96},
and \cite[Theorem~1.1]{Christensen96} tells us that
so is $\big(\ogenun\cg4,\ci^4\big)$.
The map $F\la D_4$ is a morphism in $\ci^4$ out
of an object in $\ogenun\cg4$ and must vanish, making $F$ a direct summand
of $F_4\in\ct^-_c$. Proposition~\ref{P17.21} tells us that $\ct^-_c$ is
thick, and therefore $F\in\ct^-_c$.
\eprf

\lem{L37.2063}
Let the assumptions be as in
Theorem~\ref{T37.1067}. Suppose $f:F\la F'$ is a morphism
in $\ct^-_c$ and assume $F'$ belongs to $\ct^b_c$. Then
$\cy(f)=0$ implies $f=0$.
\elem

\prf
Because $F'$ belongs to $\ct^b_c$ there must be an integer $\ell$
with $F'\in\ct^{\geq\ell}$; without loss of generality we may assume
$\ell=0$. Now $F$ belongs to $\ct^-_c$, hence there must exist a triangle
$E\stackrel g\la F\stackrel h\la D$ with $E\in\ct^c$ and $D\in\ct^{\leq-1}$.
Choose
such a triangle.

The vanishing of $\cy(f)$ means that $\Hom(E,f):\Hom(E,F)\la\Hom(E,F')$
must take $g\in\Hom(E,F)$ to zero; that means $fg=0$. But
the triangle tells us that $f:F\la F'$ must factor as
$F\stackrel h\la D\la F'$. As $D\in\ct^{\leq-1}$ and $F'\in\ct^{\geq0}$ we
have $\Hom(D,F')=0$, hence $f=0$.
\eprf

\thm{T37.2068}
Let the assumptions be as in
Theorem~\ref{T37.1067}. The restriction of the functor $\cy$ to
the subcategory $\ct^b_c$ is fully faithful, and the essential
image is the class of finite $\ct^c$--cohomological functors.
\ethm

\prf
The functor is full on all of $\ct^-_c$, and Lemma~\ref{L37.2063}
guarantees that on the subcategory $\ct^b_c$ it is faithful. It
remains to identify the essential image.
Let $F$ be an object in $\ct^-_c$, we need to show that
$\cy(F)$ is finite if and only if $F\in\ct^b_c$.

Suppose
$F\in\ct^b_c\subset\ct^+$ and
$C\in\ct^c\subset\ct^-$. We can choose an integer $\ell>0$ so that
$\Hom(\T^iC,F)=0$ for all $i\geq\ell$, which implies
that $\Hom(-,F)|_{\ct^c}$ is $C$--finite. Since this is true for every
$C\in\ct^c$ we have that $\cy(F)$ is finite.

Conversely: suppose $\cy(F)$ is finite and choose a
compact generator $G$. Because $\cy(F)=\Hom(-,F)|_{\ct^c}$
is $G$--finite there is an integer $\ell$ so that $\Hom(\T^iG,F)=0$
for all $i\geq\ell$. But then $\Hom(T,F)=0$ for
all $T\in\ogenul G{}{-\ell}=\ct_G^{\leq-\ell}$ and $F$ must belong to
$\ct^-_c\cap\ct_G^{\geq-\ell+1}\subset\ct^b_c$.
\eprf

\section{Implications of the hypothesis $\ct=\ogenun{G'}N$, with $G'\in\ct^b_c\subset\ct^b$}
\label{SS1}

We have proved Theorem~\ref{T1.-1}(i), and the time
has now come to turn our attention to Theorem~\ref{T1.-1}(ii).
And, as the reader can easily check, in
Theorem~\ref{T1.-1}(ii) we added the hypothesis that
there exist an object $G'\in\ct^b_c$ with
$\ct=\ogenun{G'}N$. This section
will focus on finding out
what this extra hypothesis buys us.

But before all else we prove 
the easy part of Theorem~\ref{T1.-1}(ii),
showing that the
images of the functors $\wt\cy$ and
$\wt\cy\circ\wi$ are contained
where the theorem asserts they should be. For
this the extra hypothesis, regarding the existence
of an object $G'\in\ct^b_c$ with
$\ct=\ogenun{G'}N$, is unnecessary.

\lem{LL3.-1}
Let $R$ be a commutative, noetherian ring. Suppose $\ct$ is
an $R$--linear triangulated
category.  Suppose $\ct$ has a single  compact generator
$G$, such that $\Hom(G,\T^iG)$ is a finite $R$--module for every $i\in\zz$,
and vanishes when $i\gg0$.

Then $\Hom(A,\T^iB)$ is a finite $R$--module whenever $A\in\ct^-_c$ and
$B\in\ct^b_c$. For fixed $A$ and $B$ it vanishes when $i\ll0$.

If $A$ belongs to $\ct^c\subset\ct^-_c$, then $\Hom(A,\T^iB)$ also vanishes
for $i\gg0$.
\elem

\prf
Fix $A$ and $B$. Since $A\in\ct^-$ and $B\in\ct^+$ there will be some
integer $m>0$ such that $\Hom(A,\T^iB)=0$ for $i<-m$.

We need to prove the finiteness for fixed $i$, and without loss we may assume
$i=0$. Shifting if necessary we may assume $B\in\ct^{\geq0}$. But $A\in\ct^-_c$
means that there must exist a triangle $E\la A\la D$ with $E\in\ct^c$ and
$D\in\ct^{\leq-2}$. In the exact sequence
\[\xymatrix@C+20pt{
\Hom(D,B)\ar[r] &\Hom(A,B)\ar[r] &\Hom(E,B)\ar[r] &\Hom(\Tm D,B)
}\]
we have $\Hom(D,B)=0=\Hom(\Tm D,B)$, hence $\Hom(A,B)\cong\Hom(E,B)$.
The finiteness of 
$\Hom(E,B)$ is contained in Lemma~\ref{L1.1.3}.

Now assume $A\in\ct^c$. The vanishing of $\Hom(A,\T^iB)$ for $i\gg0$ follows
from Lemma~\ref{L17.17}.
\eprf

And now we turn our attention to understanding categories
of the form $\ogenun{G'}N$, with $G'$ first assumed
to belong to $\ct^b$, and as we progress we will make
the more restrictive hypothesis that $G'$ belongs to
$\ct^b_c\subset\ct^b$.

On a matter of notation: we will focus mostly
on the object $G'\in\ct^b_c\subset\ct$, but on occasion
we will recall that the entire theory assumes the
existence of a compact generator $G\in\ct$. Because
the symbols $G$ and $G'$ are easy to confuse, in the
remainder of this section the compact generator
will be $H\in\ct^c$ and the object in $\ct^b_c\subset\ct$
will be called $G$.

\lem{LL1.-5}
Let $\ct$ be a triangulated category with coproducts, and let
$\big(\ct^{\leq0},\ct^{\geq0}\big)$ be a \tstr\ on $\ct$.
Let $G\in\ct^b$ be an object, and let $N>0$ be an integer.

Then there exists an integer $M>0$ so that
any map $A\la C$, with $A\in\ct^{\leq0}$ and
with $C\in \genuf GN{-m}$,  must factor
through some object $B\in\genu GN{-m,M}$.

If we further assume
that $\ct^{\geq0}$ is closed under coproducts,
then the integer $M>0$ may be chosen so that
any map $A\la C$, with $A\in\ct^{\leq0}$ and $C\in \ogenuf GN{-m}$,  will factor
through some object $B\in\ogenu GN{-m,M}$. 

In both statements we allow 
$m>0$ to be (possibly) infinite.
\elem

\prf
Because $G\in\ct^b$ we may choose an integer $K>0$ so that
$G\in\ct^{\geq-K}\cap\ct^{\leq K}$. Fix such a $K$; for any integer $n>0$
we clearly have
$\genuf GnK\subset\ct^{\geq0}$, and if
 $\ct^{\geq0}$ is
closed under coproducts then we also have
$\ogenuf GnK\subset\ct^{\geq0}$.

Next we proceed to prove the assertion of the lemma beginning with
``if we further assume''. We leave to
the reader the case with $C\in \genuf GN{-m}$.
We are given that $C\in\ogenuf GN{-m}$, and by
\cite[Corollary~1.12]{Neeman17} we have
\[\ogenuf GN{-m}=\Smr\Big[\Coprod_N^{}\big(G[-m,\infty)\big)\Big].\]
Therefore we may choose an object $C'$ with
$C\oplus C'\in\Coprod_N^{}\big(G[-m,\infty)\big)$, and any map
$f:A\la C$ obviously factors as the composite
$A\stackrel f\la C\stackrel i\la C\oplus C'\stackrel\pi\la C$, meaning $f$
factors 
through $C\oplus C'$. It therefore
suffices to show that, with $K>0$ as in the paragraph above, we have
\be  
\item
Any map $f:A\la C$, with $A\in\ct^{\leq0}$ and
$C\in\Coprod_N^{}\big(G[-m,\infty)\big)$, must factor through
some object $B\in\Coprod_N^{}\big(G[-m,(2K+1)N]\big)$.
\ee
Now we proceed by induction on $N$. If $N=1$ we are
given a map $A\la C$, with $A\in\ct^{\leq0}$ and 
$C\in \Coprod_1^{}\big(G[-m,\infty)\big)$. But
\[  
\Coprod_1^{}\big(G[-m,\infty)\big)\eq
\Coprod_1^{}\big(G[-m,K]\big)\bigoplus\Coprod_1^{}\big(G[K+1,\infty)\big)
\]
As $\Coprod_1^{}\big(G[K+1,\infty)\big)$ is contained in $\ct^{\geq1}$ and
$A\in\ct^{\leq0}$,
the map $A\la C$ must factor through
$B\in\Coprod_1^{}\big(G[-m,K]\big)$. We have proved an
improvement on (i) in the case $N=1$.

Next assume we know (i) for all integers $\leq N$, and keep
in mind that, for $N=1$,
we proved an improvement on (i) in the last paragraph. Now
let $\cs=\ct^{\leq0}$ and put
\[\begin{array}{ccccccc}
\cx&=&\Coprod_N^{}\big(G[-m,\infty)\big)&\qquad&
\ca&=&\Coprod_N^{}\big(G[-m,(2K+1)(N+1)]\big)\\
\cz&=&\Coprod_1^{}\big(G[-m,\infty)\big)&\qquad&
\cc&=&\Coprod_1^{}\big(G[-m,K]\big)    
\end{array}\]
By the induction any pair of maps $s\la x$ and  $s\la z$,
with $s\in\cs$, $x\in\cx$ and
$z\in\cz$, factor (respectively)
as $s\la a\la x$ and $s\la c\la z$, with $a\in\ca$ and
$c\in\cc$. 

Now let $d$ be an object of $\big(\Tm\cc\big)*\cs$. As
$\Tm\cc\subset\ct^{\leq 2K+1}$ and $\cs=\ct^{\leq0}$ we deduce
that $d\in\ct^{\leq 2K+1}=\T^{-2K-1}\cs$, and induction tells us
that any map $d\la x$, with $x\in\cx$, must
factor as $d\la a\la x$ with $a\in\ca$. The hypotheses
of \cite[Lemma~1.6]{Neeman17} are satisfied, hence any morphism
$s\la\cx*\cz=\Coprod_{N+1}^{}\big(G[-m,\infty)\big)$ factors
through $\ca*\cc\subset\Coprod_{N+1}^{}\big(G[-m,(2K+1)(N+1)]\big)$.
\eprf

\lem{LL1.-23}
Let $\ct$ be a triangulated category with coproducts, as well as a \tstr\
$\big(\ct^{\leq0},\ct^{\geq0}\big)$.
Let $G\in\ct^b$ be an object, and let $N>0$ be an integer.
Then there exists an integer $M>0$ so that
\be
\item
Any map $C\la A$, with $A\in\ct^{\geq0}$ and $C\in \genul GN{m}$,  must factor
through some object $B\in\genu GN{-M,m}$.
\item
Suppose $\ct^{\geq0}$ is closed under coproducts.
Then
any map $C\la A$, with $A\in\ct^{\geq0}$ and $C\in \ogenul GN{m}$,  must factor
through some object $B\in\ogenu GN{-M,m}$.
\setcounter{enumiv}{\value{enumi}}
\ee
In both statements the integer
$m>0$ is possibly infinite.
\elem

\prf
Because $G\in\ct^b$ we may choose an integer $K>0$ so that
$G\in\ct^{\geq-K}\cap\ct^{\leq K}$. Fix such a $K$. The fact that
$\ct^{\leq0}={^\perp\ct^{\geq1}}$ tells us that $\ct^{\leq0}$ is (automatically)
closed under coproducts, and hence for any integer $n>0$ we have
$\genul Gn{-K}\subset\ogenul Gn{-K}\subset\ct^{\leq0}$.

The rest of the proof is just the dual of the proof of Lemma~\ref{LL1.-5};
we sketch the proof of assertion (ii) and leave (i) to the reader. 
Any map $C\la A$, with $A\in\ct^{\geq0}$ and $C$ in
\[
\Coprod_1^{}\big(G(-\infty,m]\big)\eq
\Coprod_1^{}\big(G(-\infty,-K-1]\big)\bigoplus
\Coprod_1^{}\big(G[-K,m]\big)
\]
must factor through some $B\in\Coprod_1^{}\big(G[-K,m]\big)$. Now
proceeding dually to the proof of Lemma~\ref{LL1.-5}
one shows, inductively on $N$, that
\be
\setcounter{enumi}{\value{enumiv}}
\item
Any map $f:C\la A$, with $A\in\ct^{\geq0}$ and
$C\in\Coprod_N^{}\big(G(-\infty,m]\big)$, must factor through
some object $B\in\Coprod_N^{}\big(G[-(2K+1)N,m]\big)$.
\ee
More precisely: let $\cs=\ct^{\geq0}$ and
\[\begin{array}{ccccccc}
\cx&=&\Coprod_N^{}\big(G(-\infty,m]\big)&\qquad&
\ca&=&\Coprod_N^{}\big(G[-(2K+1)(N+1),m]\big)\\
\cz&=&\Coprod_1^{}\big(G[-\infty,m]\big)&\qquad&
\cc&=&\Coprod_1^{}\big(G[-K,m]\big)    
\end{array}\]
The hypothesis that $\ct^{\geq0}$ is closed under coproducts
gives that $\T\cc\subset\ct^{\geq-2K-1}=\T^{2K+1}\cs$, and
hence $\cs*\T\cc\subset\cs*\T^{2K+1}\cs\subset\T^{2K+1}\cs$.
Now induction coupled with \cite[Lemma~1.6]{Neeman17} guarantee that 
any map $y\la s$, with $y\in\cz*\cx$ and $s\in\cs=\ct^{\geq0}$,
must factor through
$B\in\cc*\ca\subset\Coprod_{N+1}^{}\big(G[-(2K+1)(N+1),m]\big)$.
\eprf

\lem{LL1.-27}
Let $\ct$ be a triangulated category with coproducts, as well as a \tstr\
$\big(\ct^{\leq0},\ct^{\geq0}\big)$ with
$\ct^{\geq0}$ is closed under coproducts.
Let $G\in\ct^b$ be an object, let $N>0$ be an integer,
and assume $\ct=\ogenun GN$.

Then there exists an integer $M>0$ such that, for any integers $a\leq b$,
\[
\ct^{\geq0}\subset\ogenuf GN{-M},\qquad
\ct^{\leq0}\subset\ogenul GN{M},\qquad
\ct^{\geq a}\cap\ct^{\leq b}\subset \ogenu GN{a-M,b+M}.
\]
\elem

\prf
For the given integer $N$ and object $G$, pick an integer $M$ so that
Lemmas~\ref{LL1.-5} and \ref{LL1.-23} both hold.

Let $a$ be an integer, and suppose $A$ is an object in $\ct^{\geq a}$.
The identity map $A\la A$ is a morphism from $A\in\ogenun GN$
to $A\in\ct^{\geq a}$, and Lemma~\ref{LL1.-23} (with $m=\infty$)
guarantees that the map factors
through an object $B\in\ogenuf GN{a-M}$. As $A$ is a direct summand of
$B$ it must also lie in $\ogenuf GN{a-M}$. The case $a=0$ gives
the first assertion
of the Lemma.

Now assume $A\in\ct^{\geq a}\cap\ct^{\leq b}$
with $a$ possibly equal
to $-\infty$.
By the first half of the Lemma, already proved, $A$ must lie
in $\ogenuf GN{a-M}$. The identity map $A\la A$ is a morphism from
$A\in\ct^{\leq b}$ to $A\in\ogenuf GN{a-M}$, and Lemma~\ref{LL1.-5}
guarantees that it factors through some $B\in\ogenu GN{a-M,b+M}$. Therefore
$A$, being a direct summand of $B$, must also belong to $\ogenu GN{a-M,b+M}$.
\eprf

\lem{LL1.3}
Let $\ct$ be a triangulated category with coproducts and a single compact
generator $H$. Assume $\Hom(H,\T^nH)=0$ for $n\gg0$.

Let
$\big(\ct^{\leq0},\ct^{\geq0}\big)$ be a \tstr\ on $\ct$,
in the preferred equivalence class. Given an integer $K>0$
and a collection of objects
$\{X_i\in\ct^{\leq K}\cap\ct^{\geq-K}\mid i\in\zz\}$, then the map
\[\xymatrix@C+40pt{
\ds\coprod_{i=-\infty}^\infty \T^iX_i \ar[r]^-\ph & \ds\prod_{i=-\infty}^\infty \T^iX_i
}\]
is an isomorphism.
\elem

\prf
Because the \tstr\ is in the preferred equivalence class and $H$ is compact,
Observation~\ref{O0.2.5} tells us that there is some
integer $A>0$ with $H\in\ct^{\leq A-1}$, and therefore
$\Hom(H,-)$ vanishes on $\ct^{\geq A}=\big[\ct^{\leq A-1}\big]^\perp$.
And Lemma~\ref{L17.17} allows us to assume, possibly after increasing
$A$, that $\Hom\big(H,\ct^{\leq-A}\big)=0$. 
Therefore
the functor $\Hom(H,-)$ vanishes on the union $\ct^{\leq-A}\cup\ct^{\geq A}$.
As $X_i$ is assumed to belong to $\ct^{\leq K}\cap\ct^{\geq-K}$, we have that
$\Hom\big(H,\T^nX_i\big)=0$ whenever $|n|>A+K$.

Any morphism $\T^nH\la \coprod_{i=-\infty}^\infty \T^iX_i $
will factor through a finite
subcoproduct by the compactness of $H$, and any map
$\T^nH\la \prod_{i=-\infty}^\infty \T^iX_i $ will factor through a finite subproduct
by the vanishing of $\Hom\big(H,\T^{i-n}X_i\big)$ for all but finitely many $i$.
Therefore the functor $\Hom(\T^nH,-)$ takes the map $\ph$ to an isomorphism,
for every $n\in\zz$, and as $H$ is a generator the map $\ph$
must be an isomorphism.
\eprf

\pro{PP1.2000}
Let $\ct$ be a triangulated category with coproducts, and
assume it has a compact generator $H$ with $\Hom(H,\T^iH)=0$
for $i\gg0$. Suppose $G$ is an object
in $\ct^b$, where we mean the $\ct^b$ that comes from the preferred
equivalence class of {\it t}--structures. Assume that
there exists an integer $N>0$ with $\ogenun GN=\ct$.

Suppose
$F\in\ct^-_c$ is an object such that $\Hom(F,\T^iG)=0$ for $i\gg0$.
Then $F$ is compact.
\epro

\prf
In the preferred equivalence class choose a \tstr\ 
$\big(\ct^{\leq0},\ct^{\geq0}\big)$ so that $\ct^{\geq0}$ is
closed under coproducts.
We are given that $G\in\ct^b$, hence we may choose an integer $K>0$ with
$G\in\ct^{\geq-K}\cap\ct^{\leq K}$. Because $F\in\ct^-_c$, for any integer $i>0$ we
may form a triangle $E_i\stackrel{\alpha_i}\la F\la D_i\la$ with $E\in\ct^c$ and
$D_i\in\ct^{\leq-i-K-2}$. For any $X\in\ct$ we have the exact sequence
\[\xymatrix{
\Hom(D_i,X)\ar[r] & \Hom(F,X)\ar[rr]^-{\Hom(\alpha_i,X)}& & \Hom(E_i,X)\ar[r] &\Hom(\Tm D_i,X)
}\]
and if $X\in\ct^{\geq-i-K}$ then $\Hom(D_i,X)=0=\Hom(\Tm D_i,X)$.
Thus $\Hom(\alpha_i,X)$ is an isomorphism for all $X\in\ct^{\geq-i-K}$.

Now choose an integer $A>0$ so that $\Hom(F,\T^iG)=0$ for all $i\geq A$.
Pick an integer $i\geq A$; the vanishing of $\Hom(F,\T^iG)$
implies the 
vanishing of $\Hom(E_i,\T^iG)$, and the compactness of $E_i$
implies the vanishing of of $\Hom\big(E_i,\T^i\Coprod_1^{}(G)\big)$.
But $\ct^{\geq-i-K}$ is closed under coproducts and contains
$\T^iG$, hence
$\T^i\Coprod_1^{}(G)$ is contained in $\ct^{\geq-i-K}$, and
the vanishing of $\Hom\big(E_i,\T^i\Coprod_1^{}(G)\big)$ implies the vanishing of
$\Hom\big(F,\T^i\Coprod_1^{}(G)\big)$.

Thus $\Hom(F,-)$ vanishes on $\T^i\Coprod_1^{}(G)$ for any
$i\geq A$. Now
every object $X\in\Coprod_1^{}\big(G(-\infty,-A]\big)$ can be written
as $\coprod_{i=A}^{\infty}\T^iX_i$ with
$X_i\in\Coprod_1^{}(G)\subset\ct^{\geq-K}\cap\ct^{\leq K}$. By Lemma~\ref{LL1.3}
there is an isomorphism
$\coprod_{i=A}^{\infty}\T^iX_i\cong\prod_{i=A}^{\infty}\T^iX_i$,
and hence $\Hom(F,-)$ vanishes on the category
$\Coprod_1^{}\big(G(-\infty,-A]\big)$.

It immediately follows that $\Hom(F,-)$ also vanishes on the category
$\ogenul GN{-A}$. But now Lemma~\ref{LL1.-27} establishes the existence
of some integer $m>0$ with $\ct^{\leq-m}\subset\ogenul GN{-A}$, and hence
$\Hom(F,-)$ vanishes on $\ct^{\leq-m}$. On the other hand $F$ belongs to
$\ct^-_c$, hence we may choose a triangle $E\la F\la D\la$ with $E\in\ct^c$
and $D\in\ct^{\leq-m}$, and as the map $F\la D$ must vanish we have
that $F$ is a direct summand of $E\in\ct^c$.
\eprf

Until now, the object $G$ with $\ct=\ogenun GN$ was only assumed to lie in
$\ct^b$. From now on we strengthen the hypotheses on $G$,
and assume it belongs to $\ct^b_c\subset\ct^b$.

\lem{LL1.10909}
Let $\ct$ be a triangulated category with coproducts and a single compact
generator $H$. Assume $\Hom(H,\T^nH)=0$ for $n\gg0$. Suppose $G$ is an object
in $\ct^b_c$, where we mean the $\ct^b_c$ that comes from the preferred
equivalence class of {\it t}--structures.

For any integer $N>0$, any map
$F\la Y$, with $F\in\ct^-_c$ and $Y\in\ogenu GN{a,b}$, factors through
an object in $\genu GN{a,b}$. We allow either or both of $b,N$ to be infinite.

For the sake of clarity: allowing $N$ to be infinite means that any map
$F\la Y$, with $F\in\ct^-_c$ and $Y\in\ogenu G{}{a,b}$, factors through
an object in $\genu G{}{a,b}$.
\elem

\prf
Choose a \tstr\ in the preferred equivalence class, and pick it so that
$\ct^{\geq0}$ is closed under coproducts.
By hypothesis $G$ is contained in $\ct^b$, hence we may choose and fix an
integer $K>0$ with $G\in\ct^{\geq-K}$. Then, for
any (finite or infinite) integer $N>0$, we have
\[
\genu GN{a,b}\sub \ogenu GN{a,b}\sub \ct^{\geq a-K}\ .
\]
Next we observe
\be
\item
For any object 
$F\in\ct^-_c$ we may choose triangle $\Tm D\la E\la F\la D$ with
$E\in\ct^c$ and $D\in\ct^{\leq a-K-2}$, and hence
\[
\Hom\big(\Tm D,\ct^{\geq a-K}\big)\eq0\eq
\Hom\big(D,\ct^{\geq a-K}\big)\ .
\]
Thus, for any object $Y\in\ogenu GN{a,b}$, the natural map
$\Hom(F,Y)\la\Hom(E,Y)$ is an isomorphism. 
\setcounter{enumiv}{\value{enumi}}
\ee

Now we prove the case $N<\infty$ of the lemma by induction on
$N$. Suppose first that $N=1$; then
any object $Y\in\ogenu G1{a,b}$ is a direct summand of
an object in $\Coprod_1^{}\big(G[a,b]\big)$, hence we may assume
$Y$ belongs to
$\Coprod_1^{}\big(G[a,b]\big)\subset\ogenu G1{a,b}$.
Because $F\in\ct^-_c$ we may choose a triangle $\Tm D\la E\la F\la D$ as
in (i). But $E$ is compact and $Y$ is a
coproduct, hence any map $E\la Y$ factors through a finite subcoproduct.
Therefore so does any map $F\la Y$, completing the proof in the case $N=1$.

Now assume we know the Lemma for all integers $i$ with $0<i\leq N$. Put 
\[
\cs=\ct^-_c\,,\quad\ca=\genu G1{a,b},\quad\cc=\genu GN{a,b},\quad\cx=\ogenu G1{a,b},\quad\cz=\ogenu GN{a,b}.
\]
By the hypotheses of the lemma $\cs=\ct^-_c$ is triagulated and contains
$\cc$. The inductive hypothesis,
coupled with \cite[Lemma~1.6 and Remark~1.7]{Neeman17}, give that
any map $s\la y$, with $s\in\cs$ and $y\in\cx*\cz$,
must factor through an object
$b\in\ca*\cc$. The case of the lemma where $N<\infty$ follows immediately.

To prove the statement with $N=\infty$, let $\car\subset\ogenu G{}{a,b}$ be
the full subcategory of all objects $Y\in\ogenu G{}{a,b}$ such that
any map $F\la Y$, with $F\in\ct^-_c$, factors through
an object in $\genu G{}{a,b}$. We need to prove that
$\car=\ogenu G{}{a,b}$.

By the case where $N<\infty$ we already know that $\ogenu G1{a,b}\subset\car$,
and it is obvious that any direct summand of an object in $\car$ belongs to
$\car$. It remains to prove that $\car$ is closed under coproducts
and extensions.

Therefore let $Y=\coprod_{\lambda\in\Lambda}^{}Y_\lambda$ be a coproduct
of objects in $\car$, and let $F\la Y$ be a morphism with $F\in\ct^-_c$.
By (i) we may choose a morphism $E\la F$, with $E\in\ct^c$, and such
that $\Hom(F,Z)\la\Hom(E,Z)$ is an isomorphism for any
$Z\in\ogenu G{}{a,b}$. In particular it is an isomorphism for
$Z=Y=\coprod_{\lambda\in\Lambda}^{}Y_\lambda$. But $E$ is compact and hence
the morphism $E\la\coprod_{\lambda\in\Lambda}^{}Y_\lambda$ factors through a
finite subcoproduct. And for each $\lambda$ in the resulting finite set, the
fact that $Y_\lambda\in\car$ allows us to factor $E\la Y_\lambda$ further
through an object $X_\lambda\in\genu G{}{a,b}$. We deduce
\be
\setcounter{enumi}{\value{enumiv}}
\item
The subcategory $\car\subset\ogenu G{}{a,b}$ is closed under coproducts.
\setcounter{enumiv}{\value{enumi}}
\ee

Next put 
\[
\cs=\ct^-_c\,,\quad\ca=\genu G{}{a,b},\quad\cc=\genu G{}{a,b},\quad\cx=\car,\quad\cz=\car.
\]
By the hypotheses of the lemma $\cs=\ct^-_c$ is triagulated and contains
$\cc$. By the definition of $\car$, any map $s\la x$ with $s\in\cs$
and $x\in\cx$ factors through an object $a\in\ca$, and any map
$s\la z$ with $s\in\cs$ and $z\in\cz$ factors through an object $c\in\cc$.
We are in a position to apply
\cite[Lemma~1.6 and Remark~1.7]{Neeman17}, which gives that
any map $s\la y$, with $s\in\cs$ and $y\in\cx*\cz$,
must factor through an object
$b\in\ca*\cc=\genu G{}{a,b}$. Thus
\be
\setcounter{enumi}{\value{enumiv}}
\item
The inclusion $\car*\car\subset\car$ holds.
\setcounter{enumiv}{\value{enumi}}
\ee
This completes the proof of the Lemma.
\eprf

\pro{PP1.309}
Suppose $\ct$ is a triangulated category, and assume it has a compact
generator $H$ with $\Hom(H,\T^iH)=0$ for $i\gg0$.
Let $\ct^b_c$ be the one corresponding to a preferred \tstr.

Suppose there
is an object $G\in\ct^b_c$ and an integer $N>0$ with $\ct=\ogenun GN$. Then
$\ct^b_c=\gen GN$.
\epro

\prf
Let $\big(\ct^{\leq0},\ct^{\geq0}\big)$ be a \tstr\ in the preferred
equivalence class, and choose it so that $\ct^{\geq0}$ is closed
under coproducts.

Take any object $A\in\ct^b_c$. Because $A$ lies in $\ct^b$
is must belong to $\ct^{\geq a}\cap\ct^{\leq b}$ for some integers
$a\leq b$, and Lemma~\ref{LL1.-27} guarantees that it must belong
to $\ogenu GN{a-M,b+M}$ for some $M>0$. But then the identity
$A\la A$ is a morphism from the object $A\in\ct^b_c\subset\ct^-_c$ to the
object $A\in \ogenu GN{a-M,b+M}$, and Lemma~\ref{LL1.10909} gives that it
must factor through an object $B\in\genu GN{a-M,b+M}\subset\gen GN$. As
$\gen GN$  is closed under direct summands it must contain $A$.
\eprf

\lem{LL1.30909}
Suppose $\ct$ is a
triangulated category with coporoducts, and assume
that $\ct$ has a compact generator $H$ with $\Hom(H,\T^iH)=0$
for $i\gg0$. Suppose further that there
is an object $G\in\ct^b_c$ and an integer $N>0$ with $\ct=\ogenun GN$.

Assume  $\big(\ct^{\leq0},\ct^{\geq0}\big)$ is a \tstr\
in the preferred equivalence
class. Then 
there exists an integer $M>0$ such that, for any object $Y\in\ct^-_c$
and any integer $i\in\zz$, the \tstr\
truncation morphism $Y\la Y^{\geq-i}$ may be factored as
$Y\la Y^{\geq-i-2M}\la C\la Y^{\geq-i}$ with
$C\in\genuf GN{-i-M}$.
\elem

\prf
Replacing the \tstr\ by an equivalent one if necessary, we may assume
$\ct^{\geq0}$ is closed under coproducts.
Choose an integer $M>0$
as in Lemma~\ref{LL1.-27}; in particular
$\ct^{\geq -i}\subset\ogenuf GN{-i-M}$. Assume further that
$M$ is large enough so that
$G\in\ct^b_c$ belongs to $\ct^{\geq -M}$.

We are given the map $Y\la Y^{\geq-i}$, with $Y\in\ct^-_c$ and  $Y^{\geq-i}$
in $\ct^{\geq-i}$, and by the choice of $M$ we have that $\ct^{\geq-i}$
is contained in $\ogenuf GN{-i-M}$.
Now Lemma~\ref{LL1.10909} permits
us to factor $Y\la Y^{\geq-i}$ through
an object $C\in\genuf GN{-i-M}$. This far we have composites
$Y\la C\la Y^{\geq-i}$.

But $G$ is an object of $\ct^{\geq-M}$, hence $\genuf GN{-i-M}$ must
be contained in $\ct^{\geq-i-2M}$. The map $Y\la C$ must therefore factor
canonically through the \tstr\ truncation, and we have our
factorization $Y\la Y^{\geq-i-2M}\la C\la Y^{\geq-i}$.
\eprf

\cor{CC1.310}
Suppose $\ct$ is a
triangulated category with coporoducts, and assume there
is an object $G\in\ct^b_c$ and an integer $N$ with $\ct=\ogenun GN$.
Assume further that $\ct$ has a compact generator $H$ with $\Hom(H,\T^iH)=0$
for $i\gg0$.
For any object $Y\in\ct^-_c$ we may choose an inverse sequence
$\cdots\la E_3\la E_2\la E_1\la E_0$
 so that the subsequence
$\cdots\la E_7\la E_5\la E_3\la E_1$ lies in $\ct^b_c$ while
the subsequence
$\cdots\la E_8\la E_6\la E_2\la E_0$
is a subsequence of
$\cdots\la Y^{\geq-3}\la  Y^{\geq-2}\la Y^{\geq -1}\la  Y^{\geq0}$.
\ecor

\prf
The construction of the sequence $E_i$ is just by iterating
Lemma~\ref{LL1.30909}.
\eprf

In view of Corollary~\ref{CC1.310}, the next lemma becomes interesting.

\lem{LL1.989}
Suppose $Y\in\ct^-_c$ is an object, mapping to
an inverse system 
$\cdots\la E_3\la E_2\la E_1\la E_0$ in $\ct^b_c$.
Assume this inverse system is pro-isomorphic to
$\cdots\la F_3\la F_2\la F_1\la F_0$, and what we know about $F_*$ is
that for any $n>0$ there exists an $m>0$ so that
the map $Y^{\geq-n}\la F_i^{\geq-n}$ is an isomorphism for all $i\geq m$.

When we view  $\Hom(Y,-)$ as a functor on $\ct^b_c$, it is
equal to the colimit of $\Hom(E_i,-)$.
\elem

\prf
Any object $Z\in\ct^b_c$ belongs to
$\ct^{\geq-n}$ for some $n>0$, and the sequence $\Hom(F_i,Z)$ becomes
stable and isomorphic to $\Hom(Y,Z)$ for $i\gg0$. Hence
$\Hom(Y,-)$ is the colimit of the sequence $\Hom(F_i,-)$. Therefore
it must also be the colimit of the ind-isomorphic sequence
$\Hom(E_i,-)$.
\eprf

\section{Combining the hypothesis $\ct=\ogenun G N$ with approximability}
\label{SS2}

The careful reader might have noticed that, in Section~\ref{SS1},
we didn't assume much about the triangulated category $\ct$.
Let us perhaps rephrase this more accurately: we made the
strong assumption that there exists an object $G\in\ct^b_c$
with $\ogenun GN=\ct$. But, beyond that, the only
assumption was that there exist
a compact generator $H\subset\ct$
satisfying $\Hom(H,\T^iH)=0$ for
$i\gg0$.

This is about to change, in this section we combine this
with the approximability hypotheses.

\ntn{NN2.21}
It is time to begin recalling the setup in Theorem~\ref{T1.-1}.
In the next lemmas $R$ will be a commutative ring and $\ct$ will be
an $R$--linear, weakly approximable triangulated category. Let
$\ct^b_c$ and $\ct^-_c$ be understood with respect to the
preferred equivalence class of {\it t}--structures. We will be
considering two functors
\[\xymatrix@C+40pt@R-30pt{
\ct \ar[r]^-\cy  & \Hom_R^{}\Big(\big[\ct^c\big]\op\,,\,\Mod R\Big)\\
\ct\op \ar[r]^-{\wt\cy}  & \Hom_R^{}\Big(\ct^b_c\,,\,\Mod R\Big)
}\]
The object $A\in\ct$ goes under the functors, respectively, to
\[
\cy(A)=\Hom(-,A),\qquad\wt\cy(A)=\Hom(A,-).
\]
For $\cy(A)$ the variable $(-)$ takes its values in $\ct^c$, while in
the case of $\wt\cy(A)$ the variable $(-)$ lies in $\ct^b_c$.

We will mostly be concerned with the restrictions of $\cy$ and $\wt\cy$
to the subcategory $\ct^-_c$.
\entn

\pro{PP2.5}
Let the conventions be as in Notation~\ref{NN2.21}, and assume there
is an object $G\in\ct^b_c$ and an integer $N$ with $\ct=\ogenun GN$.

Then restriction to $\ct^-_c$ of the map $\wt\cy$ is full. More
generally: any morphism $\ph:\wt\cy(b)\la\wt\cy(a)$, with $b\in\ct^-_c$ and
$a\in\ct$, is equal to $\wt\cy(f)$ for some $f:a\la b$.
\epro

\prf
 By Corollary~\ref{CC1.310} we may
construct for $b$  an inverse sequence
$\cdots\la C_3\la C_2\la C_1\la C_0$ pro-isomorphic to the sequence
$\cdots\la b^{\geq-3}\la  b^{\geq-2}\la b^{\geq -1}\la  b^{\geq0}$, and do it
in
such a way that $C_i$ all belong to $\ct^b_c$. Lemma~\ref{LL1.989} tells
us, moreover, that 
$\wt\cy(b)$ is the colimit
of $\wt\cy(C_i)$.

But $\wt\cy(C_i)$ is representable, hence the composite
$\wt\cy(C_i)\la\wt\cy(b)\stackrel\ph\la\wt\cy(a)$ is a morphism
from a representable
functor. Yoneda tells us that it must be $\wt\cy$ of a unique map
$f_i:a\la C_i$. These maps are compatible, and hence all factor
through a morphism $f:a\la z$ with $z=\holim C_i$. Because
the sequence $C_i$ is pro-isomorphic to the sequence
$b^{\geq-i}$ we have that $z=\holim b^{\geq-i}$, and Proposition~\ref{PP2.95}
gives an isomorphism $b\cong z$. We have produced a morphism
$f:a\la b$, and it's now obvious that $\wt\cy(f)=\ph$.
\eprf

\lem{LL2.21.5}
Let the conventions be as in Notation~\ref{NN2.21}, and assume there
is an object $G\in\ct^b_c$ and an integer $N$ with $\ct=\ogenun GN$.

Given two morphisms $f,g:X\la Y$ in the category $\ct^-_c$, we have
$\cy(f)=\cy(g)$ if and only if $\wt\cy(f)=\wt\cy(g)$.
\elem

\prf
Choose and fix a \tstr\ $\big(\ct^{\leq0},\ct^{\geq0}\big)$
in the preferred equivalence class.

Suppose $\cy(f)=\cy(g)$, let $S\in\ct^b_c$ be an
object, and let $h:Y\la S$ be an element of
$\wt\cy(Y)(S)=\Hom(Y,S)$. Because
$S\in\ct^b$ we may choose an integer $m>0$ with $S\in\ct^{\geq-m+1}$.
Now the object $X$ belongs to $\ct^-_c$, hence there is a triangle
$E\stackrel\alpha\la X\la D$ with $E\in\ct^c$ and $D\in\ct^{\leq-m}$.
Because $\cy(f)=\cy(g)$ the two composites
\[\xymatrix@C+30pt{
E\ar[r]^-\alpha & X \ar@<0.5ex>[r]^-{f} \ar@<-0.5ex>[r]_-{g} & Y
}\]
must be equal. Hence so are the longer composites
\[\xymatrix@C+30pt{
E\ar[r]^-\alpha & X \ar@<0.5ex>[r]^-{f} \ar@<-0.5ex>[r]_-{g} & Y
\ar[r]^-h & S
}\]
Therefore the map $h(f-g):X\la S$ must factor as $X\la D\la S$, but
as $D\in\ct^{\leq-m}$ and $S\in\ct^{\geq-m+1}$ the map $D\la S$ must
vanish. Therefore $hf=hg$. Since this is true for every $h$ we have
$\wt\cy(f)=\wt\cy(g)$.

Next suppose $\wt\cy(f)=\wt\cy(g)$, let $E\in\ct^c$ be an
object, and let $\alpha:E\la X$ be an
element in $\cy(X)(E)=\Hom(E,X)$. By
Lemma~\ref{L17.17} we may choose an integer $m>0$ with
$\Hom\big(E,\ct^{\leq-m}\big)=0$. By Lemma~\ref{LL1.30909} the
map $Y\la Y^{\geq-m+1}$ factors as $Y\stackrel h\la S\la Y^{\geq-m+1}$
with $S\in\ct^b_c$. Now
the two composites
\[\xymatrix@C+30pt{
X \ar@<0.5ex>[r]^-{f} \ar@<-0.5ex>[r]_-{g} & Y
\ar[r]^-h & S
}\]
are equal because $\wt\cy(f)=\wt\cy(g)$, hence the longer
composites
\[\xymatrix@C+30pt{
E\ar[r]^-\alpha & X \ar@<0.5ex>[r]^-{f} \ar@<-0.5ex>[r]_-{g} & Y
\ar[r]^-h & S\ar[r] & Y^{\geq-m+1}
}\]
must also be equal. In other words: the map $Y\la Y^{\geq-m+1}$ annihilates
the map $(f-g)\alpha:E\la Y$, and hence $(f-g)\alpha$ must factor as
$E\la Y^{\leq-m}\la Y$. But $\Hom\big(E,\ct^{\leq-m}\big)=0$, and we
deduce that $f\alpha=g\alpha$. As this is true for every $\alpha$ we have
$\cy(f)=\cy(g)$.
\eprf

Because the symbols $G$ and $G'$ are similar and easily confused with
each other, in Section~\ref{SS1} we
tried to avoid using both in any one statement. But in the next
lemma we cannot escape the conundrum; the letter $H$ it taken, it
means a homological functor.

\lem{LL2.22}
Let $R$ be a noetherian ring.
Let the conventions be as in  Notation~\ref{NN2.21},
except that weak approximability is no longer enough---for this
Lemma we assume $\ct$ approximable.
Suppose also that there
is an object $G\in\ct^b_c$ and an integer $N$ with $\ct=\ogenun GN$.
Assume further that there is a compact generator $G'\in\ct$
with $\Hom(G',\T^iG')$ a finite $R$--module for every $i\in\zz$.

Suppose we are given an $R$--linear cohomological functor
$H:\ct^b_c\la \Mod R$, as well as an object $Y\in\ct^-_c$ such that
$H$ is a direct smmand of $\wt\cy(Y)$.
Then there exists an object $Y'\in\ct^-_c$ and an isomorphism
$H\cong\wt\cy(Y')$.
\elem

\prf
We are given that $H$ is a direct summand of $\wt\cy(Y)$, and 
the composite $\wt\cy(Y)\la H\la\wt\cy(Y)$ is an idempotent endomorphism
$\ph:\wt\cy(Y)\la\wt\cy(Y)$. By Proposition~\ref{PP2.5} there is a morphism
$f:Y\la Y$ in the category $\ct^-_c$ with $\wt\cy(f)=\ph$.
Because $\wt\cy(f)$ is idempotent we have that $\wt\cy(f)=\wt\cy(f^2)$, and
Lemma~\ref{LL2.21.5} informs us that $\cy(f)=\cy(f^2)$. Therefore
$\cy(f)$ is idempotent, and corresponds to the projection to
a direct summand of $\cy(Y)$. The finiteness hypotheses on $\Hom(G',\T^iG')$
guarantee that $\cy(Y)$  is
a locally finite cohomological functor, hence so is any direct summand,
and by Theorem~\ref{T37.2068} there
exists an object $Y'\in\ct^-_c$ and morphisms
$Y\stackrel\alpha\la Y'\stackrel\beta\la Y$ with
$\cy(f)=\cy(\beta\alpha)$ and $\cy(\alpha\beta)=\id_{\cy(Y')}^{}=\cy(\id_{Y'}^{})$.

Lemma~\ref{LL2.21.5} informs us that $\wt\cy(f)=\wt\cy(\beta\alpha)$ and
$\wt\cy(\alpha\beta)=\wt\cy(\id_{Y'}^{})$. Thus $\wt\cy(f)$ factors as
\[\xymatrix@C+40pt{
\wt\cy(Y)  \ar[r]^{\wt\cy(\alpha)} &   \wt\cy(Y') \ar[r]^-{ \wt\cy(\beta)} &
\wt\cy(Y)
}\]
while
\[\xymatrix@C+40pt{
\wt\cy(Y')  \ar[r]^{\wt\cy(\beta)} &   \wt\cy(Y) \ar[r]^-{ \wt\cy(\alpha)} &
\wt\cy(Y')
}\]
composes to the identity.
The current lemma follows. 
\eprf

\section{The proof of Theorem~{\protect{\ref{T1.-1}}}($\mathrm{ii}$)}
\label{SS3}

The proof of the main theorems will be by applying to the current situation
the lemmas
of Section~\ref{S997}. More precisely

\ntn{NN3.1}
With the conventions of Notation~\ref{NN2.21},
assume given an
object $G\in\ct^b_c$ and an integer $N>0$ such that $\ct=\ogenun GN$.
Assume also that we have fixed a
 \tstr\ $(\ct^{\leq0},\ct^{\geq0})$ in the preferred
 equivalence class, with $\ct^{\geq0}$ closed under
 coproducts.
Let $\ca$ be the heart of the \tstr, and
$\ch:\ct\la\ca$ the standard homological functor.

The lemmas
of Section~\ref{S997} will be applied to the category $\ct\op$.
The subcategory $\cs$ of Notation~\ref{N997.1} will be
$\cs=\big[\ct^b_c\big]\op$.
\entn

The next definition and lemma are similar to what works
in Section~\ref{S37}. The
reader might wish to compare the definition below, and the lemma that
follows, with Definition~\ref{D1.29} and
Lemma~\ref{L1.1010101}.

\dfn{D3.3}
A \emph{powerful $\gen Gn$--approximating sequence} is an inverse system
$\cdots\la E_3\la E_2\la E_1\la E_0$
in $\ct$, so that
\be
\item
  Each $E_m$ belongs to $\gen Gn$.
\item
The map $\ch^i(E_{m+1})\la\ch^i(E_{m})$ is an isomorphism whenever $i\geq-m$.
\setcounter{enumiv}{\value{enumi}}
\ee
Suppose we are also given an object $F\in\ct$, together with
\be
\setcounter{enumi}{\value{enumiv}}
\item
A map from $F$ to
the approximating system $E_*$.
\item
The map in (iii) is such that
$\ch^i(F)\la\ch^i(E_m)$ is an isomorphism whenever $i\geq-m$.
\ee 
Then we declare $E_*$ to be a
\emph{powerful $\gen Gn$--approximating system for $F$.}
\edfn

\lem{LL3.5}
With the conventions of Definition~\ref{D3.3} we have
\be
\item
Given an object $F\in\ct$ and a powerful $\gen Gn$--approximating
system $E_*$ for $F$, then the (non-canonical) map $F\la\holim E_i$
is an isomorphism.
\item
Any $\gen Gn$--powerful approximating system
$\cdots\la E_3\la E_2\la E_1\la E_0$
has a subsequence which is a
powerful $\gen Gn$--approximating system of the
homotopy limit $F=\holim E_i$. Moreover
$F$ belongs to $\ct^-_c$.
\setcounter{enumiv}{\value{enumi}}
\ee
\elem

\prf
Assertion (i) is contained in Proposition~\ref{PP2.95}, and the
``moreover'' part of (ii) is contained in Lemma~\ref{LL29.-1}.

Let $L>0$ be the integer of Lemma~\ref{LL2.1}. Then Lemma~\ref{LL2.1}
says that, for any $\gen Gn$--powerful approximating system
$\cdots\la E_3\la E_2\la E_1\la E_0$ and with $F=\holim E_m$, the
maps $\ch^i(F)\la\ch^i(E_{m+L}^{})$ are isomorphisms whenever $i\geq-m$.
In other words the 
subsequence $\cdots\la E_{3+L}^{}\la E_{2+L}^{}\la E_{1+L}^{}\la E_{L}^{}$ is a
powerful approximating sequence for $F$. This proves the first
half of (ii).
\eprf

\rmk{RR3.9}
Let us now explain how
to specialize Lemma~\ref{L1.35} to the framework of this section.
Suppose we are given an object $\wh B\in\ct^-_c$ and a powerful
$\gen G{n}$--approximating
system $\mathfrak{B}_*$ for $\wh B$. 
Lemma~\ref{LL1.989} informs us that
the natural map $\colim\,\wt\cy(\mathfrak{B}_i)\la\wt\cy(\wh B)$
is an isomorphism.
Thus
\be
\item
A powerful
$\gen G{n}$--approximating
system $\mathfrak{B}_*$ for $\wh B$ is an approximating system
for $\wt\cy(\wh B)$ in the sense of
Definition~\ref{D27.109}.
Moreover: Lemma~\ref{LL3.5}(i) tells us that the
map $\wh B\la\holim\mathfrak{B}_m$ is an isomorphism,
hence a powerful $\gen Gn$--approximating system for $\wh B$, as
in Definition~\ref{D3.3}, is also an approximating system for
$\wh B$ as in Remark~\ref{R27.7}.
\setcounter{enumiv}{\value{enumi}}
\ee
Note also that in Lemma~\ref{LL3.5}(ii) we learned that any powerful
$\gen G{n}$--approximating
system $\mathfrak{B}_*$ has a subsequence which is a powerful
$\gen G{n}$--approximating
system of $\holim\mathfrak{B}_i$.

Next assume we are given
\be
\setcounter{enumi}{\value{enumiv}}
\item
A morphism $\wh\beta:\wh B\la \wh C$ in the category $\ct^-_c$.
\item
Two integers $n$ and $n'$, as well as a powerful $\gen G{n}$--approximating
system $\mathfrak{B}_*$ for $\wh B$ and a
powerful $\gen G{n'}$--approximating
system $\mathfrak{C}_*$ for $\wh C$.
\setcounter{enumiv}{\value{enumi}}
\ee
The dual of Lemma~\ref{L27.592} allows us to
choose a subsequence of
$\mathfrak{B'}_*\subset\mathfrak{B}_*$ and a map of sequences
$\beta_*:\mathfrak{B'}_*\la\mathfrak{C}_*$ compatible with
$\wh\beta:\wh B\la \wh C$.
A subsequence of a powerful $\gen G{n}$--approximating sequence is
clearly a powerful $\gen G{n}$--approximating sequence,
hence $\mathfrak{B'}_*$ is a powerful $\gen G{n}$--approximating sequence for
$\wh B$.
Now 
as in Lemma~\ref{L1.35} we extend
$\beta_*:\mathfrak{B'}_*\la\mathfrak{C}_*$ to a sequence of
triangles, in particular
for each $m>0$ this gives a morphism of triangles
\[\xymatrix@C+20pt{
 \Tm \mathfrak{B'}_{m+1}\ar[d]\ar[r]^{\Tm\beta_{m+1}}&
  \Tm\mathfrak{C}_{m+1}\ar[d]\ar[r]^{\Tm\gamma_{m+1}} &
  \mathfrak{A}_{m+1}\ar[d]\ar[r]^{\alpha_{m+1}}&
  \mathfrak{B'}_{m+1}\ar[d]\ar[r]^{\beta_{m+1}}&
  \mathfrak{C}_{m+1}\ar[d] \\
   \Tm \mathfrak{B'}_m\ar[r]^{\Tm\beta_m}
    & \Tm\mathfrak{C}_m\ar[r]^{\Tm\gamma_m} &
  \mathfrak{A}_m\ar[r]^{\alpha_m} & \mathfrak{B'}_m\ar[r]^{\beta_m}
  & \mathfrak{C}_m 
}\]
Applying the functor $\ch^i$ with $i\geq-m+1$ yields a
commutative diagram in the heart of $\ct$
where the rows are exact, and where
vertical maps away from the middle are isomorphisms.
By the 5-lemma the middle vertical map, i.e.~the map
$\ch^i(\mathfrak{A}_{m+1})\la\ch^i(\mathfrak{A}_{m})$, must also be an
isomorphism when $i\geq-m+1$.
We conclude that a subsequence of $\mathfrak{A}_*$ is a powerful
$\gen G{n'+n}$--approximating system. Put $\wh A=\holim\mathfrak{A}_*$.
By Lemma~\ref{LL3.5}(ii) the object $\wh A$ belongs to $\ct^-_c$, 
and Lemma~\ref{LL1.989} coupled with
Proposition~\ref{PP2.5} guarantee that the weak triangle
$A\stackrel u\la B\stackrel v\la C\stackrel w\la\T A$
of Lemma~\ref{L1.35} is isomorphic
to the image under $\wt\cy$ of a weak triangle
$\wh A\stackrel{\wh u}\la \wh B\stackrel{\wh v}\la\wh C\stackrel{\wh w}\la\T \wh A$
in the category $\ct^-_c$.
\ermk

\lem{LL3.11}
Suppose $H$ is a locally finite $\gen G1$--homological functor. Then there is
a surjection $\wt\cy(F)|_{\gen G1}\la H$, where $F\in\ct^-_c$ has a powerful
$\gen G1$--approximating system $\cdots\la E_3\la E_2\la E_1\la E_0$.
Moreover: the system may be chosen so that the maps $E_{m+1}\la E_m$ are
split epimorphisms.
\elem

\prf
We have that $H(\T^i G)$ is a finite $R$--module
for every
$i\in\zz$, and vanishes when $i\ll0$. For each $i$ with $H(\T^iG)\neq0$
choose a
finite number of generators $\{f_{ij},\,j\in J_i\}$
for the $R$--module $H(\T^iG)$. By
Yoneda
every $f_{ij}\in H(\T^iG)$ corresponds
to a morphism $\ph_{ij}^{}:\wt\cy(\T^i G)\la H$. Let $F$ be
defined by
\[
F \eq\coprod_{i\in\zz}\bigoplus_{j\in J_i}\T^i G\quad\cong\quad\prod_{i\in\zz}\bigoplus_{j\in J_i}\T^i G
\]
where the isomorphism of the coproduct and the product is by Lemma~\ref{LL1.3}.
Let the morphism $\ph:\wt\cy(F)\la H$
be given by
\[\xymatrix@C+20pt{
  \wt\cy(F)|_{\gen G1}\ar@{=}[r]&
  \ds\bigoplus_{i\in\zz}\bigoplus_{j\in J_i}\wt\cy(\T^iG)
  \ar[rr]^-{(\ph_{ij})} && H
}\]
where $(\ph_{ij})$ stands for the row matrix with entries $\ph_{ij}$; on
the $i,j$ summand the map is $\ph_{ij}$.
Finally: because $G\in\ct^b_c$ there is an integer $B>0$ with
$\T^BG\in\ct^{\leq-1}$. For $m>0$ we define
\[
E_m=\bigoplus_{i< m+B}\,\bigoplus_{j\in J_i}\T^i G
\]
The sum is finite by hypothesis, making
$E_m$ an object of $\gen G1$. The obvious map $E_{m+1}\la E_m$
is a split epimorphism, and in the decomposition $ F\cong E_m\oplus \wt F$ we
have that $\wt F$, being the coproduct of $\T^iG$ for $i\geq m+B$,
belongs to $\ct^{\leq-m-1}$. Therefore the map $\ch^i(F)\la\ch^i(E_m)$ is an
isomorphism if $i\geq-m$, making the $E_*$ a powerful
$\gen G1$--approximating system
for $F$.
\eprf

And now the time has come to prove the main results.

\thm{TT3.13}
Let $R$ be a noetherian, commutative ring.
Let $\ct$ be an $R$--linear, approximable triangulated category, and assume
there is a compact generator $G'\in\ct$ such that $\Hom(G',\T^iG')$ is
a finite $R$--module for all $i\in\zz$.
Let
$\ct^-_c$ and $\ct^b_c$ be the ones corresponding to the preferred
equivalence
class of {\it t}--structures, and assume there is an object $G\in\ct^b_c$
and an integer $N>0$ with $\ct=\ogenun GN$.

Then the functor $\wt\cy:\ct^-_c\la\Hom\big[\ct^b_c\,,\,\Mod R\big]$
if full, and the essential image consists
of the locally finite homological functors.
\ethm

\prf
The fact that the functor $\wt\cy$ is full was proved in
Proposition~\ref{PP2.5}, and the fact that its image is contained
in the locally finite homological functors was shown in
Lemma~\ref{LL3.-1}. What needs proof is that every locally
finite homological functor can be realized as $\wt\cy(F)$ for
some $F\in\ct^-_c$.

Suppose therefore that $H$ is a locally finite $\ct^b_c$--homological functor.
Therefore $H|_{\gen G1}$ is a locally finite $\gen G1$--homological functor,
and Lemma~\ref{LL3.11} produces an object $F_1\in\ct^-_c$, with
a powerful $\gen G1$--approximating system, and an
epimorphism $\wt\cy(F_1)|_{\gen G1}\la H|_{\gen G1}$. From
Corollary~\ref{C27.5} if follows that we may lift the natural
transformation to all of $\ct^b_c$; there is a natural transformation
$\ph_1^{}:\wt\cy(F_1)\la H$ so that
$\ph_1^{}|_{\gen G1}:\wt\cy(F_1)|_{\gen G1}\la H|_{\gen G1}$ is surjective.

Next we proceed inductively. Suppose we have constructed
$F_n\in\ct^-_c$, with a powerful $\gen Gn$--approximating system,
and a natural transformation $\ph_n^{}:\wt\cy(F_n)\la H$,
and assume
that $\ph_n^{}|_{\gen G1}:\wt\cy(F_n)|_{\gen G1}\la H|_{\gen G1}$ is surjective.
Since both $H|_{\gen G1}$ and $\wt\cy(F_n)|_{\gen G1}$ are locally finite
and the ring $R$ is noetherian, the kernel of $\ph_n^{}|_{\gen G1}$ is also
locally finite. Lemma~\ref{LL3.11} permits us to find a surjection to the
kernel: there is an object
$F'\in\ct^-_c$, with a powerful $\gen G1$--approximating system
$\cdots\la E_3\la E_2\la E_1\la E_0$
in which
all the connecting maps $E_{m+1}\la E_m$ are split epimorphisms,
and an exact sequence
$\wt\cy(F')|_{\gen G1}\la\wt\cy(F_n)|_{\gen G1}\la H|_{\gen G1}$.
Now Corollary~\ref{C27.5} allows us to lift
the map to $\ct^-_c$. We deduce:
\be
\item
There is a morphism $\alpha:F_n\la F'$ so that
the sequence below is exact
\[\xymatrix@C+40pt{
\wt\cy(F')|_{\gen G1}\ar[r]^-{\wt\cy(\alpha)|_{\gen G1}} &\wt\cy(F_n)|_{\gen G1}
\ar[r]^-{\ph_n|_{\gen G1}} & H|_{\gen G1}
}\]
\setcounter{enumiv}{\value{enumi}}
\ee
Forget for a second the exactness; the vanishing
of the composite in (i), coupled
with Lemma~\ref{L1.35}, allows us to construct
\be
\setcounter{enumi}{\value{enumiv}}
\item
With the notation of Definition~\ref{D997.203},
and working in the category $\ct\op$ and with $\cs=\big[\ct^b_c\big]\op$,
there is an object $F_{n+1}\in\ct^-_c$ with a powerful
$\gen G{n+1}$--approximating system, a weak triangle in $\ct^-_c$ of the form
$F_{n+1}\stackrel{\beta_n}\la F_n\stackrel{\alpha}\la F'\la\T F_{n+1}$,
and a morphism $\ph_{n+1}:\wt\cy(F_{n+1})\la H$ so that the
following triangle commutes
\[\xymatrix@C+40pt@R-20pt{
                     & \wt\cy(F_{n+1})\ar[dd]^-{\ph_{n+1}} \\
\wt\cy(F_n)\ar[ru]^-{\wt\cy(\beta_n)}\ar[rd]_-{\ph_n} & \\
                     & H
}\]
\ee
This inductively constructs an inverse sequence in $\ct^-_c$ of the
form $\cdots \la F_4\stackrel{\beta_3^{}}\la F_3\stackrel{\beta_2^{}}\la F_2\stackrel{\beta_1^{}}
\la F_1$,
as well as compatible maps
$\ph_n:\wt\cy(F_n)\la H$.

Now the map $\ph_1^{}:\wt\cy(F_1)\la H$ restricts to an epimorphism
on $\gen G1$ by construction,
and the exactness of the sequence in (i) coupled with
Lemma~\ref{L997.601} informs us, by induction, that
$\ph_n:\wt\cy(F_n)\la H$ restricts to an epimorphism on $\gen Gn$.
By Proposition~\ref{PP1.309} we have $\ct^b_c=\gen GN$, hence
$\ph_N^{}$ is an epimorphism.

Now apply Lemma~\ref{L997.223} to the diagram
\[\xymatrix@C+20pt{
                  & \wt\cy(F_{N})\ar[d]  \\
\wt\cy(F')\ar[r]^-{\wt\cy(\alpha)} &\wt\cy(F_{N+n}) \ar[r]^-{\wt\cy(\beta_{N+n}^{})} & \wt\cy(F_{N+n+1})\ar[r]\ar[d]^-{\ph_{n+N+1}} & \wt\cy(\Tm F')\\
  & & H
}\]
Induction on $n\geq0$ teaches us that the map
$\wt\cy(F_N)|_{\gen Gn}\la\wt\cy(F_{N+n})|_{\gen Gn}$ annihilates the
kernel of $\wt\cy(F_N)|_{\gen Gn}\la H|_{\gen Gn}$. If we put $n=N$ and
remember that $\gen GN=\ct^b_c$, we have that the map
$\wt\cy(F_N)\la\wt\cy(F_{2N})$ and the epimorphism
$\wt\cy(F_N)\la H$ have the same kernel. Thus $\wt\cy(F_N)\la\wt\cy(F_{2N})\la H$
factors as $\wt\cy(F_N)\la H\la\wt\cy(F_{2N})\la H$,
making $H$ is a direct
summand of $\wt\cy(F_{2N})$. 
Lemma~\ref{LL2.22} produces an object $Y\in\ct^-_c$
with $H=\wt\cy(Y)$.
\eprf

\thm{TT3.15}
Let the notation be as in Theorem~\ref{TT3.13}. The essential image under
$\wt\cy$ of the subcategory $\ct^c\subset\ct^-_c$ is precisely the
finite homological functors. Moreover the restriction of $\wt\cy$
to $\ct^c$ is fully faithful: it induces an
equivalence of $\big[\ct^c\big]\op$ with the category of finite
homological functors $\ct^b_c\la\mod R$.

The ``moreover'' part can even be strengthened as follows: for any
pair of objects $a\in\ct^c$ and $b\in\ct^-_c$ the natural map is an
isomorphism
\[\xymatrix@C+40pt{
\Hom(a,b)\ar[r] & \Hom\big[\wt\cy(b),\wt\cy(a)\big]
}\]
\ethm

\prf
The fact that $\wt\cy(A)$ is finite when $A\in\ct^c$ follows
from Lemma~\ref{LL3.-1}---the essential image under $\wt\cy$ of
the subcategory $\ct^c\subset\ct^-_c$ is contained in the finite
functors.

Now suppose $H:\ct^b_c\la\mod R$ is a finite homological functor.
Since finite homological functors
are locally finite Theorem~\ref{TT3.13} tells us that there exists
an object $A\in\ct^-_c$ and an isomorphism $\wt\cy(A)\cong H$. It
suffices to prove that $A\in\ct^c$. But the finiteness tells us that,
for the object $G\in\ct^b_c$ with $\ct=\ogenun GN$ of the
hypotheses of Theorem~\ref{TT3.13}, we must have
that $H^{i}(G)\cong\Hom(A,\T^iG)=0$ for $i\gg0$. By Proposition~\ref{PP1.2000}
$A$ must be compact.

It remains to prove the full faithfulness, or rather the strengthened version.
We already know the surjectivity of the map
\[\xymatrix@C+40pt{
\Hom(a,b)\ar[r] & \Hom\big[\wt\cy(b),\wt\cy(a)\big]\ ,
}\]
that was part of  Theorem~\ref{TT3.13}. Suppose therefore that we have two
morphisms $f,g:a\la b$ with $\wt\cy(f)=\wt\cy(g)$. Then Lemma~\ref{LL2.21.5}
informs us that $\cy(f)=\cy(g)$, and as $\cy(a)=\Hom(a,-)$ is representable
we deduce from Yoneda that $f=g$.
\eprf

\section{Applications: the construction of adjoints}
\label{S73}

We prove Corollary~\ref{C0.-978}, a restricted version of
which was the key tool in Jack Hall's
original, simple proof of
GAGA---see Remark~\ref{R1.-8}. Hall's later proofs of more general results,
see~\cite{Hall18,Hall22},
sidestep the representability theorems presented here.

\thm{T73.1}
Let $R$ be a noetherian, commutative ring.
Let $\ct$ be an $R$--linear triangulated category with coproducts, and assume
that it is approximable. Let $\ct^b_c\subset\ct^-_c$
be the subcategories of Definition~\ref{D0.13}, constructed
using a \tstr\ in the preferred equivalence class.
Assume the category $\ct^c$ is contained in $\ct^b_c$.
Assume further that $\ct$  has a compact
generator $G$ so that $\Hom_\ct^{}(-,G)$ is a $G$--locally finite
cohomological functor.

Let $\cl:\ct^b_c\la\cs$ be an $R$--linear triangulated functor, and
let $(\ct^{\leq0},\ct^{\geq0})$ be any \tstr\ in the preferred equivalence class. Then the functor $\cl$ has a right adjoint
if and only if the following three conditions hold:
\be
\item
  For any pair of objects $(t,s)$, with $t\in\ct^c$
  and $s\in\cs$, the $R$--module
  $\Hom\big(\cl(t),s\big)$ is finite.
\item
  For any object $s\in\cs$ there exists an integer $A>0$ with
  $\Hom\big(\cl(\ct^b_c\cap\ct^{\leq-A})\,,\,s\big)=0$.
\item
  For any object $t\in\ct^c$ and any object $s\in\cs$ there exists an integer
  $A$ so that $\Hom\big(\cl(\T^mt),s\big)=0$ for all $m\leq-A$.
\ee
\ethm

\prf
We begin by proving the necessity of the three conditions.
Assume that $\cl$ has a right adjoint
$\car:\cs\la\ct^b_c$, and we will prove that (i), (ii) and (iii) hold.

For any object
$s\in\cs$ we have that $\car(s)\in\ct^b_c$. If $t\in\ct^c$
is any object, then $\Hom_{\ct}^{}\big(t,\car(s)\big)$
is a finite $R$--module by Lemma~\ref{L1.1.3}.
The isomorphism
\[ 
\Hom_{\ct}^{}\big(t,\car(s)\big)\cong\Hom_{\cs}^{}\big(\cl(t),s\big)
\]
informs us that $\Hom_{\cs}^{}\big(\cl(t),s\big)$ is also a finite
$R$--module, giving (i).

With $s$ still an object of $\cs$,
we have that $\car(s)\in\ct^b_c\subset\ct^b$. There must
therefore 
exist an integer $A>0$ with $\car(s)\in\ct^{\geq-A}\cap\ct^{\leq A}$.
The fact that $\car(s)$ belongs to $\ct^{\geq-A}$ tells us that
$\Hom_{\ct}^{}\big(-,\car(s)\big)$ annihilates $\ct^{\leq-A-1}$, and hence its subcategory $\ct^b_c\cap\ct^{\leq-A-1}$. But
by adjunction $\Hom(-,s)$ must annihilate
$\cl\big(\ct^b_c\cap\ct^{\leq-A-1}\big)$,
proving (ii).

Now let $t\in\ct^c$ be an object.
Remark~\ref{R0.14.7} combined with Lemma~\ref{L17.17}
allow us to choose an integer
$B>0$ with $\Hom\big(t,\T^B\ct^{\leq0}\big)=0$.
But now $\T^m\car(s)$ belongs to $\ct^{\leq 0}$ for all
$m\geq A$, and hence
$\Hom\big(t,\T^{m+B}\car(s)\big)=0$ for all $m\geq A$.
And now the isomorphism
\[
\Hom_{\ct}^{}\big(t\,,\,\T^{m+B}\car(s)\big)\cong
\Hom_{\cs}^{}\big(\cl(\T^{-m-B}t)\,,\,s\big)
\]
gives the vanishing of 
$\Hom_{\cs}^{}\big(\cl(\T^{-m-B}t)\,,\,s\big)$ for all
$m\geq A$,
completing the proof of (iii).

Now for the sufficiency: assume the three conditions
hold, and we need to produce the right adjoint $\car$.
For any pair of objects $t\in\ct^c$, $s\in\cs$ and any integer $m\in\zz$, from
(i) we learn that $\Hom\big(\cl(\T^mt),s\big)$ is a finite $R$--module.
Now (ii) and (iii) guarantee that it vanishes
whenever $m\gg0$ or $m\ll0$. Thus $\Hom\big(\cl(-),s\big)$ is a
finite $\ct^c$--cohomological functor. The assignment
taking $s\in\cs$ to the functor $\Hom\big(\cl(-),s\big)$ is
a functor from $\cs$ to the category of finite $\ct^c$--cohomological
functors; by Theorem~\ref{T37.1067} we can lift it through the
equivalence of categories $\cy$.
There is a functor $\car:\cs\la\ct^b_c$ so that,
for all objects $t\in\ct^c$ and all objects $s\in\cs$, we have
a natural isomorphism
\[\xymatrix@C+20pt{
\Hom\big(\cl(t),s\big)\ar[r]^-{\ph} & \Hom\big(t,\car(s)\big)\ .
}\]
Fix $t'\in\ct^b_c$ and consider the following composite, which is
natural in $t\in\ct^c$, $t'\in\ct^b_c$
\[\xymatrix@C+20pt{
  \Hom(t,t')\ar[r]^-{\cl}&\Hom\big(\cl(t),\cl(t')\big)\ar[r]^-{\ph}
  &\Hom\big(t,\car\cl(t')\big)\ .
}\]
We have objects $t',\car\cl(t')\in\ct^b_c$ and a natural
transformation $\cy(t')\la \cy\big(\car\cl(t')\big)$,
and Theorem~\ref{T37.2068} allows us to express it as
$\cy(\alpha_{t'}^{})$ for a unique morphism $\alpha_{t'}^{}:t'\la \car\cl(t')$
in the category $\ct^b_c$. We leave it to the reader to check that $\alpha_{t'}^{}$
is natural in $t'$; it gives a natural transformation $\alpha:\id\la \car\cl$.

For a general $t'\in\ct^b_c$ we have no idea how to compute
$\alpha_{t'}^{}$; it is a mysterious morphism
in $\ct^b_c$ that comes from a representability
theorem applied to $\cy(t')\la\cy\big(\car\cl(t')\big)$.
However: when $t'$ belongs to $\ct^c\subset\ct^b_c$
then $\cy(t')=\Hom(-,t')$ is a representable
functor on $\ct^c$, and Yoneda's lemma tells us how
to compute $\alpha_{t'}^{}$; it is the image of
$1\in\Hom(t',t')$ under
$\cy(t')\la\cy\big(\car\cl(t')\big)$. But this is
explicit: it is the image of $1\in\Hom(t',t')$
under the composite
\[\xymatrix@C+20pt{
  \Hom(t',t')\ar[r]^-{\cl}&\Hom\big(\cl(t'),\cl(t')\big)\ar[r]^-{\ph}
  &\Hom\big(t',\car\cl(t')\big)\ ,
}\]
which agrees with the image of
$1:\cl(t')\la\cl(t')$ under the map
$\ph:\Hom\big(\cl(t'),\cl(t')\big)\la\Hom\big(t',\car\cl(t')\big)$.

Now we define a natural
transformation $\psi:\Hom\big(\cl(-),-\big)\la \Hom\big(-,\car(-)\big)$.
For objects $t\in\ct^b_c$, $s\in\cs$ the map is
\[\xymatrix@C+10pt{
  \Hom\big(\cl(t),s\big)\ar[r]^-{\car}&\Hom\big(\car\cl(t),\car(s)\big)
  \ar[rr]^-{\Hom\big(\alpha_t^{},\car(s)\big)}
  &&\Hom\big(t,\car(s)\big)\ .
}\]
Once again: for a general $t\in\ct^b_c$ we don't understand
this map, but when we restrict $t$ to lie in
$\ct^c\subset\ct^b_c$ then we can compute. If
$\beta:\cl(t)\la s$ is a morphism in $\cs$, then
the fact that
$\ph$ is a natural transformation
gives the commutativity of the square
\[\xymatrix@C+40pt{
\Hom\big(\cl(t),\cl(t)\big)
  \ar[r]^-{\Hom\big(\cl(t),\beta\big)}\ar[d]_\ph & 
\Hom\big(\cl(t),s\big)\ar[d]^\ph \\
\Hom\big(t,\car\cl(t)\big)
  \ar[r]^-{\Hom\big(t,\car(\beta)\big)} & 
\Hom\big(t,\car(s)\big) \\
}\]
and, computing the image of $1:\cl(t)\la\cl(t)$ under
the two equal composites, we obtain the first equality
in
\[
\ph(\beta)=\car(\beta)\circ\alpha_t^{}=\psi(\beta)\ ,
\]
where the second equality is by the definition of
the natural transformation $\psi$.

Hence
when restricted to $t\in\ct^c\subset\ct^b_c$ the map $\psi$ agrees with
$\ph$ and is an isomorphism. It suffices to prove that $\psi$ is
an isomorphism for all $t\in\ct^b_c$ and all $s\in\cs$. 

Fix $t\in\ct^b_c$ and $s\in\cs$. By (ii) we can choose an integer $A>0$ with
$\Hom\big(\cl(\ct^b_c\cap\ct^{\leq-A})\,,\,s\big)=0$. Because $\car(s)$ belongs
to $\ct^b_c\subset\ct^+$ we may choose an integer
$A'>0$ so that $\Hom\big(\ct^{\leq-A'},\car(s)\big)=0$. Now take $m\geq1+\max(A,A')$, and choose a triangle $\Tm d\la e\la t\la d$ with $e\in\ct^c$ and
$d\in\ct^{\leq-m}$. Because $t\in\ct^b_c$ and $e\in\ct^c\subset\ct^b_c$
we have that $d\in\ct^b_c\cap\ct^{\leq-m}$.
Consider the commutative diagram with exact rows
\[\xymatrix@C+5pt{
  \Hom\big(\cl(d),s\big)\ar[r]\ar[d]& \Hom\big(\cl(t),s\big)\ar[r]^-a\ar[d]^b&
  \Hom\big(\cl(e),s\big)\ar[r]\ar[d]^c&\Hom\big(\cl(\Tm d),s\big)\ar[d] \\
  \Hom\big(d,\car(s)\big)\ar[r] &\Hom\big(t,\car(s)\big)\ar[r]^{a'} &
  \Hom\big(e,\car(s)\big)\ar[r] &\Hom\big(\Tm d,\car(s)\big)
}\]
By our choice of $m$ we know that
\[
\Hom\big(\cl(d),s\big)=0=\big(\cl(\Tm d),s\big),\quad
\Hom\big(d,\car(s)\big)=0=\Hom\big(\Tm d,\car(s)\big).
\]
Hence $a,a'$ are isomorphisms. But $c$ is an isomorphism by the compactness of
$e$, and therefore $b$ is an isomorphism.
\eprf

\appendix

\section{A criterion for checking that a triangulated functor is an equivalence}
\label{SA}

\lem{LA.1}
Let $\cl:\cu\la\cs$ be a triangulated functor with right
adjoint $\car:\cs\la\cu$. Suppose $P$ is a class of objects
in $\cu$ satisfying $P=\T P$, and such that  
\be
\item
  $P^\perp=\{0\}$, which means
  that if  $u\in\cu$ is an object and $\Hom(P,u)=0$ then $u=0$.
\item
$\cl(P)^\perp=\{0\}$, meaning
that if
$s\in\cs$ is an object and $\Hom\big(\cl(P),s\big)=0$ then $s=0$.
\item
The map $\Hom(p,u)\la\Hom\big(\cl(p),\cl(u)\big)$
is an isomorphism for objects $p\in P$ and $u\in\cu$.
\ee
Then $\cl$ and $\car$ are quasi-inverses.
\elem

\prf
Let $\eta:\id\la\car\cl$ and $\e:\cl\car\la\id$ be (respectively)
the unit and counit
of adjunction---it suffices to prove that $\eta$ and $\e$ are isomorphisms.

Let us begin with $\eta$. For objects $p\in P$ and $u\in\cu$ the natural
maps
\[\xymatrix@C+20pt{
\Hom_\cu^{}(p,u)\ar[r]^-\alpha & \Hom_\cs^{}\big(\cl(p),\cl(u)\big)
\ar[r]^-\beta & \Hom_\cu^{}\big(p,\car\cl(u)\big)
}\]
are both isomorphisms, $\alpha$ by (iii) and $\beta$ by the adjunction.
Hence the composite, which is the map
$\Hom(p,\eta):\Hom_\cu^{}(p,u)\la\Hom_\cu^{}\big(p,\car\cl(u)\big)$, must be
an isomorphism. Thus $\Hom(p,-)$ annihilates the mapping
cone of $\eta:u\la\car\cl(u)$, and by (iii) $\eta$ must be an
isomorphism.

We have proved that $\eta$ is an isomorphism, and the fact that the composite
$\car\stackrel{\eta\car}\la\car\cl\car\stackrel{\car\e}\la\car$
is the identity tells us that $\car\e:\car\cl\car\la\car$
must also be an isomorphism. Hence for any $p\in P$ and any object $s\in\cs$
we have that $\Hom_\cu^{}(p,-)$ takes that map
$\car\e(s):\car\cl\car(s)\la\car(s)$ to an isomorphism,
and adjunction tells us that $\Hom_\cs^{}\big(\cl(p),-\big)$ must
take the map $\e(s):\cl\car(s)\la s$ to an isomorphism. Applying (ii) to the
mapping cone of $\e(s)$, for every
$s\in\cs$, we deduce that $\e$ must be an isomorphism.
\eprf

\exm{EA.3}
Let $X$ be a scheme proper over the field $\C$
of complex numbers, and let $X^{\text{\rm an}}$ be the 
analytification of $X$. The category $\ct=\Dqc(X)$ is approximable and
$\C$--linear.
Now let $\cl:\ct^b_c\la\cs$ be the analytification functor
$\cl:\dcoh(X)\la\dcoh(X^{\text{\rm an}})$.
Then the hypotheses of
Theorem~\ref{T73.1} are satisfied, hence the functor $\cl$ has
a right adjoint $\car$.

Next we apply Lemma~\ref{LA.1}: for every closed point
$x\in X$ choose
a nonzero perfect complex $p(x)$ supported at $x$, and we set
\[
P\eq\Big\{\T^np(x)\,\,\left|\,\,\,n\in\zz,\text{ while } x\in X\text{ is a closed point}\Big\}\ .\right.
\]
It's an easy exercise
to show that this choice of $P$ satisfies the hypotheses of Lemma~\ref{LA.1},
hence
the functor $\cl$ must be an equivalence.

The idea at the core of the argument above is due to Jack
Hall, the technical problem he faced was that the representability
theorems available to him were less powerful than Theorem~\ref{T73.1}.
The reader is also referred to Serre~\cite{Serre56} for the
first proof of the version of GAGA in the couple of
paragraphs above,
and to Hall~\cite{Hall18,Hall22} for generalizations
of his core idea that go in a direction different from the one
of this article.
\eexm

\rmk{RA.5}
The reader should note that the proof
in Example~\ref{EA.3} depends on substantial
structural theorems about $\Dqc(X)$, but all we need to know
about the category $\cs=\dcoh(X^{\text{\rm an}})$
and the functor $\cl:\dcoh(X)\la\dcoh(X^{\text{\rm an}})$
is the minimal data that goes into showing that the hypotheses of
Theorem~\ref{T73.1} and Lemma~\ref{LA.1} are satisfied.
This is the reason that all the GAGA theorems in algebraic geometry
are special cases of general theorems as in  Hall~\cite[Theorems~A and B]{Hall22};
in fact the $P\subset\dcoh(X)$ of the proof of Example~\ref{EA.3}
works for all of them.
\ermk

\bibliographystyle{amsplain}
\bibliography{stan}

\providecommand{\bysame}{\leavevmode\hbox to3em{\hrulefill}\thinspace}
\providecommand{\MR}{\relax\ifhmode\unskip\space\fi MR }
\providecommand{\MRhref}[2]{%
  \href{http://www.ams.org/mathscinet-getitem?mr=#1}{#2}
}
\providecommand{\href}[2]{#2}
\begin{thebibliography}{10}

\bibitem{Adams71}
J.~Frank Adams, \emph{A variant of {E. H. Brown's} representability theorem},
  Topology \textbf{10} (1971), 185--198.

\bibitem{Alonso-Jeremias-Souto03}
Leovigildo Alonso~Tarr{\'{\i}}o, Ana Jerem{\'{\i}}as~L{\'o}pez, and
  Mar{\'{\i}}a~Jos{\'e} Souto~Salorio, \emph{Construction of {$t$}-structures
  and equivalences of derived categories}, Trans. Amer. Math. Soc. \textbf{355}
  (2003), no.~6, 2523--2543 (electronic).

\bibitem{Aoki20}
Ko~Aoki, \emph{Quasiexcellence implies strong generation}, J. Reine Angew.
  Math. \textbf{780} (2021), 133--138.

\bibitem{BeiBerDel82}
Alexander~A. Be{\u\i}linson, Joseph Bernstein, and Pierre Deligne,
  \emph{Analyse et topologie sur les {\'e}spaces singuliers}, Ast{\'e}risque,
  vol. 100, Soc. Math. France, 1982 (French).

\bibitem{BenZvi-Nadler-Preygel16}
David Ben-Zvi, David Nadler, and Anatoly Preygel, \emph{Integral transforms for
  coherent sheaves}, J. Eur. Math. Soc. (JEMS) \textbf{19} (2017), no.~12,
  3763--3812.

\bibitem{BondalvandenBergh04}
Alexei~I. Bondal and Michel Van~den Bergh, \emph{Generators and
  representability of functors in commutative and noncommutative geometry},
  Mosc. Math. J. \textbf{3} (2003), no.~1, 1--36, 258.

\bibitem{Bondarko22}
Mikhail~V. Bondarko, \emph{Producing ``new'' semi-orthogonal decompositions in
  arithmetic geometry}, arXiv:2203.07315.

\bibitem{Bondarko-Vostokov20}
Mikhail~V. Bondarko and S.~V. Vostokov, \emph{On weakly negative subcategories,
  weight structures, and (weakly) approximable triangulated categories},
  Lobachevskii J. Math. \textbf{41} (2020), 151–--159.

\bibitem{Burke-Neeman-Pauwels18}
Jesse Burke, Amnon Neeman, and Bregje Pauwels, \emph{Gluing approximable
  triangulated categories}, Forum Math. Sigma \textbf{11} (2023), Paper No.
  e110, 18.

\bibitem{Christensen96}
J.~Daniel Christensen, \emph{Ideals in triangulated categories: Phantoms,
  ghosts and skeleta}, Advances in Mathematics \textbf{136} (1998), 284--339.

\bibitem{Christensen-Keller-Neeman99}
J.~Daniel Christensen, Bernhard Keller, and Amnon Neeman, \emph{Failure of
  {Brown} representability in derived categories}, Topology \textbf{40} (2001),
  1339--1361.

\bibitem{Hall18}
Jack Hall, \emph{G{AGA} theorems}, J. Math. Pures Appl. (9) \textbf{175}
  (2023), 109--142. \MR{4598930}

\bibitem{Hall22}
\bysame, \emph{{GAGA} theorems}, J. Math. Pures Appl. \textbf{175} (2023),
  109–142.

\bibitem{Hall-Rydh13}
Jack Hall and David Rydh, \emph{Perfect complexes on algebraic stacks}, Compos.
  Math. \textbf{153} (2017), no.~11, 2318--2367.

\bibitem{Illusie71B}
Luc Illusie, \emph{Existence de r{\'e}solutions globales}, Th\'eorie des
  intersections et th\'eor\`eme de {R}iemann-{R}och, Springer-Verlag, Berlin,
  1971, S\'eminaire de G\'eom\'etrie Alg\'ebrique du Bois-Marie 1966--1967 (SGA
  6, Expos{\'e} II), pp.~160--221. Lecture Notes in Mathematics, Vol. 225.

\bibitem{Illusie71A}
\bysame, \emph{G{\'e}n{\'e}ralit{\'e}s sur les conditions de finitude dans les
  cat{\'e}gories d{\'e}riv{\'e}es}, Th\'eorie des intersections et th\'eor\`eme
  de {R}iemann-{R}och, Springer-Verlag, Berlin, 1971, S\'eminaire de
  G\'eom\'etrie Alg\'ebrique du Bois-Marie 1966--1967 (SGA 6, Expos{\'e} I),
  pp.~78--159. Lecture Notes in Mathematics, Vol. 225.

\bibitem{Karmazyn-Kuznetsov-Shinder22}
Joseph Karmazyn, Alexander Kuznetsov, and Evgeny Shinder, \emph{Derived
  categories of singular surfaces}, J. Eur. Math. Soc. (JEMS) \textbf{24}
  (2022), no.~2, 461--526.

\bibitem{Lipman-Neeman07}
Joseph Lipman and Amnon Neeman, \emph{Quasi-perfect scheme maps and boundedness
  of the twisted inverse image functor}, Illinois J. Math. \textbf{51} (2007),
  209--236.

\bibitem{Neeman96}
Amnon Neeman, \emph{The {Grothendieck} duality theorem via {Bousfield's}
  techniques and {Brown} representability}, Jour. Amer. Math. Soc. \textbf{9}
  (1996), 205--236.

\bibitem{Neeman4}
\bysame, \emph{On a theorem of {Brown} and {Adams}}, Topology \textbf{36}
  (1997), 619--645.

\bibitem{Neeman17}
\bysame, \emph{Strong generators in {$\dperf X$} and {$\dcoh(X)$}}, Ann. of
  Math. (2) \textbf{193} (2021), no.~3, 689--732.

\bibitem{Neeman22A}
\bysame, \emph{Bounded {$t$}-structures on the category of perfect complexes},
  Acta Math. \textbf{233} (2024), no.~2, 239--284.

\bibitem{Rouquier08}
Rapha{\"e}l Rouquier, \emph{Dimensions of triangulated categories}, J. K-Theory
  \textbf{1} (2008), no.~2, 193--256.

\bibitem{Saorin-Stovicek23}
Manuel Saor\'{\i}n and Jan \v{S}t'ov\'{\i}\v{c}ek, \emph{{$t$}-structures with
  {G}rothendieck hearts via functor categories}, Selecta Math. (N.S.)
  \textbf{29} (2023), no.~5, Paper No. 77, 73.

\bibitem{Serre56}
Jean-Pierre Serre, \emph{G\'eom\'etrie alg\'ebrique et g\'eom\'etrie
  analytique}, Ann. Inst. Fourier, Grenoble \textbf{6} (1955--1956), 1--42.

\bibitem{stacks-project}
The {Stacks Project Authors}, \emph{\itshape stacks project},
  http://stacks.math.columbia.edu.

\bibitem{ThomTro}
Robert~W. Thomason and Thomas~F. Trobaugh, \emph{Higher algebraic {K--theory}
  of schemes and of derived categories}, The Grothendieck Festschrift ( a
  collection of papers to honor Grothendieck's 60'th birthday), vol.~3,
  Birkh{\"a}user, 1990, pp.~247--435.

\end{thebibliography}

\end{document}